\documentclass[12pt,a4paper]{amsart} 
\usepackage{amscd,amssymb,amsmath, pictex}

\newcounter{myenumi}
\newenvironment{Enumerate}
{\begin{list} {\textnormal{(\arabic{myenumi})}}
{\usecounter{myenumi}\setlength{\leftmargin}{1.5em}\setlength{\labelwidth}{2.5em}}}
{\end{list}}

\newenvironment{Itemize}
{\begin{list}{$\bullet$}{\setlength{\leftmargin}{2em}\setlength{\labelwidth}{1em}}}
{\end{list}}

\newtheorem{thm}{Theorem}[section]
\newtheorem{cor}[thm]{Corollary}
\newtheorem{lem}[thm]{Lemma}
\newtheorem{prop}[thm]{Proposition}

\theoremstyle{definition}
\newtheorem{dfn}[thm]{Definition}

\theoremstyle{remark}
\newtheorem{rmk}[thm]{Remark}
\newtheorem{rmks}[thm]{Remarks}
\newtheorem{example}[thm]{Example}
\newtheorem{examples}[thm]{Examples}

\newtheorem{ntn}[thm]{Notation}

\def\id{\operatorname{id}}
\def\coker{\operatorname{coker}}
\def\ker{\operatorname{ker}}
\def\Im{\operatorname{Im}}

\numberwithin{equation}{section}
\newcommand{\F}{\mathcal{F}}
\newcommand{\e}{\epsilon}
\def\cs#1{{
\ensuremath{\mathbin{\stackrel{#1}{\leftharpoondown}}}}}

\newcommand{\cocycle}{\mathfrak{s}}

\newcommand{\field}[1]{\mathbf{#1}}
\newcommand{\CC}{\field{C}}
\newcommand{\NN}{\field{N}}
\newcommand{\QQ}{\field{Q}}
\newcommand{\ZZ}{\field{Z}}
\newcommand{\TT}{\field{T}}

\newcommand{\LAMBDA}{\mathbf{\Lambda}}

\newcommand{\Oo}{\mathcal{O}}

\newcommand{\Mm}{\mathcal{M}}

\newcommand{\lsp}{\operatorname{span}}
\newcommand{\clsp}{\overline{\lsp}}

\newcommand{\tgrphlim}{{
    \renewcommand{\leftarrow}{\leftharpoondown}
    \varprojlim
}}

\title{$C^*$-algebras associated to coverings of $k$-graphs}
\author{Alex Kumjian}
\address{Alex Kumjian\\ Department of Mathematics (084)\\ University
of Nevada\\ Reno NV 89557-0084\\ USA} \email{alex@unr.edu}
\author{David Pask}
\address{David Pask and Aidan Sims\\ School of Mathematics and
Applied Statistics  \\
University of Wollongong\\
NSW  2522\\
AUSTRALIA} \email{dpask, asims@uow.edu.au}
\author{Aidan Sims}

\thanks{This research was supported by the Australian Research Council.}

\begin{document}

\begin{abstract}
A covering of $k$-graphs (in the sense of Pask-Quigg-Raeburn) induces an
embedding of universal $C^*$-algebras. We show how to build a
$(k+1)$-graph whose universal algebra encodes this embedding. More
generally we show how to realise a direct limit of $k$-graph algebras
under embeddings induced from coverings as the universal algebra of a
$(k+1)$-graph. Our main focus is on computing the $K$-theory of the
$(k+1)$-graph algebra from that of the component $k$-graph algebras.

Examples of our construction include a realisation of the Kirchberg
algebra $\mathcal{P}_n$ whose $K$-theory is opposite to that of
$\mathcal{O}_n$, and a class of A$\TT$-algebras that can naturally be
regarded as higher-rank Bunce-Deddens algebras.
\end{abstract}

\maketitle

\section{Introduction}

A directed graph $E$ consists of a countable collection $E^0$ of
vertices, a countable collection $E^1$ of edges, and maps $r, s : E^1
\to E^0$ which give the edges their direction; the edge $e$ points
from $s(e)$ to $r(e)$. Following the convention established in
\cite{r}, the associated graph algebra $C^*(E)$ is the universal
$C^*$-algebra generated by partial isometries $\{s_e : e \in E^1\}$
together with mutually orthogonal projections $\{p_v : v \in E^0\}$
such that $p_{s(e)} = s^*_e s_e$ for all $e \in E^1$, and $p_v \ge
\sum_{e \in F} s_e s^*_e$ for all $v \in E^0$ and finite $F \subset
r^{-1}(v)$, with equality when $F = r^{-1}(v)$ is finite and
nonempty.

Graph algebras, introduced in \cite{EW, KPRR}, have been studied
intensively in recent years because much of the structure of $C^*(E)$
can be deduced from elementary features of $E$. In particular, graph
$C^*$-algebras are an excellent class of models for Kirchberg
algebras, because it is easy to tell from the graph $E$ whether
$C^*(E)$ is simple and purely infinite \cite{KPR}. Indeed, a
Kirchberg algebra can be realised up to Morita equivalence as a graph
$C^*$-algebra if and only if its $K_1$-group is torsion-free
\cite{Sz}.  It is also true that every AF algebra can be realised up
to Morita equivalence as a graph algebra; the desired graph is a
Bratteli diagram for the AF algebra in question (see \cite{D} or
\cite{JT}). However, this is the full extent to which graph algebras
model simple classifiable $C^*$-algebras due to the following
dichotomy: if $E$ is a directed graph and $C^*(E)$ is simple, then
$C^*(E)$ is either AF or purely infinite (see
\cite[Corollary~3.10]{KPR}, \cite[Remark~5.6]{BPRS}).

Higher-rank graphs, or $k$-graphs, and their $C^*$-algebras
were originally developed by the first two authors \cite{kp} to
provide a graphical framework for the higher-rank Cuntz-Krieger
algebras of Robertson and Steger \cite{rst}. A $k$-graph
$\Lambda$ is a kind of $k$-dimensional graph, which one can
think of as consisting of vertices $\Lambda^0$ together with
$k$ collections of edges $\Lambda^{e_1}, \dots, \Lambda^{e_k}$
which we think of as lying in $k$ different dimensions. As an
aid to visualisation, we often distinguish the different types
of edges using $k$ different colours.

Higher-rank graphs and their $C^*$-algebras are generalisations of
directed graphs and their algebras. Given a directed graph $E$, its path
category $E^*$ is a $1$-graph, and the $1$-graph $C^*$-algebra $C^*(E^*)$
as defined in \cite{kp} is canonically isomorphic to the graph algebra
$C^*(E)$ as defined in \cite{KPRR}. Furthermore, every $1$-graph arises
this way, so the class of graph algebras and the class of $1$-graph
algebras are one and the same. For $k \ge 2$, there are many $k$-graph
algebras which do not arise as graph algebras. For example, the original
work of Robertson and Steger on higher-rank Cuntz-Krieger algebras
describes numerous $2$-graphs $\Lambda$ for which $C^*(\Lambda)$ is a
Kirchberg algebra and $K_1(C^*(\Lambda))$ contains torsion.

Recent work of Pask, Raeburn, R{\o}rdam and Sims has shown that one can
also realise a substantial class of A$\TT$-algebras as $2$-graph algebras,
and that one can tell from the $2$-graph whether or not the resulting
$C^*$-algebra is simple and has real-rank zero \cite{PRRS}. The basic idea
of the construction in \cite{PRRS} is as follows. One takes a Bratteli
diagram in which the edges are coloured red, and replaces each vertex with
a blue simple cycle (there are technical restrictions on the relationship
between the lengths of the blue cycles and the distribution of the red
edges joining them, but this is the gist of the construction). The
resulting $2$-graph is called a \emph{rank-2 Bratteli diagram}. The
associated $C^*$-algebra is A$\TT$ because the $C^*$-algebra of a simple
cycle of length $n$ is isomorphic to $M_n(C(\TT))$ \cite{aHR}. The results
of \cite{PRRS} show how to read off from a rank-2 Bratteli diagram the
$K$-theory, simplicity or otherwise, and real-rank of the resulting A$\TT$
algebra.

The construction explored in the current paper is motivated by the
following example of a rank-2 Bratteli diagram. For each $n \in \NN$, let
$L_{2^n}$ be the simple directed loop graph with $2^n$ vertices labelled
$0, \dots, 2^n -1$ and $2^n$ edges $f_0, \dots, f_{2^n - 1}$, where $f_i$
is directed from the vertex labelled $i+1\ ({\rm mod}\ 2^n)$ to the vertex
labelled $i$. We specify a rank-2 Bratteli diagram $\Lambda(2^\infty)$ as
follows. The $n^{\rm th}$ level of $\Lambda(2^\infty)$ consists of a
single blue copy of $L_{2^{n-1}}$ $(n = 1, 2, \cdots)$. For $0 \le i \le
2^n - 1$, there is a single red edge from the vertex labelled $i$ at the
$(n+1)^{\rm st}$ level to the vertex labelled $i\ ({\rm mod}\ 2^n)$ at the
$n^{\rm th}$ level. The $C^*$-algebra of the resulting $2$-graph is Morita
equivalent to the Bunce-Deddens algebra of type $2^\infty$, and this was
one of the first examples of a $2$-graph algebra which is simple but
neither purely infinite nor AF (see \cite[Example~6.7]{PRRS}).

\medskip

The purpose of this paper is to explore the observation that the growing
blue cycles in $\Lambda(2^\infty)$ can be thought of as a tower of
\emph{coverings} of $1$-graphs (roughly speaking, a covering is a locally
bijective surjection --- see Definition~\ref{dfn:covering}), where the red
edges connecting levels indicate the covering maps.

\setcounter{footnote}{1} In Section~2, we describe how to
construct $(k+1)$-graphs from coverings. In its simplest form,
our construction takes $k$-graphs $\Lambda$ and $\Gamma$ and a
covering map $p : \Gamma \to \Lambda$, and produces a $(k +
1)$-graph $\Lambda \cs{p} \Gamma$ in which each edge in the
$(k+1)^{\rm st}$ dimension points from a vertex $v$ of $\Gamma$
to the vertex $p(v)$ of $\Lambda$ which it covers\footnote{In
its full generality, our construction is more complicated (see
Proposition~\ref{prp:multiple covering}), enabling us to
recover the important example of the irrational rotation
algebras discussed in \cite{PRRS}. To keep technical detail in
this introduction to a minimum, we discuss only the basic
construction here.}. Building on this construction, we show how
to take an infinite tower of coverings $p_n : \Lambda_{n+1} \to
\Lambda_n$, $n \ge 1$, and produce an infinite $(k+1)$-graph
$\tgrphlim(\Lambda_n, p_n)$ with a natural inductive structure
(Corollary~\ref{cor:tower graph}).

The next step, achieved in Section~\ref{sec:C*algs}, is to determine how
the universal $C^*$-algebra of $\Lambda \cs{p} \Gamma$ relates to those of
$\Lambda$ and $\Gamma$. We show that $C^*(\Lambda \cs{p} \Gamma)$ is
Morita equivalent to $C^*(\Gamma)$ and contains an isomorphic copy of
$C^*(\Lambda)$ (Proposition~\ref{prp:embedding}). We then show that given
a system of coverings $p_n : \Lambda_{n+1} \to \Lambda_n$, the
$C^*$-algebra $C^*(\tgrphlim(\Lambda_n, p_n))$ is Morita equivalent to a
direct limit of the $C^*(\Lambda_n)$ (Theorem~\ref{thm:direct limit}).

In Section~\ref{sec:simplicity}, we use results of \cite{RobSi}
to characterise simplicity of $C^*(\tgrphlim(\Lambda_n, p_n))$,
and we also give a sufficient condition for this $C^*$-algebra
to be purely infinite. In Section~\ref{sec:K-th}, we show how
various existing methods of computing the $K$-theory of the
$C^*(\Lambda_n)$ can be used to compute the $K$-theory of
$C^*(\tgrphlim(\Lambda_n, p_n))$. Our results boil down to
checking that each of the existing $K$-theory computations for
the $C^*(\Lambda_n)$ is natural in the appropriate sense. Given
that $K$-theory for higher-rank graph $C^*$-algebras has proven
quite difficult to compute in general (see \cite{e}), our
$K$-theory computations are an important outcome of the paper.

We conclude in Section~\ref{sec:examples} by exploring some
detailed examples which illustrate the covering-system
construction, and show how to apply our $K$-theory calculations
to the resulting higher-rank graph $C^*$-algebras. For integers
$3 \le n < \infty$, we obtain a $3$-graph algebra realisation
of Kirchberg algebra $\mathcal{P}_n$ whose $K$-theory is
opposite to that of $\mathcal{O}_n$ (see
Section~\ref{subsec:P_n}). We also obtain, using $3$-graphs, a
class of simple A$\TT$-algebras with real-rank zero which
cannot be obtained from the rank-2 Bratteli diagram
construction of \cite{PRRS} (see Section~\ref{subsec:Delta_2}),
and which we can describe in a natural fashion as higher-rank
analogues of the Bunce-Deddens algebras. These are, to our
knowledge, the first explicit computations of $K$-theory for
infinite classes of $3$-graph algebras.

\section{Covering systems of $k$-graphs}\label{sec:construction}

For $k$-graphs we adopt the conventions of \cite{kp, pqr, rsy};
briefly, a \emph{$k$-graph} is a countable small category $\Lambda$
equipped with a functor $d :\Lambda \to \NN^k$ satisfying the
\emph{factorisation property}: for all $\lambda \in \Lambda$ and
$m,n \in \NN^k$ such that $d( \lambda )=m+n$ there exist unique $\mu
,\nu \in \Lambda$ such that $d(\mu)=m$, $d(\nu)=n$, and $\lambda=\mu
\nu$. When $d(\lambda )=n$ we say $\lambda$ has \emph{degree} $n$.
By abuse of notation, we will use $d$ to denote the degree functor
in every $k$-graph in this paper; the domain of $d$ is always clear
from context.

The standard generators of $\NN^k$ are denoted $e_1, \dots, e_k$,
and for $n \in \NN^k$ and $1 \le i \le k$ we write $n_i$ for the
$i^{\rm th}$ coordinate of $n$.

If $\Lambda$ is a $k$-graph, the \emph{vertices} are the morphisms of
degree 0. The factorisation property implies that these are precisely the
identity morphisms, and so can be identified with the objects. For $\alpha
\in\Lambda$, the \emph{source} $s(\alpha)$ is the domain of $\alpha$, and
the \emph{range} $r(\alpha)$ is the codomain of $\alpha$ (strictly
speaking, $s(\alpha)$ and $r(\alpha)$ are the identity morphisms
associated to the domain and codomain of $\alpha$).

For $n \in \NN^k$, we write $\Lambda^n$ for $d^{-1}(n)$. In particular,
$\Lambda^0$ is the vertex set. For $u,v\in\Lambda^0$ and $E \subset
\Lambda$, we write $u E := E \cap r^{-1}(u)$ and $E v := E \cap
s^{-1}(v)$. For $n \in \NN^k$, we write
\[
\Lambda^{\le n} := \{\lambda \in \Lambda : d(\lambda) \le n , s(\lambda)
\Lambda^{e_i}  = \emptyset \text{ whenever } d(\lambda) + e_i \le n\} .
\]
We say that $\Lambda$ is \emph{connected} if the equivalence relation on
$\Lambda^0$ generated by $\{(v,w) \in \Lambda^0 \times \Lambda^0 :
v\Lambda w \not= \emptyset\}$ is the whole of $\Lambda^0 \times
\Lambda^0$. A \emph{morphism} between $k$-graphs is a degree-preserving
functor.

We say that $\Lambda$ is \emph{row-finite} if $v\Lambda^n$ is finite for
all $v \in \Lambda^0$ and $n \in \NN^k$. We say that $\Lambda$ is
\emph{locally convex} if whenever $1 \le i < j \le k$, $e \in
\Lambda^{e_i}$, $f \in \Lambda^{e_j}$ and $r(e) = r(f)$, we can extend
both $e$ and $f$ to paths $ee'$ and $ff'$ in $\Lambda^{e_i + e_j}$.

We next introduce the notion of a covering of one $k$-graph by another.
For a more detailed treatment of coverings of $k$-graphs, see \cite{pqr}.

\begin{dfn}\label{dfn:covering}
A \emph{covering} of a $k$-graph $\Lambda$ is a surjective $k$-graph
morphism $p :\Gamma \to \Lambda$ such that for all $v\in \Gamma^0$, $p$
maps $\Gamma v$ 1-1 onto $\Lambda p(v)$ and $v\Gamma$ 1-1 onto
$p(v)\Lambda$.  A covering $p :\Gamma \to\Lambda$ is \emph{connected} if
$\Gamma$, and hence also $\Lambda$, is connected.

A covering $p : \Gamma \to \Lambda$ is \emph{finite} if $p^{-1} (v)$ is
finite for all $v \in \Lambda^0$.
\end{dfn}

\begin{rmks}
\begin{Enumerate}
\item A covering $p : \Gamma \to \Lambda$ has the unique path
lifting property: for every $\lambda \in \Lambda$ and $v \in \Gamma^0$
with $p(v)=s(\lambda)$ there exists a unique $\gamma$ such that $p (
\gamma ) = \lambda$ and $s ( \gamma ) = v$; likewise, if $p(v) = r (
\lambda )$ there is a unique $\zeta$ such that $p ( \zeta ) = \lambda$ and
$r(\zeta) = v$.

\item If $\Lambda$ is connected then surjectivity of $p$ is
implied by the unique path-lifting property.

\item If there is a fixed integer $n$ such that $|p^{-1}(v)| = n$ for all $v
\in \Lambda^0$, $p$ is said to be an \emph{$n$-fold covering}. If
$\Gamma$ is connected, then $p$ is automatically an $n$-fold covering
for some $n$.
\end{Enumerate}
\end{rmks}

\begin{ntn}\label{ntn:S_m}
For $m \in \NN\setminus\{0\}$, we write $S_m$ for the group of
permutations of the set $\{1, \dots, m\}$. We denote both composition of
permutations in $S_m$, and the action of a permutation in $S_m$ on an
element of $\{1, \dots, m\}$ by juxtaposition; so for $\phi,\psi \in S_m$,
$\phi\psi \in S_m$ is the permutation $\phi \circ \psi$, and for $\phi \in
S_m$ and $j \in \{1, \dots, m\}$, $\phi j \in \{1, \dots, m\}$ is the
image of $j$ under $\phi$. When convenient, we regard $S_m$ as (the
morphisms of) a category with a single object.
\end{ntn}

\begin{dfn}
Fix $k,m \in \NN\setminus\{0\}$, and let $\Lambda$ be a $k$-graph. A
\emph{cocycle} $\cocycle : \Lambda \to S_m$ is a functor $\lambda \mapsto
\cocycle(\lambda)$ from the category $\Lambda$ to the category $S_m$. That
is, whenever $\alpha,\beta \in \Lambda$ satisfy $s(\alpha) = r(\beta)$ we
have $\cocycle(\alpha)\cocycle(\beta) = \cocycle(\alpha\beta)$.
\end{dfn}

We are now ready to describe the data needed for our construction.

\begin{dfn}\label{dfn:covering system}
A \emph{covering system of $k$-graphs} is a quintuple
$(\Lambda, \Gamma, p, m, \cocycle)$ where $\Lambda$ and
$\Gamma$ are $k$-graphs, $p : \Lambda \to \Gamma$ is a
covering, $m$ is a nonzero positive integer, and $\cocycle :
\Gamma \to S_m$ is a cocycle. We say that the covering system
is \emph{row finite} if the covering map $p$ is finite and both
$\Lambda$ and $\Gamma$ are row finite. When $m = 1$ and
$\cocycle$ is the identity cocycle, we drop references to $m$
and $\cocycle$ altogether, and say that $(\Lambda, \Gamma, p)$
is a covering system of $k$-graphs.
\end{dfn}

Given a covering system $(\Lambda, \Gamma, p, m, \cocycle)$ of
$k$-graphs, we will define a $(k+1)$-graph $\Lambda \cs{p,
\cocycle} \Gamma$ which encodes the covering map. Before the
formal statement of this construction, we give an intuitive
description of $\Lambda \cs{p,\cocycle} \Gamma$. The idea is
that $\Lambda \cs{p,\cocycle} \Gamma$ is a $(k+1)$-graph
containing disjoint copies $\imath(\Lambda)$ and
$\jmath(\Gamma)$ of the $k$-graphs $\Lambda$ and $\Gamma$ in
the first $k$ dimensions. The image $\jmath(v)$ of a vertex $v
\in \Gamma$ is connected to the image $\imath(p(v))$ of the
vertex it covers in $\Lambda$ by $m$ parallel edges $e(v,1),
\dots, e(v,m)$ of degree $e_{k+1}$. Factorisations of paths in
$\Lambda \cs{p,\cocycle} \Gamma$ involving edges $e(v,l)$ of
degree $e_{k+1}$ are determined by the unique path-lifting
property and the cocycle $\cocycle$.

It may be helpful on the first reading to consider the case
where $m = 1$ so that $\cocycle$ is necessarily trivial. To
state the result formally, we first establish some notation.

\begin{ntn}\label{ntn:N^k embeddings}
Fix $k > 0$. For $n \in \NN^k$ we denote by $(n, 0_1)$ the
element $\sum^k_{i=1} n_i e_i \in \NN^{k+1}$ and for $m \in
\NN$, we denote by $(0_k, m)$ the element $m e_{k+1} \in
\NN^{k+1}$. We write $(\NN^k, 0_1)$ for $\{(n,0_1) : n \in
\NN^k\}$ and $(0_k, \NN)$ for $\{(0_k,m) : m \in \NN\}$.

Given a $(k+1)$-graph $\Xi$, we write $\Xi^{(0_k,\NN)}$ for
$\{\xi \in \Xi : d(\xi) \in (0_k, \NN)\}$, and we write
$\Xi^{(\NN^k, 0_1)}$ for $\{\xi \in \Xi : d(\xi) \in (\NN^k,
0_1)\}$. When convenient, we regard $\Xi^{(0_k,\NN)}$ as a
$1$-graph and $\Xi^{(\NN^k,0_1)}$ as a $k$-graph, ignoring the
distinctions between $\NN$ and $(0_k, \NN)$ and between $\NN^k$
and $(\NN^k, 0_1)$.
\end{ntn}

\begin{prop}\label{prp:Lambda cs Gamma}
Let $(\Lambda,\Gamma,p,m,\cocycle)$ be a covering system of
$k$-graphs. There is a unique $(k+1)$-graph $\Lambda \cs{p,\cocycle}
\Gamma$ such that:
\begin{Enumerate}
\item
there are injective functors $\imath : \Lambda \to \Lambda \cs{p,\cocycle}
\Gamma$ and $\jmath : \Gamma \to \Lambda \cs{p,\cocycle} \Gamma$ such that
$d(\imath(\alpha)) = (d(\alpha), 0_1)$ and $d(\jmath(\beta)) = (d(\beta),
0_1)$ for all $\alpha \in \Lambda$ and $\beta \in \Gamma$;
\item
$\imath(\Lambda) \cap \jmath(\Gamma) = \emptyset$ and $\imath(\Lambda)
\cup \jmath(\Gamma) = \{\tau \in \Lambda \cs{p,\cocycle} \Gamma :
d(\tau)_{k+1} = 0\}$;
\item
there is a bijection $e : \Gamma^0 \times \{1, \dots, m\} \to (\Lambda
\cs{p,\cocycle} \Gamma)^{e_{k+1}}$;
\item
$s(e(v,l)) = \jmath(v)$ and $r(e(v,l)) = \imath(p(v))$ for all $v \in
\Gamma^0$ and $1 \le l \le m$; and
\item
$e(r(\lambda),l) \jmath(\lambda) = \imath(p(\lambda)) e(s(\lambda),
\cocycle(\lambda)^{-1} l)$ for all $\lambda \in \Gamma$ and $1 \le l \le
m$.
\end{Enumerate}
If the covering system $(\Lambda,\Gamma,p,m,\cocycle)$ is row finite,
then $\Lambda \cs{p, \cocycle} \Gamma$ is row finite. Moreover,
$\Lambda$ is locally convex if and only if $\Gamma$ is locally
convex, and in this case $\Lambda \cs{p, \cocycle} \Gamma$ is also
locally convex.
\end{prop}

\begin{ntn} \label{ntn:m=1}
If $m = 1$ so that $\cocycle$ is necessarily trivial, we drop
all reference to $\cocycle$. We denote $\Lambda \cs{p,\cocycle}
\Gamma$ by $\Lambda \cs{p} \Gamma$, and write $(\Lambda \cs{p}
\Gamma)^{e_{k + 1}} = \{e(v) : v \in \Gamma^0\}$. In this case,
the factorisation property is determined by the unique
path-lifting property alone.
\end{ntn}

The main ingredient in the proof of Proposition~\ref{prp:Lambda cs Gamma}
is the following fact from \cite[Remark~2.3]{fs} (see also
\cite[Section~2]{rsy}).

\begin{lem}\label{lem:FS lemma}
Let $E_1, \dots, E_k$ be $1$-graphs with the same vertex set $E^0$. For
distinct $i,j \in \{1, \dots, k\}$, let $E_{i,j} := \{(e,f) \in E_i^1
\times E_j^1 : s(e) = r(f)\}$, and write $r\big((e,f)\big) = r(e)$ and
$s\big((e,f)\big) = s(f)$. For distinct $h,i,j \in \{1, \dots, k\}$, let
$E_{h,i,j} := \{(e,f,g) \in E^1_h \times E^1_i \times E^1_j : (e,f) \in
E_{h,i}, (f,g) \in E_{i,j}\}$.

Suppose we have bijections $\theta_{i,j} : E_{i,j} \to E_{j,i}$
such that $r \circ \theta_{i,j} = r$, $s \circ \theta_{i,j} =
s$ and $\theta_{i,j} \circ \theta_{j,i} = \id$, and such that
\begin{equation} \label{associativity}
(\theta_{i,j} \times \id)(\id \times \theta_{h,j})(\theta_{h,i} \times
\id) = (\id \times \theta_{h,i})(\theta_{h,j} \times \id)(\id \times
\theta_{i,j})
\end{equation}
as bijections from $E_{h,i,j}$ to $E_{j,i,h}$.

Then there is a unique $k$-graph $\Lambda$ such that $\Lambda^0 = E^0$,
$\Lambda^{e_i} = E_i^1$ for $1 \le i \le k$, and for distinct $i,j \in
\{1, \dots, k\}$ and $(e,f) \in E_{i,j}$, the pair $(f',e') \in E_{j,i}$
such that $(f',e') = \theta_{i,j}(e,f)$ satisfies $ef = f'e'$ as morphisms
in $\Lambda$.
\end{lem}

\begin{rmk} \label{keyfact}
Every $k$-graph arises in this way: Given a $k$-graph $\Lambda$, let $E^0
:= \Lambda^0$, and $E^1_i := \Lambda^{e_i}$ for $1 \le i \le k$, and
define $r,s : E^1_i \to E^0$ by restriction of the range and source maps
in $\Lambda$. Define bijections $\theta_{i,j} : E_{i,j} \to E_{j,i}$ via
the factorisation property: $\theta_{i,j}(e,f)$ is equal to the unique
pair $(f',e') \in E_{j,i}$ such that $ef = f'e'$ in $\Lambda$. Then
condition~(\ref{associativity}) holds by the associativity of the category
$\Lambda$, and the uniqueness assertion of Lemma~\ref{lem:FS lemma}
implies that $\Lambda$ is isomorphic to the $k$-graph obtained from the
$E_i$ and the $\theta_{i,j}$ using Lemma~\ref{lem:FS lemma}.
\end{rmk}

Lemma~\ref{lem:FS lemma} tells us how to describe a $k$-graph
pictorially. As in \cite{rsy, PRRS}, the \emph{skeleton} of a
$k$-graph $\Lambda$ is the directed graph $E_\Lambda$ with vertices
$E_\Lambda^0 = \Lambda^0$, edges $E_\Lambda^1 = \bigcup^k_{i=1}
\Lambda^{e_i}$, range and source maps inherited from $\Lambda$, and
edges of different degrees in $\Lambda$ distinguished using $k$
different colours in $E_\Lambda$: in this paper, we will often refer
to edges of degree $e_1$ as ``blue" and edges of degree $e_2$ as
``red." Lemma~\ref{lem:FS lemma} implies that the skeleton
$E_\Lambda$ together with the factorisation rules $f g = g' f'$ where
$f,f' \in \Lambda^{e_i}$ and $g,g' \in \Lambda^{e_j}$ completely
specify $\Lambda$. In practise, we draw $E_\Lambda$ using solid,
dashed and dotted edges to distinguish the different colours, and
list the factorisation rules separately.

\begin{proof}[Proof of Proposition~\ref{prp:Lambda cs Gamma}]
The idea is to apply Lemma~\ref{lem:FS lemma} to obtain the $(k+1)$-graph
$\Lambda \cs{p,\cocycle} \Gamma$. We first define sets $E^0$ and $E^1_i$
for $1 \le i \le k+1$. As a set, $E^0$ is a copy of the disjoint union
$\Lambda^0 \sqcup \Gamma^0$. We denote the copy of $\Lambda^0$ in $E^0$ by
$\{\imath(v) : v \in \Lambda^0\}$ and the copy of $\Gamma^0$ in $E^0$ by
$\{\jmath(w) : w \in \Gamma^0\}$ where as yet the $\imath(v)$ and
$\jmath(w)$ are purely formal symbols. So
\[
E^0 = \{\imath(v) : v \in \Lambda^0\} \sqcup \{\jmath(w) : w \in
\Gamma^0\}.
\]
For $1 \le i \le k$, we define, in a similar fashion,
\[
E^1_i := \{\imath(f) : f \in \Lambda^{e_i}\} \sqcup \{\jmath(g) : g \in
\Gamma^{e_i}\}
\]
to be a copy of the disjoint union $\Lambda^{e_i} \sqcup \Gamma^{e_i}$. We
define $E^1_{k+1}$ to be a copy of $\Gamma^0 \times \{1, \dots, n\}$ which
is disjoint from $E^0$ and each of the other $E^1_i$, and use formal
symbols $\{e(v,l) : v \in \Gamma^0, 1 \le l \le m\}$ to denote its
elements. For $1 \le i \le k$, define range and source maps $r,s : E^1_i
\to E^0$ by $r(\imath(f)) := \imath(r(f))$, $s(\imath(f)) :=
\imath(s(f))$, $r(\jmath(g)) := \jmath(r(g))$ and $s(\jmath(g)) :=
\jmath(s(g))$. Define $r,s : E^1_{k+1} \to E^0$ as in
Proposition~\ref{prp:Lambda cs Gamma}(4).

For distinct $i,j \in \{1, \dots, k+1\}$, define $E_{i,j}$ as in
Lemma~\ref{lem:FS lemma}. Define bijections $\theta_{i,j} : E_{i,j} \to
E_{j,i}$ as follows:
\begin{Itemize}
\item
For $1 \le i,j \le k$ and $(e,f) \in E_{i,j}$, we must have either $e
= \imath(a)$ and $f = \imath(b)$ for some composable pair $(a,b) \in
\Lambda^{e_i} \times_{\Lambda^0} \Lambda^{e_j}$, or else $e =
\jmath(a)$ and $f = \jmath(b)$ for some composable pair $(a,b) \in
\Gamma^{e_i} \times_{\Gamma^0} \Gamma^{e_j}$. If $e = \imath(a)$ and
$f = \imath(b)$, the factorisation property in $\Lambda$ yields a
unique pair $b' \in \Lambda^{e_j}$, $a' \in \Lambda^{e_i}$ such that
$ab = b'a'$, and we then define $\theta_{i,j}(e,f) =
(\imath(b'),\imath(a'))$. If $e = \jmath(a)$ and $f = \jmath(b)$, we
define $\theta_{i,j}(e,f)$ similarly using the factorisation property
in $\Gamma$.
\item
For $1 \le i \le k$, and $(e,f) \in E_{k+1,i}$, we have $f = \jmath(b)$
and $e = e(r(b),l)$ for some $b \in \Gamma^{e_i}$ and $1 \le l \le m$.
Define $\theta_{k+1,i}(e,f) := (\imath(p(b)), e(s(f),
\cocycle(f)^{-1}l))$.
\item
For $1 \le i \le k$, to define $\theta_{i, k+1}$, first note that if
$(f', e') = \theta_{k+1, i}(e,f)$, then $e' = e(w,l)$ for some $w \in
\Gamma^0$ and $l \in \{1, \dots, m\}$ such that $p(w) = s(f')$, $f$
is the unique lift of $f'$ such that $s(f) = \jmath(w)$, and $e =
e(r(f), \cocycle(f)l)$. It follows that $\theta_{k+1, i}$ is a
bijection and we may define $\theta_{i, k+1} := \theta_{k+1,
i}^{-1}$.
\end{Itemize}

Since $\Lambda$ and $\Gamma$ are $k$-graphs, the maps $\theta_{i,j}$, $1
\le i,j \le k$ are bijections with $\theta_{j,i} = \theta_{i,j}^{-1}$, and
we have $\theta_{i, k+1} = \theta_{k+1,i}^{-1}$ by definition, so to
invoke Lemma~\ref{lem:FS lemma}, we just need to establish
equation~\eqref{associativity}.

Equation~\eqref{associativity} holds when $h,i,j \le k$ because
$\Lambda$ and $\Gamma$ are both $k$-graphs. Suppose one of
$h,i,j = k+1$. Fix edges $f_h \in E^1_h$, $f_i \in E^1_i$ and
$f_j \in E^1_j$. First suppose that $h = k+1$; so $f_h =
e(r(f_i), l)$ for some $l$, and $f_i$ and $f_j$ both belong to
$\jmath(\Gamma)$. Apply the factorisation property for $\Gamma$
to obtain $f'_j$ and $f'_i$ such that $f'_i \in E^1_i$, $f'_j
\in E^1_j$ and $f'_j f'_i = f_i f_j$. We then have
$\theta_{i,j}(f_i, f_j) = (f_j', f_i')$. If we write $\tilde p$
for the map from $\{\jmath(f) : f \in \bigcup^k_{i=1}
\Gamma^{e_i}\}$ to $\{\imath(f) : f \in \bigcup^k_{i=1}
\Lambda^{e_i}\}$ given by $\tilde p(\jmath(\lambda)) :=
\imath(p(\lambda))$, then the properties of the covering map
imply that $\theta_{i,j}(\tilde p(f_i), \tilde p(f_j)) =
(\tilde p(f_j'), \tilde p(f_i)')$. Now
\begin{align}
(\theta_{i,j} \times \id)(\id \times \theta_{h,j})
& (\theta_{h,i} \times \id)(f_h, f_i, f_j)\nonumber\\
& = (\theta_{i,j} \times \id)(\id \times \theta_{h,j})(\tilde
p(f_i),e(s(f_i), \cocycle(f_i)^{-1}l), f_j) \nonumber\\
& = (\theta_{i,j} \times \id)(\tilde p(f_i), \tilde p(f_j), e(s(f_j),
\cocycle(f_j)^{-1}(\cocycle(f_i)^{-1})l)) \nonumber\\
& = (\tilde p(f_j'), \tilde p(f_i'), e(s(f_j), \cocycle(f_i f_j)^{-1}
l),\label{eq:factorising}
\end{align}
where, in the last equality,
$\cocycle(f_j)^{-1}\cocycle(f_i)^{-1} = \cocycle(f_i f_j)^{-1}$
by the cocycle property. On the other hand,
\begin{align*}
(\id \times \theta_{h,i})(\theta_{h,j} \times \id)&(\id \times
\theta_{i,j})(f_h, f_i, f_j) \\
&= (\id \times \theta_{h,i})(\theta_{h,j} \times \id)(f_h, f_j',
f_i') \\
&= (\id \times \theta_{h,i})(\tilde p(f_j'), e(s(f_j),
\cocycle(f_j')^{-1} l), f_i') \\
&= (\tilde p(f_j'), \tilde p(f_i'),
\cocycle(f_i')^{-1}(\cocycle(f_j')^{-1} l)) \\
&= (\tilde p(f_j'), \tilde p(f_i'), \cocycle(f_j'f_i')^{-1}l).
\end{align*}
Since $f_j' f_i' = f_i f_j$, this establishes~\eqref{associativity} when
$h = k+1$ and $1 \le i,j \le k$. Similar calculations
establish~\eqref{associativity} when $i = k+1$ and when $j = k+1$.

By Lemma~\ref{lem:FS lemma}, there is a unique $(k+1)$-graph
$\Lambda \cs{p,\cocycle} \Gamma$ with $(\Lambda \cs{p,\cocycle}
\Gamma)^0 = E^0$, $(\Lambda \cs{p,\cocycle} \Gamma)^{e_i} =
E^1_i$ for all $i$ and with commuting squares determined by the
$\theta_{i,j}$. Since the $\theta_{i,j}$, $1 \le i,j \le k$
agree with the factorisation properties in $\Gamma$ and
$\Lambda$, the uniqueness assertion of Lemma~\ref{lem:FS lemma}
applied to paths consisting of edges in $E^1_1 \cup \dots \cup
E^1_k$ shows that $\imath$ and $\jmath$ extend uniquely to
injective functors from $\Lambda$ and $\Gamma$ to $(\Lambda
\cs{p,\cocycle} \Gamma)^{(\NN^k,0_1)}$ which satisfy
Proposition~\ref{prp:Lambda cs Gamma}(2). Assertions
(3)~and~(4) of Proposition~\ref{prp:Lambda cs Gamma} follow
from the definition of $E^1_{k+1}$, and the last assertion~(5)
is established by factorising $\lambda$ into edges from the
$E^1_i$, $1 \le i \le k$ and then performing calculations
like~\eqref{eq:factorising}.

Now suppose that $p$ is finite. Then $\Gamma$ is row-finite if and
only if $\Lambda$ is, and in this case, $\Lambda \cs{p, \cocycle}
\Gamma$ is also row-finite because $p$ is locally bijective and $m <
\infty$. That $p$ is locally bijective shows that $\Lambda$ is
locally convex if and only if $\Gamma$ is. Suppose that $\Gamma$ is
locally convex. Fix $1 \le i < j \le k+1$, $a \in (\Lambda \cs{p,
\cocycle} \Gamma)^{e_i}$ and $b \in (\Lambda \cs{p, \cocycle}
\Gamma)^{e_j}$ with $r(a) = r(b)$. If $j < k+1$ then $a$ and $b$ can
be extended to paths of degree $e_i + e_j$ because $\Lambda$ and
$\Gamma$ are locally convex. If $j = k+1$, then $b = e(v,l)$ for some
$v \in \Gamma^0$ and $1 \le l \le m$. Let $a'$ be the lift of $a$
such that $r(a') = s(v)$, then $b a'$ and $a e(s(a'), l)$ extend $a$
and $b$ to paths of degree $e_i + e_j$. It follows that $\Lambda
\cs{p} \Gamma$ is locally convex.
\end{proof}

\begin{cor}\label{cor:tower graph}
Fix $N \ge 2$ in $\NN \cup \{\infty\}$. Let $(\Lambda_n, \Lambda_{n+1},
p_n, m_n, \cocycle_n)^{N-1}_{n=1}$ be a sequence of covering systems of
$k$-graphs. Then there is a unique $(k+1)$-graph $\mathbf{\LAMBDA}$ such
that $\LAMBDA^{e_i} = \bigsqcup^N_{n=1} \Lambda_n^{e_i}$ for $1 \le i \le
k$, $\LAMBDA^{e_{k+1}} = \bigsqcup^{N-1}_{n=1} (\Lambda_n \cs{p_n,
\cocycle_n} \Lambda_{n+1})^{e_{k+1}}$, and such that range, source and
composition are all inherited from the $\Lambda_n \cs{p_n, \cocycle_n}
\Lambda_{n+1}$.

If each $(\Lambda_n, \Lambda_{n+1}, p_n, m_n, \cocycle_n)$ is row-finite
then $\LAMBDA$ is row-finite. If each $\Lambda_n$ is locally convex, so is
$\LAMBDA$, and if each $\Lambda_n$ is connected, so is $\LAMBDA$.
\end{cor}
\begin{proof}
For the first part we just apply Lemma~\ref{lem:FS lemma}; the hypotheses
follow automatically from the observation that if $h,i,j$ are distinct
elements of $\{1, \dots, k+1\}$ then each path of degree $e_h + e_i + e_j$
lies in some $\Lambda_n \cs{p_n, \cocycle_n} \Lambda_{n+1}$, and these are
all $(k+1)$-graphs by Proposition~\ref{prp:Lambda cs Gamma}.

The arguments for row-finiteness, local convexity and connectedness are
the same as those in Proposition~\ref{prp:Lambda cs Gamma}.
\end{proof}

\begin{ntn}
When $N$ is finite, the $(k+1)$-graph $\LAMBDA$ of the previous
lemma will henceforth be denoted $\Lambda_1 \cs{p_1,
\cocycle_1} \cdots \cs{p_{N-1},\cocycle_{N-1}} \Lambda_N$. If
$N = \infty$, we instead denote $\LAMBDA$ by $\tgrphlim
(\Lambda_n; p_n, \cocycle_n)$.
\end{ntn}

\subsection
{Matrices of covering systems}\label{subsec:matrix of cs}

In this subsection, we generalise our construction to allow for a
different covering system $(\Lambda_j, \Gamma_i, p_{i,j}, m_{i,j},
\cocycle_{i,j})$ for each pair of connected components $\Lambda_j \subset
\Lambda$ and $\Gamma_i \subset \Gamma$. The objective is to recover the
example of the irrational rotation algebras \cite[Example~6.5]{PRRS}.

\begin{dfn}
Fix nonnegative integers $c_\Lambda, c_\Gamma \in \NN \setminus\{0\}$. A
\emph{matrix of covering systems} $(\Lambda_j, \Gamma_i, m_{i,j}, p_{i,j},
\cocycle_{i,j})^{c_\Gamma, c_\Lambda}_{i,j=1}$ consists of:
\begin{Enumerate}
\item
$k$-graphs $\Lambda$ and $\Gamma$ which decompose into connected
components $\Lambda = \bigsqcup_{j = 1, \dots, c_\Lambda} \Lambda_j$ and
$\Gamma = \bigsqcup_{i = 1, \dots, c_\Gamma} \Gamma_i$;
\item
a matrix $(m_{i,j})^{c_\Gamma, c_\Lambda}_{i,j=1} \in M_{c_\Gamma,
c_\Lambda}(\NN)$ with no zero rows or columns; and
\item
for each $i,j$ such that $m_{i,j} \not= 0$, a covering system $(\Lambda_i,
\Gamma_j, p_{i,j}, m_{i,j}, \cocycle_{i,j})$ of $k$-graphs.
\end{Enumerate}
\end{dfn}

\begin{prop}\label{prp:multiple covering}
Fix nonnegative integers $c_\Lambda, c_\Gamma \in \NN \setminus\{0\}$ and
a matrix of covering systems $(\Lambda_j, \Gamma_i, m_{i,j}, p_{i,j},
\cocycle_{i,j})^{c_\Gamma, c_\Lambda}_{i,j=1}$. Then there is a unique
$(k+1)$-graph
\[\textstyle
\big(\bigsqcup\Lambda_j\big) \cs{p, \cocycle} \big(\bigsqcup\Gamma_i\big)
\]
such that
\[\textstyle
\Big(\big(\bigsqcup\Lambda_j\big) \cs{p, \cocycle}
\big(\bigsqcup\Gamma_i\big)\Big)^{e_{k+1}} = \bigsqcup_{i,j} (\Lambda_j
\cs{p_{i,j}, \cocycle_{i,j}} \Gamma_i)^{e_{k+1}},
\]
each $\big(\big(\bigsqcup\Lambda_j\big) \cs{p, \cocycle}
\big(\bigsqcup\Gamma_i\big)\big)^{e_l}$ for $1 \le l \le k$ is equal to
$\Lambda^{e_l} \sqcup \Gamma^{e_l}$ and the commuting squares are
inherited from the $\Lambda_j \cs{p_{i,j}, \cocycle_{i,j}} \Gamma_i$.

If each $(\Lambda_i, \Gamma_j, p_{i,j}, m_{i,j},
\cocycle_{i,j})$ is row finite then
$\big(\bigsqcup\Lambda_j\big) \cs{p, \cocycle}
\big(\bigsqcup\Gamma_i\big)$ is row finite. If $\Lambda$ and
$\Gamma$ are locally convex, so is
$\big(\bigsqcup\Lambda_j\big) \cs{p, \cocycle}
\big(\bigsqcup\Gamma_i\big)$.
\end{prop}
\begin{proof}
We apply Lemma~\ref{lem:FS lemma}; since the commuting squares are
inherited from the $\Lambda_j \cs{p_{i,j}, \cocycle_{i,j}} \Gamma_i$, they
satisfy the associativity condition~\eqref{associativity} because each
$\Lambda_j \cs{p_{i,j}, \cocycle_{i,j}} \Gamma_i$ is a $(k+1)$-graph.
\end{proof}

\begin{cor}\label{cor:pasting}
Fix $N \ge 2$ in $\NN \cup \{\infty\}$. Let $(c_n)^N_{n=1}
\subset \NN\setminus\{0\}$ be a sequence of positive integers.
For $1 \le n < N$, let $(\Lambda_{n,j}, \Lambda_{n+1, i},
p^n_{i,j}, m^n_{i,j}, \cocycle^n_{i,j})^{c_{n+1}, c_n}_{i,j =
1}$ be a matrix of covering systems. Then there is a unique
$(k+1)$-graph $\LAMBDA$ such that $\LAMBDA^{e_i} =
\bigcup^N_{n=1} \bigcup^{c_n}_{j=1}\Lambda_{n,j}^{e_i}$ for $1
\le i \le k$, $\LAMBDA^{e_{k+1}} = \bigcup^{N-1}_{n=1}
\big(\big(\bigsqcup^{c_n}_{j=1} \Lambda_{n,j}\big) \cs{p^n,
\cocycle^n} \big(\bigsqcup^{c_{n+1}}_{i=1}
\Lambda_{n+1,i}\big)\big)^{e_{k+1}}$, and range, source and
composition are inherited from the $\big(\bigsqcup^{c_n}_{j=1}
\Lambda_{n,j}\big) \cs{p^n, \cocycle^n}
\big(\bigsqcup^{c_{n+1}}_{i=1} \Lambda_{n+1,i}\big)$.

If each $(\Lambda_{n,j}, \Lambda_{n+1,i}, p^n_{i,j}, m^n_{i,j},
\cocycle^n_{i,j})$ is row finite, then $\LAMBDA$ is row finite.
If each $\Lambda_n$ is locally convex, so is $\LAMBDA$.
\end{cor}

\begin{example}[The Irrational Rotation algebras]\label{eg:IR algs}
Fix $\theta \in [0,1]\setminus\QQ$. Let $[a_1, a_2, \dots]$ be the simple
continued fraction expansion of $\theta$. For each $n$, let $c_n = 2$, let
$\phi_n := \big(\begin{smallmatrix} a_n & 1 \\ 1
&0\end{smallmatrix}\big)$, and let $m^n := (m^n_{i,j})^2_{i,j=1}$ be the
matrix product $\phi_{T(n+1)} \cdots \phi_{T(n) + 1}$ where $T(n) :=
n(n+1)/2$ is the $n^{\rm th}$ triangular number. Of all the integers
$m^n_{i,j}$ obtained this way, only $m^1_{1,2}$ is equal to zero, so the
matrices $m^n$ have no zero rows or columns. Whenever $m^n_{i,j} \not= 0$,
let $\cocycle^n_{i,j}$ be the permutation of the set $\{1, \dots,
m^n_{i,j}\}$ given by $\cocycle^n_{i,j}l = l+1$ if $1 \le l < m^n_{i,j}$,
and $\cocycle^n_{i,j}m^n_{i,j} = 1$.

Let $\Lambda_{n,i}$, $n \in \NN \setminus \{0\}$, $i = 1,2$ be mutually
disjoint copies of the $1$-graph $T_1$ whose skeleton consists of a single
vertex and a single directed edge. For each $n$, let $\Lambda_n$ be the
$1$-graph $\Lambda_{n,1} \sqcup \Lambda_{n,2}$ so that for each $n$,
$(\Lambda_{n,j}, \Lambda_{n+1,i}, p^n_{i,j}, m^n_{i,j},
\cocycle^n_{i,j})^2_{i,j = 1}$ is a matrix of covering systems.

\begin{figure}[ht]
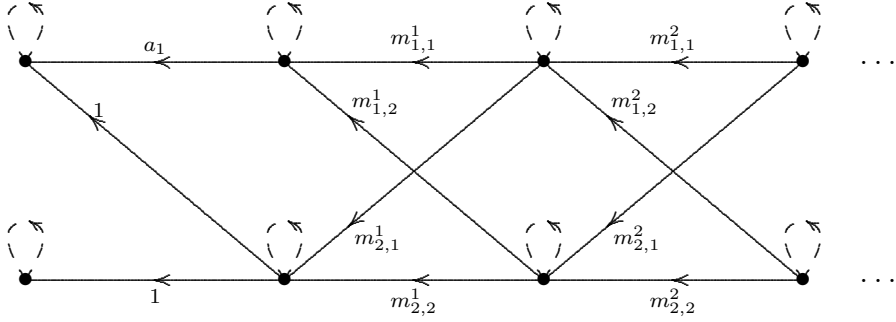

\(
\beginpicture \setcoordinatesystem units <1.65em,1.4em>
    \setquadratic \setdashes <0.417em>
    \put{$\bullet$} at 0 0
    \put{
        \plot 0 0 -0.34 .96 0 1.33 /
        \plot 0 0 0.34  .96 0 1.33 /
        \arrow <0.15cm> [0.25,0.75] from 0.1  1.314 to 0.08 1.328
    }[b] at -0.1 0
    \put{$\bullet$} at 0 5
    \put{
        \plot 0 0 -0.34 .96 0 1.33 /
        \plot 0 0 0.34  .96 0 1.33 /
        \arrow <0.15cm> [0.25,0.75] from 0.1  1.314 to 0.08 1.328
    }[b] at  -0.1 5
    \put{$\bullet$} at 5 0
    \put{
        \plot 0 0 -0.34 .96 0 1.33 /
        \plot 0 0 0.34  .96 0 1.33 /
        \arrow<0.15cm> [0.25,0.75] from 0.1  1.314 to 0.08 1.328
    }[b] at 4.9 0
    \put{$\bullet$} at 5 5
    \put{
        \plot 0 0 -0.34 .96 0 1.33 /
        \plot 0 0 0.34  .96 0 1.33 /
        \arrow <0.15cm> [0.25,0.75] from 0.1  1.314 to 0.08 1.328
    }[b] at 4.9 5
    \put{$\bullet$} at 10 0
    \put{
        \plot 0 0 -0.34 .96 0 1.33 /
        \plot 0 0 0.34  .96 0 1.33 /
        \arrow <0.15cm> [0.25,0.75] from 0.1  1.314 to 0.08 1.328
    }[b] at 9.9 0
    \put{$\bullet$} at 10 5
    \put{
        \plot 0 0 -0.34 .96 0 1.33 /
        \plot 0 0 0.34  .96 0 1.33 /
        \arrow <0.15cm> [0.25,0.75] from 0.1  1.314 to 0.08 1.328
    }[b] at 9.9 5
    \put{$\bullet$} at 15 0
    \put{
        \plot 0 0 -0.34 .96 0 1.33 /
        \plot 0 0 0.34  .96 0 1.33 /
        \arrow <0.15cm> [0.25,0.75] from 0.1  1.314 to 0.08 1.328
    }[b] at 14.9 0
    \put{$\bullet$} at 15 5
    \put{
        \plot 0 0 -0.34 .96 0 1.33 /
        \plot 0 0 0.34  .96 0 1.33 /
        \arrow <0.15cm> [0.25,0.75] from 0.1  1.314 to 0.08 1.328
    }[b] at 14.9 5
    \put{$\dots$} at 16.5 0
    \put{$\dots$} at 16.5 5
    \setlinear \setsolid
    \plot 0 0 15 0 /
    \arrow <0.5em> [0.25,0.75] from 2.51 0 to 2.49 0
    \arrow <0.5em> [0.25,0.75] from 7.51 0 to 7.49 0
    \arrow <0.5em> [0.25,0.75] from 12.51 0 to 12.49 0
    \plot 0 5 15 5 /
    \arrow <0.5em> [0.25,0.75] from 2.51 5 to 2.49 5
    \arrow <0.5em> [0.25,0.75] from 7.51 5 to 7.49 5
    \arrow <0.5em> [0.25,0.75] from 12.51 5 to 12.49 5
    \plot 0 5 5 0 /
    \plot 5 0 10 5 /
    \plot 5 5 10 0 /
    \plot 10 0 15 5 /
    \plot 10 5 15 0 /
    \arrow <0.5em> [0.25,0.75] from 1.26 3.74 to 1.24 3.76
    \put{$\scriptstyle a_1$} [b] at 2.5 5.2
    \put{$\scriptstyle 1$} [t] at 2.5 -0.2
    \put{$\scriptstyle 1$} [bl] at 1.3 3.75
    \arrow <0.5em> [0.25,0.75] from 6.26 1.26 to 6.24 1.24
    \arrow <0.5em> [0.25,0.75] from 6.26 3.74 to 6.24 3.76
    \put{$\scriptstyle m^1_{1,1}$} [b] at 7.5 5.2
    \put{$\scriptstyle m^1_{2,2}$} [t] at 7.5 -0.2
    \put{$\scriptstyle m^1_{1,2}$} [bl] at 6.3 3.75
    \put{$\scriptstyle m^1_{2,1}$} [tl] at 6.35 1.25
    \arrow <0.5em> [0.25,0.75] from 11.26 1.26 to 11.24 1.24
    \arrow <0.5em> [0.25,0.75] from 11.26 3.74 to 11.24 3.76
    \put{$\scriptstyle m^2_{1,1}$} [b] at 12.5 5.2
    \put{$\scriptstyle m^2_{2,2}$} [t] at 12.5 -0.2
    \put{$\scriptstyle m^2_{1,2}$} [bl] at 11.3 3.75
    \put{$\scriptstyle m^2_{2,1}$} [tl] at 11.35 1.25
\endpicture
\) \caption{A tower of coverings with multiplicities}\label{fig:IR}
\end{figure}
Modulo relabelling the generators of $\NN^2$, the $2$-graph
$\tgrphlim\left(\bigsqcup_{j=1}^{c_n} \Lambda_{n,j}; p^n_{i,j},
\cocycle^n_{i,j}\right)$ obtained from this data as in
Corollary~\ref{cor:pasting} is precisely the rank-2 Bratteli
diagram of \cite[Example~6.5]{PRRS} whose $C^*$-algebra is
Morita equivalent to the irrational rotation algebra
$A_\theta$. Figure~\ref{fig:IR} is an illustration of its
skeleton (parallel edges drawn as a single edge with a label
indicating the multiplicity). The factorisation rules are all
of the form $f g = \sigma(g)f'$ where $f$ and $f'$ are the
dashed loops at either end of a solid edge in the diagram, and
$\sigma$ is a transitive permutation of the set of edges with
the same range and source as $g$.

More generally, Section~7 of \cite{PRRS} considers in some detail the
structure of the $C^*$-algebras associated to rank-2 Bratteli diagrams
with length-1 cycles. All such rank-2 Bratteli diagrams can be recovered
as above from Corollary~\ref{cor:pasting}.
\end{example}

\section{$C^*$-algebras associated to covering systems of
$k$-graphs}\label{sec:C*algs}

In this section, we describe how a covering system
$(\Lambda,\Gamma,p,m,\cocycle)$ induces an inclusion of
$C^*$-algebras $C^*(\Lambda) \hookrightarrow M_m(C^*(\Gamma))$
and hence a homomorphism of $K$-groups $K_*(C^*(\Lambda)) \to
K_*(C^*(\Gamma))$. The main result of the section is
Theorem~\ref{thm:direct limit} which shows how to use these
maps to compute the $K$-theory of $C^*(\tgrphlim(\Lambda_n;
p_n, \cocycle_n))$ from the data in a sequence $(\Lambda_n,
\Lambda_{n+1}, p_n, m_n, \cocycle_n)^\infty_{n=1}$ of covering
systems.

The following definition of the Cuntz-Krieger algebra of a row-finite
locally convex $k$-graph $\Lambda$ is taken from \cite[Definition
3.3]{rsy}.

Given a row-finite, locally convex $k$-graph $( \Lambda, d )$, a
Cuntz-Krieger $\Lambda$-family is a collection $\{t_\lambda : \lambda \in
\Lambda \}$ of partial isometries satisfying
\begin{itemize}
\item[(CK1)] $\{t_v : v \in \Lambda^0\}$ is a collection of mutually
orthogonal projections; \item[(CK2)] $t_\lambda t_\mu = t_{\lambda \mu}$
whenever $s(\lambda) = r(\mu)$; \item[(CK3)] $t^*_\lambda t_\lambda =
t_{s(\lambda)}$ for all $\lambda \in \Lambda$; and
\item[(CK4)] $t_v = \sum_{\lambda \in v \Lambda^{\le n}} t_\lambda
t^*_\lambda$ for all $v \in \Lambda^0$ and $n \in \NN^k$.
\end{itemize}
The \emph{Cuntz-Krieger algebra} $C^*(\Lambda)$ is the $C^*$-algebra
generated by a Cuntz-Krieger $\Lambda$-family $\{s_\lambda : \lambda \in
\Lambda \}$ which is universal in the sense that for every Cuntz-Krieger
$\Lambda$-family $\{t_\lambda : \lambda \in \Lambda \}$ there is a unique
homomorphism $\pi_t$ of $C^*( \Lambda )$ satisfying $\pi_t(s_\lambda) =
t_\lambda$ for all $\lambda \in \Lambda$.

\begin{rmks}
If $\Lambda$ has no sources (that is $v\Lambda^n \not= \emptyset$ for all
$v \in \Lambda^0$ and $n \in \NN^k$), then $\Lambda$ is automatically
locally convex, and the definition of $C^* ( \Lambda )$ given above
reduces to the original definition \cite[Definition 1.5]{kp}.

By \cite[Theorem 3.15]{rsy} there is a Cuntz-Krieger $\Lambda$-family
$\{t_\lambda : \lambda \in \Lambda\}$ such that $t_\lambda \neq 0$ for all
$\lambda \in \Lambda$. The universal property of $C^*(\Lambda)$ therefore
implies that the generating partial isometries $\{s_\lambda : \lambda \in
\Lambda\} \subset C^*(\Lambda)$ are all nonzero.
\end{rmks}

Let $\Xi$ be a $k$-graph. The universal property of $C^*(\Xi)$
gives rise to an action $\gamma$ of $\TT^k$ on $C^*(\Xi)$,
called the \emph{gauge-action} (see, for example \cite[\S
4.1]{rsy}), such that $\gamma_z(s_\xi) = z^{d(\xi)} s_\xi$ for
all $z \in \TT^k$ and $\xi \in \Xi$.

\begin{prop} \label{prp:embedding}
Let $(\Lambda, \Gamma, p, m, \cocycle)$ be a row-finite
covering system of locally convex $k$-graphs. Let
$\gamma_\Lambda$ and $\gamma_\Gamma$ denote the gauge actions
of $\TT^k$ on $C^*(\Lambda)$ and $C^*(\Gamma)$, and let
$\gamma$ denote the gauge action of $\TT^{k+1}$ on $C^*(\Lambda
\cs{p,\cocycle} \Gamma)$.
\begin{Enumerate}
\item The inclusions $\imath : \Lambda \to \Lambda
\cs{p,\cocycle} \Gamma$ and $\jmath : \Gamma \to \Lambda \cs{p,\cocycle}
\Gamma$ induce embeddings of $C^* ( \Lambda )$ and $C^* ( \Gamma )$ in
$C^* ( \Lambda \cs{p,\cocycle} \Gamma )$ characterised by
\[
\imath_*(s_\alpha) = s_{\imath(\alpha)} \text{ and } \jmath_*(s_\beta) =
s_{\jmath(\beta)}\quad\text{ for $\alpha \in \Lambda$ and $\beta \in
\Gamma$.}
\]

\item The sum $\sum_{v \in \jmath (\Gamma^0)} s_v$ converges
strictly to a full projection $Q \in \mathcal{M}(C^*(\Lambda
\cs{p,\cocycle} \Gamma))$, and the range of $\jmath_*$ is $Q C^*(\Lambda
\cs{p,\cocycle} \Gamma) Q$.

\item For $1 \le i \le m$, the sum $\sum_{v \in \Gamma^0}
s_{e(v,i)}$ converges strictly to a partial isometry $V_i \in \mathcal{M}
(C^*(\Lambda \cs{p,\cocycle} \Gamma))$. The sum $\sum_{v \in
\imath(\Lambda^0)} s_v$, converges strictly to the full projection $P :=
\sum^m_{i=1} V_i V^*_i \in \mathcal{M}(C^*(\Lambda \cs{p,\cocycle}
\Gamma))$. Moreover, $\imath_*$ is a nondegenerate homomorphism into $P
C^*(\Lambda \cs{p,\cocycle} \Gamma) P$.

\item 
There is an isomorphism $\phi : M_m(C^*(\Gamma)) \to P C^*(\Lambda
\cs{p,\cocycle} \Gamma)P$ such that
\[
\phi\big((a_{i,j})^m_{i,j=1}\big) = \sum^m_{i,j=1} V_i \jmath_*(a_{i,j})
V_j^*.
\]

\item There is an embedding $\iota_{p,\cocycle} : C^*(\Lambda)
    \to M_m(C^*(\Gamma))$ such that $\phi \circ
    \iota_{p,\cocycle} = \imath_*$. The embedding $\iota_{p,\cocycle}$ is
    equivariant in $\gamma_\Lambda$ and the action $\id_m
    \otimes \gamma_\Gamma$ of $\TT^k$ on $M_m(C^*(\Gamma))$ by
    coordinate-wise application of $\gamma_\Gamma$.

\item If we identify $K_*(C^*(\Gamma))$ with
    $K_*(M_m(C^*(\Gamma)))$, then the induced homomorphism
    $(\iota_{p,\cocycle})_*$ may be viewed as a map from
    $K_*(C^*(\Lambda)) \to K_*(C^*(\Gamma))$. When applied to
    $K_0$-classes of vertex projections, this map satisfies
\[
(\iota_{p,\cocycle})_*([s_v]) = \sum_{p(u) = v} m \cdot [s_u] \in
K_0(C^*(\Gamma)).
\]
\end{Enumerate}
\end{prop}

The proofs of the last three statements require the following general
Lemma. This is surely well-known but we include it for completeness.

\begin{lem} \label{lem:matrix alg}
Let $A$ be a $C^*$-algebra, let $q \in \mathcal{M}(A)$ be a
projection, and suppose that $v_1 , \ldots , v_n \in
\mathcal{M} (A)$ satisfy $v_i^* v_j = \delta_{i,j} q$ for $1
\le i,j \le n$. Then $p = \sum_{i=1}^n v_i v_i^*$ is a
projection and $pAp \cong M_n (qAq)$.
\end{lem}

\begin{proof}
That $v^*_i v_j = \delta_{i,j} q$ implies that the $v_i$ are partial
isometries with mutually orthogonal range projections $v_i v^*_i$. Hence
$p$ is a projection in $\Mm(A)$. Define a map $\phi$ from $pAp$ to
$M_n(qAq)$ as follows: for $a \in pAp$ and $1 \le i,j \le n$, let $a_{i,j}
:= v_i^* a v_j$, and define $\phi(a)$ to be the matrix $\phi(a) =
(a_{i,j})^n_{i,j=1}$.

It is straightforward to check using the properties of the $v_i$ that
$\phi$ is a $C^*$-homomorphism. It is an isomorphism because the
homomorphism $\psi : M_n(qAq) \to pAp$ defined by
\[
\psi\big((a_{i,j})_{i,j=1}^n\big) := \sum_{i,j=1}^n v_i a_{ij} v_j^* \in
qAq
\]
is an inverse for $\phi$.
\end{proof}

\begin{proof}[Proof of Proposition~\ref{prp:embedding}]
(1) The collection $\{ s_{\imath ( \lambda )} : \lambda \in
\Lambda \}$ forms a Cuntz-Krieger $\Lambda$-family in $C^* (
\Lambda \cs{p,\cocycle} \Gamma )$, and so by the universal
property of $C^* ( \Lambda )$ induces a homomorphism $\imath_*
: C^* ( \Lambda ) \to C^* ( \Lambda \cs{p,\cocycle} \Gamma )$.
For $z \in \TT^k$, write $(z,1)$ for the element $(z_1, \dots,
z_k, 1) \in \TT^{k+1}$. Recall that $\gamma$ denotes the gauge
action of $\TT^{k+1}$ on $C^*(\Lambda \cs{p,\cocycle} \Gamma)$.
Then the action $z \mapsto \gamma_{(z,1)}$ of $\TT^k$ on $C^*
(\Lambda \cs{p,\cocycle} \Gamma)$ satisfies
\[
\imath_* ((\gamma_\Lambda)_z(a)) = \gamma_{(z,1)} (\imath_*(a))
\]
for all $a \in C^* ( \Lambda )$ and $z \in \TT^k$. Since
$\imath_* ( s_v ) = s_{\imath (v)} \neq 0$ for all $v \in
\Lambda^0$ it follows from the gauge-invariant uniqueness
theorem \cite[Theorem 2.1]{kp} that $\imath_*$ is injective. A
similar argument applies to $\jmath_*$.

(2) As the projections $s_v$, $v \in \jmath ( \Gamma^0 )$ are mutually
orthogonal, a standard argument shows that the sum $\sum_{v \in
\jmath(\Gamma^0)} s_v$ converges to a projection $Q$ in the multiplier
algebra (see \cite[Lemma 2.1]{r}). The range of $\jmath_*$ is equal to $Q
C^*(\Lambda \cs{p,\cocycle} \Gamma) Q$ because $\jmath(\Gamma^0) (\Lambda
\cs{p,\cocycle} \Gamma) \jmath(\Gamma^0) = \jmath(\Gamma)$. To see that
$Q$ is full, it suffices to show that every generator of $C^*(\Lambda
\cs{p,\cocycle} \Gamma)$ belongs to the ideal $I(Q)$ generated by $Q$. So
let $\alpha \in \Lambda \cs{p,\cocycle} \Gamma$. Either $s(\alpha) \in
\jmath(\Gamma^0)$ or $s(\alpha) \in \imath(\Lambda^0)$. If $s(\alpha) \in
\jmath(\Gamma^0)$, then $s_\alpha = s_\alpha Q \in I(Q)$. On the other
hand, if $s(\alpha) \in \imath(\Lambda^0)$, the Cuntz-Krieger relation
ensures that
\[\textstyle
s_\alpha = \sum_{p(w) = s(\alpha)} \sum^m_{i=1} s_\alpha s_{e(w,i)} Q
s^*_{e(w,i)},
\]
which also belongs to $I(Q)$.

(3) For fixed $i$, the partial isometries $s_{e(v,i)}$ have
mutually orthogonal range projections and mutually orthogonal
source projections. Hence an argument similar that of
\cite[Lemma~2.1]{r} shows that $\sum_{v \in \Gamma^0}
s_{e(v,i)}$ converges strictly to a multiplier $V_i \in
\Mm(C^*(\Lambda \cs{p,\cocycle} \Gamma))$. A simple calculation
shows that $V^*_i V_j = \delta_{i,j} Q$ for all $i,j$. Hence
each $V_i$ is a partial isometry, and $P$ is full because $Q$
is full. The homomorphism $\imath_*$ is nondegenerate because
the net
\[\textstyle
\Big(\imath_*\big(\sum_{v \in F} s_v\big)\Big)_{F \subset \Lambda^0\text{
finite}}
\]
converges strictly to $P \in \mathcal{M}(C^*(\Lambda \cs{p,\cocycle}
\Gamma))$.

(4) This follows directly from Part~(3) and Lemma~\ref{lem:matrix alg}.

(5) We define $\iota_{p,\cocycle} := \phi^{-1} \circ \imath_*$. For the
gauge-equivariance, recall that $\imath_*$ (respectively $\jmath_*$) are
equivariant in $\gamma|_{(\TT^k,1)}$ and $\gamma_\Lambda$ (respectively
$\gamma_\Gamma$). By definition, $\phi$ is equivariant in $(\id \otimes
\gamma)$ and $\gamma_{(\TT^k,1)} \circ \jmath_*$. The equivariance of
$\iota_{p,\cocycle}$ follows.

(6) By~(CK4), for $v \in \Lambda^0$ we have $s_{\imath(v)} = \sum_{f \in
v(\Lambda \cs{p,\cocycle} \Gamma)^{e_{k+1}}} s_f s^*_f$, so the
$K_0$-class $[s_{\imath(v)}]$ is equal to $\sum_{f \in v(\Lambda
\cs{p,\cocycle} \Gamma)^{e_{k+1}}} [s_f s^*_f]$. We can write $v(\Lambda
\cs{p,\cocycle} \Gamma)^{e_{k+1}}$ as the disjoint union
\[
v(\Lambda \cs{p,\cocycle} \Gamma)^{e_{k+1}} = \bigsqcup_{p(u) = v} \{e(u,
i) : 1 \le i \le m\}.
\]
In $K_0(C^*(\Lambda \cs{p,\cocycle} \Gamma))$, we have $[s_{e(u, i)}
s^*_{e(u, i)}] = [s^*_{e(u, i)} s_{e(u, i)}] = [s_{\jmath(u)}]$, and the
result follows.
\end{proof}

\begin{ntn}\label{ntn:m=1_C*}
As in Notation~\ref{ntn:m=1}, when $m = 1$ so that $\cocycle$
is trivial, we continue to drop references to $\cocycle$ at the
level of $C^*$-algebras. So Proposition~\ref{prp:embedding}(5)
gives an inclusion $\iota_p : C^*(\Lambda) \to C^*(\Gamma)$ and
the induced homomorphism of $K$-groups obtained from
Proposition~\ref{prp:embedding}(6) is denoted $(\iota_p)_* :
K_*(C^*(\Lambda)) \to K_*(C^*(\Gamma))$. This homomorphism
satisfies
\[
(\iota_p)_*([s_v]) = \sum_{p(u) = v} [s_u].
\]
\end{ntn}

When no confusion is likely to occur, we will suppress the maps
$\imath$, $\jmath$, $\imath_*$ and $\jmath_*$ and regard
$\Lambda$ and $\Gamma$ as subsets of $\Lambda \cs{p,\cocycle}
\Gamma$ and $C^*(\Lambda)$ and $C^*(\Gamma)$ as
$C^*$-subalgebras of $C^*(\Lambda \cs{p,\cocycle} \Gamma)$.

\begin{rmk}
\begin{Enumerate}
\item The isomorphism $\phi$ of
    Proposition~\ref{prp:embedding}(4) extends to an
    isomorphism $\tilde\phi : M_{m+1}(C^*(\Gamma)) \to
    C^*(\Lambda \cs{p,\cocycle} \Gamma)$ which takes the block
    diagonal matrix $\left(\begin{smallmatrix} 0_{m \times m} &
    0_{m \times 1} \\ 0_{1 \times m} &
    a\end{smallmatrix}\right)$ to $\jmath_*(a)$. To see this,
    let $V_, \dots, V_m$ be as in Proposition~\ref{prp:Lambda cs
    Gamma}(3), let $V_{m+1} = Q$, and apply
    Lemma~\ref{lem:matrix alg}.

\item If $m = 1$ then $\phi$ is an isomorphism of $C^*(\Gamma)$
onto $P C^*(\Lambda \cs{p} \Gamma) P$, and $\iota_p : C^*(\Lambda)
\hookrightarrow C^*(\Gamma)$ satisfies
\[\textstyle
\iota_p(s_\lambda) = \sum_{p(\tilde\lambda) = \lambda} s_{\tilde\lambda}.
\]
\end{Enumerate}
\end{rmk}

Fix $N \ge 2$ in $\NN$. Let $(\Lambda_n, \Lambda_{n+1}, p_n, m_n,
\cocycle_n)^{N-1}_{n=1}$ be a sequence of row-finite covering systems of
locally convex $k$-graphs. Recall that in Corollary~\ref{cor:tower graph}
we obtained from such data a $(k+1)$-graph $\Lambda_1 \cs{p_1,\cocycle_1}
\cdots \cs{p_{N-1},\cocycle_{N-1}} \Lambda_N$, which for convenience we
will denote $\LAMBDA_N$ (the subscript is unnecessary here, but will be
needed in Proposition~\ref{prp:tower matrix inclusion}). We now examine
the structure of $C^*(\LAMBDA_N)$ using Proposition~\ref{prp:embedding}.

\begin{prop} \label{prp:tower algebra}
Continue with the notation established in the previous paragraph. For each
$v \in \Lambda_N^0$, list $\LAMBDA_N^{N e_{k+1}} v$ as $\{\alpha(v,i) : 1
\le i \le M\}$ where $M = m_1 m_2 \cdots m_{N-1}$.
\begin{Enumerate}
\item
For $1 \le n \le N$, the sum $\sum_{v \in \Lambda_n^0} s_v$ converges
strictly to a full projection $P_n \in \Mm(C^*(\LAMBDA_N))$.
\item
For $1 \le i \le M$, the sum $\sum_{v \in \Lambda^0_N} s_{\alpha(v, i)}$
converges strictly to a partial isometry $V_i \in \Mm(C^*(\LAMBDA_N))$
such that $V_i^* V_i = P_N$.
\item We have $\sum_{i = 1}^M V_i V^*_i = P_1$, and there is an
isomorphism
\[
\phi : M_{M}(C^*(\Lambda_N)) \to P_1 C^*(\LAMBDA_N) P_1
\]
such that $\phi((a_{i,j})^M_{i,j=1}) = \sum^M_{i,j=1}V_i a_{i,j} V^*_j$.
\end{Enumerate}
\end{prop}
\begin{proof}
Calculations like those in parts (2)~and~(3) of
Proposition~\ref{prp:embedding} show that the sums defining the $P_n$ and
the $V_i$ converge in the multiplier algebra of $C^*(\LAMBDA_N)$ and that
each $P_n$ is full.

Since distinct paths in $\LAMBDA_N^{N e_{k+1}}$ have orthogonal range
projections and since paths in $\LAMBDA_N^{N e_{k+1}}$ with distinct
sources have orthogonal source projections, each $V^*_i V_i = P_N$, and
$\sum^M_{i=1} V_i V^*_i = P_1$.

One checks as in Proposition~\ref{prp:embedding}(1) that the inclusions
$\imath_n : \Lambda_n \hookrightarrow \LAMBDA_N$ induce inclusions
$(\imath_n)_* : C^*(\Lambda_n) \hookrightarrow P_n C^*(\LAMBDA_N)P_n$, and
in particular that $(\imath_N)_* : C^*(\Lambda_N) \to P_N
C^*(\LAMBDA_N)P_N$ is an isomorphism. The final statement follows from
Lemma~\ref{lem:matrix alg}.
\end{proof}

We now describe the inclusions of the corners determined by $P_1$ as $N$
increases. To do this, we first need some notation. Given a $C^*$-algebra
$A$, and positive integers $m,n$, we denote by $\pi_{m,n} \otimes \id_A :
M_m(M_n(A)) \to M_{mn}(A)$ the canonical isomorphism which takes the
matrix $a = \big((a_{i,j,j',i'})^n_{j,j' = 1}\big)^m_{i,i' = 1}$ to the
matrix $\pi(a)$ satisfying
\[
\pi(a)_{j + n(i-1), j' + n(i'-1)} =
a_{i,j,j',i'}\quad\text{for $1 \le i,i' \le m$, $1 \le j,j' \le n$.}
\]
Given $C^*$-algebras $A$ and $B$, a positive integer $m$, and a
$C^*$-homomorphism $\psi : A \to B$, we write $\id_m \otimes
\psi : M_m(A) \to M_m(B)$ for the $C^*$-homomorphism
\[
(\id_m \otimes \psi)\big((a_{i,j})^m_{i,j=1}\big) =
\big(\psi(a_{i,j})\big)^m_{i,j = 1}.
\]
Finally, given a matrix algebra $M_m(A)$ over a $C^*$-algebra $A$, and
given $1 \le i,i' \le m$ and $a \in A$, we write $\theta_{i,i'}a$ for the
matrix
\[
\big(\theta_{i,i'}a)_{j,j'} =
\begin{cases}
a &\text{ if $j=i$ and $j' = i'$} \\
0 &\text{ otherwise.}
\end{cases}
\]

\begin{prop} \label{prp:tower matrix inclusion}
Fix $N \ge 2$ in $\NN$. Let $(\Lambda_n, \Lambda_{n+1}, p_n, m_n,
\cocycle_n)^{N}_{n=1}$ be a sequence of row-finite covering systems of
locally convex $k$-graphs. We view the $(k+1)$-graph $\LAMBDA_N :=
\Lambda_1 \cs{p_1, \cocycle_1} \cdots \cs{p_{N-1},\cocycle_{N-1}}
\Lambda_N$ as a subcategory of $\LAMBDA_{N+1} := \Lambda_1 \cs{p_1,
\cocycle_1} \cdots \cs{p_{N},\cocycle_{N}} \Lambda_{N+1}$ and likewise
regard $C^*(\LAMBDA_N)$ as a $C^*$-subalgebra of $C^*(\LAMBDA_{N+1})$. In
particular, we view $P_1 = \sum_{v \in \Lambda_1^0} s_v$ as a projection
in both $\Mm(C^*(\LAMBDA_N))$ and $\Mm(C^*(\LAMBDA_{N+1}))$.

Let $M := m_1m_2 \dots m_{N-1}$, and let $\phi_N : M_M(C^*(\Lambda_N)) \to
P_1(C^*(\LAMBDA_N) P_1$ and $\phi_{N+1} : M_{Mm_N}(C^*(\Lambda_{N+1})) \to
P_1 C^*(\LAMBDA_{N+1}) P_1$ be the isomorphisms obtained from
Proposition~\ref{prp:tower algebra}. Then the following diagram commutes.
\[
    \beginpicture
    \setcoordinatesystem units <2.5em,2.5em>
    \put{$P_1 C^*(\LAMBDA_N) P_1$}[cb] at 0 2
    \put{$P_1 C^*(\LAMBDA_{N+1})P_1 $}[cb] at 10 2
    \put{$M_M(C^*(\Lambda_N))$}[ct] at 0 0
    \put{$M_{M m_N}(C^*(\Lambda_{N+1}))$}[ct] at 10 0
    \arrow <0.6em> [0.1, 0.3] from 0 0.1 to 0 1.9
    \put{$\phi_N$}[l] at 0.1 1
    \arrow <0.6em> [0.1, 0.3] from 10 0.1 to 10 1.9
    \put{$\phi_{N+1}$}[l] at 10.1 1
    \arrow <0.6em> [0.1, 0.3] from 1.5 2.2 to 8.5 2.2
    \circulararc -180 degrees from 1.5 2.2 center at 1.5 2.3
    \put{$\subseteq$}[b] at 4.8 2.3
    \arrow <0.6em> [0.1, 0.3] from 1.3 -0.2 to 8.1 -0.2
    \put{$(\pi_{M,m_N} \otimes \id_{C^*(\Lambda_{N+1})}) \circ (\id_M \otimes \iota_{p_N,\cocycle_N})$}[b] at 4.6
-0.1
    \endpicture
\]
\end{prop}
\begin{proof}
As in Proposition~\ref{prp:tower algebra}, write $\LAMBDA_N^{N e_{k+1}} =
\{\alpha(v,i) : v \in \Lambda_N^0, i \in \{1, \cdots, M\}\}$. For $i = 1,
\dots, M$, let $V_i := \sum_{v \in \Lambda_N^0} s_{\alpha(v,i)}$. For $j =
1, \dots, m_N$, let
\[
W_j := \sum_{w \in \Lambda_{N+1}^0}\sum^M_{i=1} s_{\alpha(p_N(w),i)}
s_{e(w,j)}.
\]
For $(i,j) \in \{1, \dots, M\} \times \{1, \dots, m_N\}$, let $U_{j +
m_N(i-1)} := \sum_{u \in \Lambda_{N+1}^0} s_{\alpha(p_N(u),i) e(u,j)}$. In
what follows, we suppress canonical inclusion maps, and regard
$C^*(\Lambda_N)$ as a subalgebra of $C^*(\LAMBDA_N)$, and both
$C^*(\LAMBDA_N)$ and $C^*(\Lambda_{N+1})$ as subalgebras of
$C^*(\LAMBDA_{N+1})$. The corner $P_1 C^*(\LAMBDA_N) P_1$ is equal to the
closed span of elements of the form $V_i a V^*_{i'}$ where $a \in
C^*(\Lambda_N)$ and $i,i' \in \{1, \dots, M\}$, and $P_1
C^*(\LAMBDA_{N+1}) P_1$ is equal to the closed span of elements of the
form $U_{l} b U^*_{l'}$ where $b \in C^*(\Lambda_{N+1})$, $l,l' \in \{1,
\dots, M m_N\}$.

We have $\phi_N\big((a_{i,i'})^M_{i,i'=1}\big) = \sum_{i,i' =
1}^M V_i a_{i,i'} V^*_{i'}$ by definition. The isomorphism
$\phi_{N+1}$ from $M_{Mm_N}(C^*(\Lambda_{N+1}))$ to $P_1
C^*(\LAMBDA_{N+1})P_1$ described in Proposition~\ref{prp:tower
algebra} satisfies
\[\textstyle
\phi_{N+1}\big(\sum_{l,l' = 1}^{Mm_N} U_l b_{l,l'} U^*_{l'}\big) =
\big(b_{l,l'}\big)^{Mm_N}_{l,l'=1}.
\]

The Cuntz-Krieger relations show that $V_i V^*_{i'} W_j W^*_{j'} = U_{j +
m_N(i-1)} U^*_{j' + m_N(i'-1)} = W_j W^*_{j'} V_i V^*_{i'}$ for $1 \le
i,i' \le M$, $1 \le j,j' \le m_N$, and this decomposition of the matrix
units $U_{l} U^*_{l'}$ implements $\pi_{M, m_N}$. Hence $\phi_{N+1} \circ
(\pi_{M, m_N} \otimes \id_{C^*(\Lambda_{N+1})})$ satisfies
\begin{equation}\label{eq:whatamess}
\begin{split}
\phi_{N+1} \circ (\pi_{M, m_N} \otimes \id_{C^*(\Lambda_{N+1})})
&\Big(\big(\big(b_{i,j,j',i'}\big)^{m_N}_{j,j'=1}\big)^M_{i,i'=1}\Big) \\
&\textstyle= \sum^M_{i,i'=1} \sum^{m_N}_{j,j'=1} U_{j + m_N(i-1)}
b_{i,j,j',i'} U^*_{j'+m_N(i'-1)}.
\end{split}
\end{equation}

The Cuntz-Krieger relations also show that $V_i = \sum^{m_N}_{j=1} W_j
W^*_j V_i$ for all $i$, and hence $V_i a V^*_{i'} = \sum_j U_{j +
m_N(i-1)} W^*_j a W_j U^*_{j + m_N(i'-1)}$ for all $a \in P_1
C^*(\LAMBDA_N) P_1$. One now checks that for $\lambda \in \Lambda_N$, we
have
\[\textstyle
W^*_j s_\lambda W_j = \sum_{p_N(\lambda') = \lambda}
s^*_{e(r(\lambda'),j)} s_{e(r(\lambda), \cocycle_N(\lambda') j)}
s_{\lambda'},
\]
and hence that $V_i s_\lambda V^*_{i'} = \sum_j \sum_{p_N(\lambda') =
\lambda} U_{\cocycle_N(\lambda')j + m_N(i-1)} s_{\lambda'} U^*_{j +
m_N(i'-1)}$. Recall that $\theta_{i,i'}s_\lambda \in M_M(C^*(\Lambda_N))$
denotes the matrix
\[
\big(\theta_{i,i'}s_\lambda)_{j,j'} =
\begin{cases}
s_\lambda &\text{ if $j=i$ and $j' = i'$} \\
0 &\text{ otherwise.}
\end{cases}
\]
Then $V_i s_\lambda V^*_{i'} = \phi_N\big(\theta_{i,i'}s_\lambda)$ by
definition of $\phi_N$, so
\[\textstyle
\phi_N(\theta_{i,i'}s_\lambda) = \sum_j \sum_{p_N(\lambda') = \lambda}
U_{\cocycle_N(\lambda')j + m_N(i-1)} s_{\lambda'} U^*_{j + m_N(i'-1)}.
\]

Since $(\id_M \otimes \iota_{p_N, \cocycle_N})(\theta_{i,i'}s_\lambda) =
\theta_{i,i'} \sum_{p_N(\lambda') = \lambda} s_{\lambda'}$, we may
therefore apply~\eqref{eq:whatamess} to see that
\[
\phi_N(\theta_{i,i'}s_\lambda) = \phi_{N+1} \circ (\pi_{M,m_N} \otimes
\id_{C^*(\Lambda_{N+1})}) \circ (\id_M \otimes \iota_{p_N, \cocycle_N})
(\theta_{i,i'}s_\lambda).
\]
Since elements of the form $\theta_{i,i'}s_\lambda$ generate
$M_M(C^*(\Lambda_N))$ this proves the result.
\end{proof}

\begin{thm}\label{thm:direct limit}
Let $(\Lambda_n, \Lambda_{n+1}, p_n, m_n, \cocycle_n)^\infty_{n=1}$ be a
sequence of row-finite coverings of locally convex $k$-graphs. For each
$n$, let $\LAMBDA_n := \Lambda_1\cs{p_1,\cocycle_1} \cdots
\cs{p_{n-1},\cocycle_{n-1}} \Lambda_n$, identify $\LAMBDA_n$ with the
corresponding subset of $\tgrphlim(\Lambda_n; p_n, \cocycle_n)$, and
likewise identify $C^*(\LAMBDA_n)$ with the corresponding $C^*$-subalgebra
of $C^*(\tgrphlim(\Lambda_n; p_n, \cocycle_n))$. Then
\begin{equation}\label{eq:dir lim expr}
\textstyle C^*(\tgrphlim(\Lambda_n; p_n, \cocycle_n)) =
\overline{\bigcup^\infty_{n=1} C^*(\LAMBDA_n)}.
\end{equation}
Let $P_1 := \sum_{v \in \Lambda_1^0} s_v$, and for each $n$, let $M_n :=
m_1 m_2 \cdots m_{n-1}$. Then $P_1$ is a full projection in each
$\Mm(C^*(\LAMBDA_n))$, and we have
\begin{equation}\label{eq:corner dir lim}
P_1 C^*(\tgrphlim(\Lambda_n; p_n, \cocycle_n)) P_1 \cong
\varinjlim\big(M_{M_n}(C^*(\Lambda_n)), \id_{M_n} \otimes
\iota_{p_n,\cocycle_n}\big).
\end{equation}
In particular,
\begin{align*}
K_*(C^*(\tgrphlim(\Lambda_n; p_n, \cocycle_n)))
    &= K_*(P_1 C^*(\tgrphlim(\Lambda_n; p_n, \cocycle_n)) P_1) \\
    &\cong \varinjlim (K_*(C^*(\Lambda_n)),
(\iota_{p_{n},\cocycle_{n}})_*).
\end{align*}
\end{thm}
\begin{proof}
For the duration of the proof, let $\LAMBDA := \tgrphlim(\Lambda_n; p_n,
\cocycle_n)$. We have $C^*(\LAMBDA) = \clsp\{s_\mu s^*_\nu : \mu,\nu \in
\LAMBDA\}$, so for the first statement, we need only show that
\[\textstyle
\lsp\{s_\mu s^*_\nu : \mu,\nu \in \LAMBDA\} \subset \bigcup^\infty_{n=1}
C^*(\LAMBDA_n).
\]
To see this we simply note that for any finite $F \subset \LAMBDA$, the
integer $N := \max\{n \in \NN : s(F) \cap \Lambda_n^0 \not= \emptyset\}$
satisfies $F \subset \LAMBDA_N$.

Since $P_1$ is full in each $C^*(\LAMBDA_n)$ by
Proposition~\ref{prp:embedding}(3), it is full in $C^*(\LAMBDA)$
by~\eqref{eq:dir lim expr}. Equation~\ref{eq:corner dir lim} follows from
Proposition~\ref{prp:tower matrix inclusion}. The final statement then
follows from continuity of the $K$-functor.
\end{proof}

\begin{rmk}\label{rmk:restricted-gauge}
 Note that if we let $\gamma$ denote the restriction of the gauge
action to $P_1 C^*(\tgrphlim(\Lambda_n; p_n, \cocycle_n)) P_1$ then
$\gamma_{(1, \cdots, 1, z)}$ is trivial for all $z \in \TT$. Indeed, if
$s_\mu s_\nu^*$ is a nonzero element $P_1 C^*(\tgrphlim(\Lambda_n; p_n,
\cocycle_n)) P_1$, then $d(\mu)_{n+1} = d(\nu)_{n+1}$.  So $\gamma$ may be
regarded as an action by $\TT^k$ rather than $\TT^{k+1}$.
\end{rmk}

We can extend Theorem~\ref{thm:direct limit} to the situation
of matrices of covering systems as discussed in
Section~\ref{subsec:matrix of cs} as follows.

\begin{prop}\label{prp:disconnected covering algebra}
Resume the notation of Corollary~\ref{cor:pasting}. Each
$C^*(\Lambda_n)$ is canonically isomorphic to $\bigoplus_{j =
1}^{c_n} C^*(\Lambda_{n,j})$. There are homomorphisms
$(\iota_n)_* : K_*(C^*(\Lambda_n)) \to K_*(C^*(\Lambda_{n+1}))$
such that the partial homomorphism which maps the $j^{\rm th}$
summand in $K_*(C^*(\Lambda_n))$ to the $i^{\rm th}$ summand in
$K_*(C^*(\Lambda_{n+1}))$ is equal to $0$ if $m^n_{i,j} = 0$,
and is equal to $(\iota_{p^n_{i,j}, \cocycle^n_{i,j}})_*$
otherwise. The sum $\sum_{v \in \Lambda_1^0} s_v$ converges
strictly to a full projection $P_1 \in \Mm(C^*(\LAMBDA))$.
Furthermore,
\[
K_*(P_1 C^*(\LAMBDA) P_1) \cong \varinjlim\Big(\bigoplus^{c_n}_{j=1}
K_*(C^*(\Lambda_{n,j})), (\iota_n)_*\Big).
\]
\end{prop}
\begin{proof}
For each $\lambda \in \Lambda_n = \bigsqcup^{c_n}_{j=1} \Lambda_{n,j}$,
define a partial isometry $t_\lambda \in \bigoplus^{c_n}_{j=1}
C^*(\Lambda_{n,j})$ by $t_\lambda := (0,\dots,0,s_\lambda,0,\dots,0)$ (the
nonzero term is in the $j^{\rm th}$ coordinate when $\lambda \in
\Lambda_{n,j}$). These nonzero partial isometries form a Cuntz-Krieger
$\Lambda_n$-family consisting of nonzero partial isometries. The universal
property of $C^*(\Lambda_n)$ gives a homomorphism $\pi^n_t :
C^*(\Lambda_n) \to \bigoplus^{c_n}_{j=1} C^*(\Lambda_{n,j})$ which
intertwines the direct sum of the gauge actions on the
$C^*(\Lambda_{n,j})$ and the gauge action on $C^*(\Lambda_n)$. The
gauge-invariant uniqueness theorem \cite[Theorem~3.4]{kp}, and the
observation that each generator of each summand in $\bigoplus^{c_n}_{j=1}
C^*(\Lambda_{n,j})$ is nonzero and belongs to the image of $\pi^n_t$
therefore shows that $\pi^n_t$ is an isomorphism.

The individual covering systems $(\Lambda_{n,j},
\Lambda_{n+1,i}, p^n, m^n, \cocycle^n)$ induce inclusions
$\iota_{p^n_{i,j},\cocycle^n_{i,j}} : C^*(\Lambda_{n,j}) \to
M_{m^n_{i,j}}(C^*(\Lambda_{n+1,i}))$ as in
Proposition~\ref{prp:embedding}(5). We therefore obtain
homomorphisms $(\iota_{p^n_{i,j},\cocycle^n_{i,j}})_* :
K_*(C^*(\Lambda_{n,j})) \to K_*(C^*(\Lambda_{n+1,i}))$. The
statement about the partial homomorphisms of $K$-groups then
follows from the properties of the isomorphism $K_*(\bigoplus_i
A_i) \cong \bigoplus_i K_*(A_i)$ for $C^*$-algebras $A_i$.

The final statement can then be deduced from arguments similar to those of
Theorem~\ref{thm:direct limit}.
\end{proof}

\section{Simplicity and pure infiniteness}\label{sec:simplicity}
Theorem~3.1 of \cite{RobSi} gives a necessary and sufficient condition for
simplicity of the $C^*$-algebra of a row-finite $k$-graph with no sources.
Specifically, $C^*(\Lambda)$ is simple if and only if $\Lambda$ is cofinal
and every vertex of $\Lambda$ receives an aperiodic infinite path (see
below for the definitions of cofinality and aperiodicity).  In this
section we present some means of deciding whether $\tgrphlim(\Lambda_n;
p_n, \cocycle_n)$ is cofinal (Lemma~\ref{lem:cofinal tower}), and whether
an infinite path in $\tgrphlim(\Lambda_n; p_n, \cocycle_n)$ is aperiodic
(Lemma~4.3). We also present a condition under which
$C^*(\tgrphlim(\Lambda_n; p_n, \cocycle_n))$ is purely infinite
(Proposition~\ref{prp:pi tower}).

We begin by recalling the notation and definitions required to make sense
of the hypotheses of \cite[Theorem~3.1]{RobSi}. For more detail, see
Section~2 of \cite{rsy}.

\begin{ntn}\label{ntn:inf paths}
We write $\Omega_k$ for the $k$-graph such that $\Omega_k^q := \{(m, n)
\in \NN^k \times \NN^k : n-m = q\}$ for each $q \in \NN^k$, with $r(m,n) =
(m,m)$ and $s(m,n) = (n,n)$. We identify $\Omega_k^0 = \{(m,m) : m \in
\NN^k\}$ with $\NN^k$. An infinite path in a $k$-graph $\Xi$ is a graph
morphism $x : \Omega_k \to \Xi$, and we denote the image $x(0)$ of the
vertex $0 \in \Omega_k^0$ by $r(x)$. We write $\Xi^\infty$ for the
collection of all infinite paths in $\Xi$, and for $v \in \Xi^0$ we denote
by $v\Xi^\infty$ the collection $\{x \in \Xi^\infty : r(x) = v\}$. For $x
\in \Xi^\infty$ and $q \in \NN^k$, there is a unique infinite path
$\sigma^q(x) \in \Xi^\infty$ such that $\sigma^q(x)(m,n) = x(m+q, n+q)$
for all $m\le n \in \NN^k$.
\end{ntn}

\begin{dfn}\label{dfn:aperiodic+cofinal}
We say that a row-finite $k$-graph $\Xi$ with no sources is
\emph{aperiodic} if for each vertex $v \in \Xi^0$ there is an infinite
path $x \in v\Xi^\infty$ such that $\sigma^q(x) \not= \sigma^{q'}(x)$ for
all $q \not= q' \in \NN^k$. We say that $\Xi$ is \emph{cofinal} if for
each $v \in \Xi^0$ and $x \in \Xi^\infty$ there exists $m \in \NN^k$ such
that $v \Xi x(m) \not= \emptyset$.
\end{dfn}

We continue to make use in the following of the notation
established earlier (see Notation~\ref{ntn:N^k embeddings}) for
the embeddings of $\NN^k$ and of $\NN$ in $\NN^{k+1}$.

If $y$ is an infinite path in the $(k+1)$-graph $\Xi$, we write
$\alpha_y$ for the infinite path in $\Xi^{(0_k,\NN)}$ defined
by $\alpha_y(p,q) := y((0_k, p),(0_k,q))$ for $p \le q \in
\NN$, and we write $x_y$ for the infinite path in
$\Xi^{(\NN^k,0_1)}$ defined by $x_y(p,q) := y((p,0_1),
(q,0_1))$ where $p \le q \in \NN^k$.

\begin{prop}\label{prp:aperiodic tower}
Let $(\Lambda_n, \Lambda_{n+1}, p_n, m_n, \cocycle_n)^\infty_{n=1}$ be a
sequence of row-finite covering systems of $k$-graphs with no sources. For
$a,b \in \NN^{k+1}$, an infinite path $y \in \big(\tgrphlim(\Lambda_n;
p_n,\cocycle_n)\big)^\infty$ satisfies $\sigma^a(y) = \sigma^b(y)$ if and
only if $x_{\sigma^a(y)} = x_{\sigma^b(y)}$ and $\alpha_{\sigma^a(y)} =
\alpha_{\sigma^b(y)}$.
\end{prop}
\begin{proof}
The ``only if'' implication is trivial. For the ``if''
implication, note that the factorisation property implies that
an infinite path $z$ of $\tgrphlim(\Lambda_n; p_n,\cocycle_n)$
is uniquely determined by $x_z$ and the paths
$\alpha_{\sigma^{(n,0_1)}(z)}$, $n \in \NN^k$. So it suffices
to show that each $\alpha_{\sigma^{(n,0_1)}(z)}$ is uniquely
determined by $x_z(0_k,n)$ and $\alpha_z$. Fix $n \in \NN^k$
and let $\lambda := x_z(0_k,n) = z(0_{k+1},(n,0_1))$. Fix $i
\in \NN$. We will show that
$\alpha_{\sigma^{(n,0_1)}(z)}(0_1,i)$ is uniquely determined by
$\alpha_z(0_1,i)$ and $\lambda$. Let $v = r(z)$, and let $N \in
\NN$ be the element such that $v \in \Lambda^0_N$. For $1 \le j
\le i$, let $w_j = \alpha_z(i) \in \Lambda^0_{N+j}$, and let $1
\le l_j \le m_{N+j-1}$ be the integer such that
$\alpha_z(j-1,j) = e(w_i, l_j)$. We have $p_N(w_1) = v$, and
$p_{N+j-1}(w_j) = w_{j-1}$ for $2 \le j \le i$. For each $j$,
let $\lambda_j$ be the unique lift of $\lambda$ such that
$r(\lambda_j) = w_j$. By definition of the $(k+1)$-graph
$\tgrphlim(\Lambda_n; p_n,\cocycle_n)$, the path
\[
\lambda e(s(\lambda_1),\cocycle(\lambda_1)^{-1}l_1)
e(s(\lambda_2),\cocycle(\lambda_2)^{-1}l_2) \dots e(s(\lambda_i),
\cocycle(\lambda_i)^{-1} l_i) = \alpha_z(0_1,i) \lambda_i
\]
is the unique minimal common extension of $\lambda$ and
$\alpha_z(0_1,i)$ in $\tgrphlim(\Lambda_n; p_n,\cocycle_n)$.
Hence
\[
\alpha_{\sigma^{(n,0_1)}(z)}(0_1,i)
=e(s(\lambda_1),\cocycle(\lambda_1)^{-1}l_1)
e(s(\lambda_2),\cocycle(\lambda_2)^{-1}l_2) \dots e(s(\lambda_i),
\cocycle(\lambda_i)^{-1} l_i)
\]
which is uniquely determined by $\lambda$ and
$\alpha_z(0_1,i)$.
\end{proof}

\begin{cor}\label{cor:aperiodic components}
Let $(\Lambda_n, \Lambda_{n+1}, p_n, m_n, \cocycle_n)^\infty_{n=1}$ be a
sequence of row-finite covering systems of $k$-graphs with no sources.
Suppose that $\Lambda_n$ is aperiodic for some $n$. Then so is $\tgrphlim
(\Lambda_n; p_n, \cocycle_n)$.
\end{cor}
\begin{proof}
Since each vertex in $\Lambda_n$ receives an aperiodic path in
$\Lambda_n$, Proposition~\ref{prp:aperiodic tower}, guarantees that each
vertex in $\Lambda_n$ receives an aperiodic path in $\tgrphlim (\Lambda_n;
p_n, \cocycle_n)$. Since the $p_n$ are coverings, it follows that every
vertex of $\tgrphlim(\Lambda_n; p_n, \cocycle_n)$ receives an infinite
path of the form $\lambda y$ or of the form $\sigma^p(y)$ where $y$ is an
aperiodic path with range in $\Lambda_n$. If $y$ is aperiodic, then
$\lambda y$ is aperiodic for any $\lambda$ and $\sigma^a(y)$ is aperiodic
for any $a$ and the result follows.
\end{proof}

\begin{lem}\label{lem:LPF analogue}
Let $(\Lambda_n, \Lambda_{n+1}, p_n, m_n, \cocycle_n)^\infty_{n=1}$ be a
sequence of row-finite covering systems of $k$-graphs with no sources. Fix
$y \in \big(\tgrphlim(\Lambda_n; p_n, \cocycle_n)\big)^\infty$, with $y(0)
\in \Lambda_n$ and $a,b \in \NN^{k+1}$. Let $\tilde a$ and $\tilde b$
denote the elements of $\NN^k$ determined by the first $k$ coordinates of
$a$ and $b$. For each $m \ge n$, let $v_m$ and $i_m$ be the unique pair
such that $\alpha_y(m, m+1) = e(v_m, i_m)$. For each $m \ge n$, let
$\mu_m$ and $\nu_m$ be the unique lifts of $x_y(0,\tilde a)$ and
$x_y(0,\tilde b)$ such that $r(\mu_m) = r(\nu_m) = v_m$. Then
$\alpha_{\sigma^a(y)} = \alpha_{\sigma^b(y)}$ if and only if the following
three conditions hold:
\begin{Enumerate}
\item
$a_{k+1} = b_{k+1}$;
\item
$s(\mu_{m}) = s(\nu_{m})$ for all $m \ge n$; and
\item
$\cocycle_m(\mu_m)i_{m} = \cocycle_m(\nu_m)i_m$ for all $m \ge n$.
\end{Enumerate}
\end{lem}
\begin{proof}
We have $\alpha_{\sigma^a(y)}(m, m+1) = e(s(\mu_{m + a_{k+1}}),
\cocycle_m(\mu_{m + a_{k+1}})i_{m + a_{k+1}})$ for all $m$, and likewise
for $b$ and $\nu$.
\end{proof}

\begin{rmk}
Lemma~5.4 of~\cite{PRRS} implies that an infinite path in a rank-2
Bratteli diagram $\Lambda$ is aperiodic if and only if the factorisation
permutations of its red coordinate-paths are of unbounded order.
Lemma~\ref{lem:LPF analogue} is the analogue of this result for general
systems of coverings. To see the analogy, note that in a rank-2 Bratteli
diagram, every $x_y$ is of the form $\lambda\lambda\lambda\dots$ for some
blue cycle $\Lambda$, so that condition~(3) fails for all $a \not= b$
precisely when the orders of the permutations $\cocycle_m(\mu_m)$ grow
arbitrarily large with $m$.
\end{rmk}

\begin{lem}\label{lem:cofinal tower}
Let $(\Lambda_n, \Lambda_{n+1}, p_n, m_n, \cocycle_n)^\infty_{n=1}$ be a
sequence of row-finite coverings of $k$-graphs with no sources. If
infinitely many of the $\Lambda_n$ are cofinal, then $\tgrphlim(\Lambda_n;
p_n,\cocycle_n)$ is also cofinal.
\end{lem}
\begin{proof}
Fix a vertex $v$ and an infinite path $z \in
(\tgrphlim(\Lambda_n; p_n,\cocycle_n))^\infty$. Let $n_1, n_2
\in \NN$ be the elements such that $v \in \Lambda^0_{n_1}$ and
$r(z) \in \Lambda^0_{n_2}$. Choose $N \ge n_1, n_2$ such that
$\Lambda_N$ is cofinal. Fix $w \in \Lambda^0_N$ such that $p_n
\circ p_{n+1} \circ \dots \circ p_{N-1}(w) = v$; so $v
(\tgrphlim(\Lambda_n; p_n,\cocycle_n)) w \not= \emptyset$. We
have $x_{\sigma^{(0_k, N-n_2)}(z)} \in \Lambda^\infty_N$, and
since $\Lambda_N$ is cofinal, it follows that $w \Lambda_N
x_{\sigma^{(0_k, N-n_2)}(z)}(q) \not= \emptyset$ for some $q
\in \NN^k$. Since $x_{\sigma^{(0_k, N-n_2)}(z)}(q) = z(q,
N-n_2)$, this completes the proof.
\end{proof}

As in \cite{Si2}, we say that a path $\lambda$ in a $k$-graph $\Lambda$ is
a \emph{cycle with an entrance} if $s(\lambda) = r(\lambda)$, and there
exists $\mu \in r(\lambda)\Lambda$ with $d(\mu) \le d(\lambda)$ and
$\lambda(0, d(\mu)) \not = \mu$.

\begin{prop}\label{prp:pi tower}
Let $(\Lambda_n, \Lambda_{n+1}, p_n, m_n, \cocycle_n)^\infty_{n=1}$ be a
sequence of row-finite coverings of $k$-graphs with no sources. There
exists $n$ such that $\Lambda_n$ contains a cycle with an entrance if and
only if every $\Lambda_n$ contains a cycle with an entrance. Moreover, if
$C^*(\tgrphlim(\Lambda_n;p_n,\cocycle_n))$ is simple and $\Lambda_1$
contains a cycle with an entrance, then $C^*(\mathbf{\Lambda})$ is purely
infinite.
\end{prop}
\begin{proof}
That the presence of a cycle with an entrance in $\Lambda_1$ is equivalent
to the presence of a cycle with an entrance in every $\Lambda_n$ is a
consequence of the properties of covering maps. Now the result follows
from \cite[Proposition~8.8]{Si2}
\end{proof}

\section{K-Theory}\label{sec:K-th}
In this section, we consider the $K$-theory of $C^*(\Lambda
\cs{p,\cocycle} \Gamma)$. Specifically, we show how the
homomorphism from $K_*(C^*(\Lambda))$ to $K_*(C^*(\Gamma))$
obtained from Proposition~\ref{prp:embedding} behaves with
respect to existing calculations of $K$-theory for various
classes of higher-rank graph $C^*$-algebras. We will use these
results later to compute the $K$-theory of
$C^*(\tgrphlim(\Lambda_n; p_n, \cocycle_n))$ for a number of
sequences of covering systems.

Throughout this section, given a $k$-graph $\Lambda$, we view the ring
$\ZZ\Lambda^0$ as the collection of finitely supported functions $f :
\Lambda^0 \to \ZZ$. For $v \in \Lambda^0$, we denote the point-mass at $v$
by $\delta_v$. Given a finite covering $p : \Gamma \to \Lambda$ of
row-finite $k$-graphs, we define $p^* : \ZZ \Lambda^0 \to \ZZ \Gamma^0$ by
$p^*(\delta_v) = \sum_{p(u)=v} \delta_w$; equivalently, $p^*(f)(w) =
f(p(w))$.

\subsection
{Coverings of $1$-graphs and the Pimsner-Voiculescu exact sequence}
\label{sec:k=1K-theory} It is shown in \cite{pr, RSz} how to compute the
$K$-theory of a graph $C^*$-algebra using the Pimsner-Voiculescu exact
sequence. In this subsection, we show how this calculation interacts with
the inclusion of $C^*$-algebras arising from a covering of $1$-graphs.

The $K$-theory computations for arbitrary graph $C^*$-algebras \cite{DT2,
BHRS} are somewhat more complicated than for the $C^*$-algebras of
row-finite graphs with no sources. Moreover, every graph $C^*$-algebra is
Morita equivalent to the $C^*$-algebra of a row-finite graph with no
sources \cite{DT2}. We therefore restrict out attention here to the
simpler setting.

\begin{thm} \label{thm:k=1K-theory}
Let $(E^*, F^*, p,m,\cocycle)$ be a row-finite covering system of
$1$-graphs with no sources. Let $A,B$ be the vertex connectivity matrices
of the underlying graphs $E$ and $F$ respectively. Then the diagram
\begin{equation}\label{eq:big mumma}
\parbox[c]{0.9\textwidth}{\hfill
\beginpicture
    \setcoordinatesystem units <2em,2.5em>
    \setplotarea x from 1 to 16, y from -0.3 to 2.5
    \put{$0$} at 1 2
    \put{$K_1 ( C^* ( E^* ) )$} at 4 2
    \put{$\ZZ E^0$} at 7 2
    \put{$\ZZ E^0$} at 10 2
    \put{$K_0 ( C^* ( E^* ) )$} at 13 2
    \put{$0$} at 16 2
    \put{$0$} at 16 0
    \put{$K_0 ( C^* ( F^* ) )$} at 13 0
    \put{$\ZZ F^0$} at 10 0
    \put{$\ZZ F^0$} at 7 0
    \put{$K_1 ( C^* ( F^* ) )$} at 4 0
    \put{$0$} at 1 0
    \put{$\scriptstyle 1-A^t$} at 8.5 2.2
    \put{$\scriptstyle 1-B^t$} at 8.5 0.2
    \put{$\scriptstyle ( \iota_{p,\cocycle} )_*$}[l] at 4.2 1
    \put{$\scriptstyle m \cdot p^*$}[l] at 7.2 1
    \put{$\scriptstyle m \cdot p^*$}[l] at 10.2 1
    \put{$\scriptstyle ( \iota_{p,\cocycle} )_*$}[l] at 13.2 1
    \arrow <0.25cm> [0.1,0.3] from 1.5 2 to 2.6 2
    \arrow <0.25cm> [0.1,0.3] from 5.3 2 to 6.4 2
    \arrow <0.25cm> [0.1,0.3] from 7.6 2 to 9.4 2
    \arrow <0.25cm> [0.1,0.3] from 10.5 2 to 11.6 2
    \arrow <0.25cm> [0.1,0.3] from 14.4 2 to 15.5 2
    \arrow <0.25cm> [0.1,0.3] from 1.5 0 to 2.6 0
    \arrow <0.25cm> [0.1,0.3] from 5.3 0 to 6.4 0
    \arrow <0.25cm> [0.1,0.3] from 7.6 0 to 9.4 0
    \arrow <0.25cm> [0.1,0.3] from 10.5 0 to 11.6 0
    \arrow <0.25cm> [0.1,0.3] from 14.4 0 to 15.5 0
    \arrow <0.25cm> [0.1,0.3] from 4 1.5 to 4 0.5
    \arrow <0.25cm> [0.1,0.3] from 7 1.5 to 7 0.5
    \arrow <0.25cm> [0.1,0.3] from 10 1.5 to 10 0.5
    \arrow <0.25cm> [0.1,0.3] from 13 1.5 to 13 0.5
\endpicture
\hfill\hbox to 0pt{}}
\end{equation}
commutes and the rows are exact.
\end{thm}

The proof of this theorem occupies the remainder of
Section~\ref{sec:k=1K-theory}. We fix, for the duration, a finite covering
$p : F^* \to E^*$ of row-finite $1$-graphs with no sources, a multiplicity
$m$ and a cocycle $\cocycle : F^* \to S_m$.

It is relatively straightforward to prove that the right-hand two squares
of~\eqref{eq:big mumma} commute and that the rows are exact.

\begin{lem}\label{lem:right 2}
Resume the notation of Theorem~\ref{thm:k=1K-theory}. We have $(1 - B^t)
p^* = p^* (1 - A^t)$, the right-hand two squares of~\eqref{eq:big mumma}
commute, and the rows are exact.
\end{lem}
\begin{proof}
For the first statement, consider a generator $\delta_v \in \ZZ
E^0$. We have
\[
\big(p^* \circ (1 - A^t)\big)(\delta_v)
 = p^*(\delta_v - \sum_{e \in vE^1} \delta_{s(e)})
 = \sum_{p(u) = v} \delta_u - \sum_{e \in vE^1} \sum_{p(f) =
e} \delta_{s(f)}.
\]
On the other hand,
\[
\big((1 - B^t) \circ p^*\big)(\delta_v)
 = (1 - B^t) \sum_{p(u) = v} \delta_u
 = \sum_{p(u) = v} \Big(\delta_u - \sum_{f \in uF^1}
\delta_{s(f)}\Big).
\]
Since $p$ is a covering the double-sums occurring in these two equations
each contain exactly one term for each edge $f \in F^1$ such that $p(r(f))
= v$, and it follows that the two are equal.

Multiplying by $m$ throughout the above calculation shows that the middle
square of~\eqref{eq:big mumma} commutes.

The identification of $K_0(C^*(E^*))$ with $\coker(1 - A^t)$ takes the
class of the projection $s_v \in C^*(E^*)$ to the class of the
corresponding generator $\delta_v \in \ZZ E^0$ (see \cite{r}). That the
right-hand square commutes then follows from
Proposition~\ref{prp:embedding}(6).

Exactness of the rows is precisely the computation of $K$-theory for
$1$-graph $C^*$-algebras \cite{c, pr, RSz}.
\end{proof}

It remains to prove that the left-hand square of~\eqref{eq:big mumma}
commutes. The strategy is to assemble the eight-term commuting diagrams
which describe the $K$-theory of each of $C^*(E^*)$ and $C^*(F^*)$ (see
equation~\eqref{eq:graph alg K-th} below) into a sixteen-term diagram, one
face of which is the left-hand square of~\eqref{eq:big mumma}. We then
focus on the cube in the sixteen-term diagram which contains left-hand
square of~\eqref{eq:big mumma} as one of its faces, and show that the
remaining five faces of this cube commute. A diagram-chase then
establishes that the sixth face commutes as well. The majority of the work
involved goes into defining the connecting maps needed to write down the
sixteen-term diagram in the first place. The proof that the various
squares in it commute is then relatively straightforward.

To begin, we recall how one shows that the rows
of~\eqref{eq:big mumma} are exact. Let $E^* \times_d
\mathbf{Z}$ be the skew-product of $E^*$ by the length functor
$d$ (see \cite[Section~5]{kp}). Let $\gamma$ be the gauge
action of $\TT$ on $C^*(E^*)$ satisfying $\gamma_z(s_e) = zs_e$
for $e \in E^1$ and $z \in \TT$. Let $(i_{\TT}, i_{C^*(E^*)})$
be the universal covariant representation of $(C^*(E^*), \TT,
\gamma)$ in the crossed product $C^*(E^*) \times_\gamma \TT$.
By \cite[Lemma~3.1]{RSz}, there is an isomorphism
\begin{equation}\label{eq:skew-prod iso}
\psi_E : C^*(E^* \times_d \ZZ) \to C^*(E^*) \times_\gamma \TT
\quad\text{satisfying}\quad \psi_E(s_{(\lambda,n)}) = i_\TT(z)^n
i_{C^*(E^*)}(s_\lambda).
\end{equation}

The $C^*$-algebra $C^*(E^* \times_d \ZZ)$ is AF with $K_0$-group
$\varinjlim(\ZZ E^0, A^t)$ (see \cite{pr, RSz}). Hence one may apply the
dual Pimsner-Voiculescu sequence \cite[Section~10.6]{Bla} to the crossed
product algebra $C^*(E^*) \times_\gamma \TT$ to show that the top row
of~\eqref{eq:big mumma} is exact (the bottom row is the same after
replacing $E$ with $F$).

From the point of view of coverings, the skew-product graph $E^* \times_d
\ZZ$ and its $C^*$-algebra are more natural to work with than the crossed
product $C^*(E^*) \times_\gamma \TT$. Before proving that the final square
of~\eqref{eq:big mumma} commutes, we therefore detail first how coverings
$p : F^* \to E^*$ interact with the isomorphisms $\psi_E : C^*(E^*
\times_d \ZZ) \to C^*(E^*) \times_\gamma \TT$.

\begin{lem} \label{lem:skew covering}
With the above notation, let $E^* \times_d \ZZ$ and $F^* \times_d \ZZ$ be
the skew-product graphs by the length functors $d$, and let $\psi_E$ and
$\psi_F$ be the isomorphisms described in~\eqref{eq:skew-prod iso}. Let
$\gamma_E$ and $\gamma_F$ denote the gauge actions of $\TT$ on $C^*(E^*)$
and $C^*(F^*)$.
\begin{Enumerate}
\item
the formulae $\tilde{p}(\lambda, n) := (p(\lambda), n)$ and
$\tilde{\cocycle}(\lambda,n) := \cocycle(\lambda)$ determine a covering
$\tilde{p} : F^* \times_d \ZZ  \to E^* \times_d \ZZ$ and a cocycle
$\tilde{\cocycle} : F^* \times_d \ZZ \to S_m$.
\item
the inclusion $\iota_{p,\cocycle} : C^*(E^*) \to M_m(C^*(F^*))$ is
equivariant in $\gamma_E$ and $\id_m \otimes \gamma_F$, and induces an
inclusion $\widetilde{\iota_{p,\cocycle}} : C^*(E^*) \times_{\gamma_E} \TT
\to M_m(C^*(F^*)) \times_{\id_m \otimes \gamma_F} \TT$.
\item
The following diagram commutes.
\[
\beginpicture
    \setcoordinatesystem units <2.5em, 2.5em>
    \put{$C^*(E^* \times_d \ZZ)$}[B] at 0 1
    \put{$M_m(C^*(F^* \times_d \ZZ))$}[B] at 4.4 1
    \put{$C^*(E^*) \times_{\gamma_E} \TT$}[t] at 0 0
    \put{$M_m(C^*(F^*)) \times_{\id_m \otimes \gamma_F} \TT$}[t] at
4.4 0
    \arrow <0.6em> [0.1,0.3] from 0 0.8 to 0 0.2
    \put{$\scriptstyle \psi_E$}[r] at -0.1 0.5
    \arrow <0.6em> [0.1,0.3] from 4.4 0.8 to 4.4 0.2
    \put{$\scriptstyle \id_m \otimes \psi_F$}[r] at 4.3 0.5
    \arrow <0.6em> [0.1,0.3] from 1.2 -0.2 to 2.3 -0.2
    \arrow <0.6em> [0.1,0.3] from 1.2 1.15 to 2.8 1.15
    \put{$\scriptstyle \iota_{\tilde p, \tilde\cocycle}$}[b] at 2
1.25
    \put{$\scriptstyle \widetilde{\iota_{p,\cocycle}}$}[b] at 1.85
-0.1
\endpicture
\]
\end{Enumerate}
\end{lem}
\begin{proof}
(1) It is straightforward to check that $\tilde p$ is a covering. To see
that $\tilde\cocycle$ is a cocycle, note that $(\mu,m)$ and $(\nu,n)$ are
composable in the skew-product precisely when $\mu$ and $\nu$ are
composable, and $n = m - d(\nu)$. So for $i \in \{1, \dots, m\}$ we may
calculate
\[
\tilde\cocycle(\mu,m)(\tilde\cocycle(\nu,m - d(\nu))i)
 = \cocycle(\mu)(\cocycle(\nu)i) = \cocycle(\mu\nu)i
 = \tilde\cocycle(\mu\nu,m-d(\nu))i.
\]

(2) That $\iota_{p,\cocycle}$ is equivariant in $\gamma_E$ and $\id_m
\otimes \gamma_F$ follows from Proposition~\ref{prp:Lambda cs Gamma}(4).
That it induces the desired inclusion $\widetilde{\iota_{p,\cocycle}}$ of
crossed-products follows from the universal properties of the
crossed-product algebras.

(3) That the diagram commutes follows from a simple calculation using the
definitions of the maps involved.
\end{proof}

\begin{proof}[Proof of Theorem~\ref{thm:k=1K-theory}]
Lemma~\ref{lem:right 2} establishes everything except that the left-hand
square in the diagram~\eqref{eq:big mumma} commutes. To establish this
last claim, recall from \cite[Theorem~3.2]{RSz} (see also \cite{pr}) that
there is a homomorphism $\phi_E : \ZZ E^0 \to K_0(C^*(E^*)
\times_{\gamma_E} \TT)$ satisfying $\phi_E(\delta_v) =
[i_\TT(1)i_{C^*(E^*)}(s_v)]$. Moreover, the rows of the following
commutative diagram are exact and the left- and right-most vertical maps
are isomorphisms (see \cite[Lemma 7.15]{r}, \cite{pr}).
\begin{equation}\label{eq:graph alg K-th}
\parbox[c]{0.9\textwidth}{\hfill
\beginpicture
    \small \setcoordinatesystem units <2.1em,3em>
    \put{$0$}[B] at -2.4 1
    \arrow <0.6em> [0.15,0.45] from -2.2 1.1 to -1.75 1.1
    \put{$\ker(1 - A^t)$}[B] at -.5 1
    \arrow <0.6em> [0.15,0.45] from 0.7 1.1 to 2.5 1.1
    \put{$\ZZ E^0$}[B] at 3 1
    \arrow <0.6em> [0.15,0.45] from 3.5 1.1 to 7.5 1.1
    \put{$\scriptstyle 1 - A^t$}[b] at 5.5 1.2
    \put{$\ZZ E^0$}[B] at 8 1
    \arrow <0.6em> [0.15,0.45] from 8.5 1.1 to 10.0 1.1
    \put{$\coker(1 - A^t)$}[B] at 11.6 1
    \arrow <0.6em> [0.15,0.45] from 13.1 1.1 to 13.5 1.1
    \put{$0$}[B] at 13.7 1
    %
    \put{$0$}[B] at -2.4 -0.25
    \arrow <0.6em> [0.15,0.45] from -2.2 -0.15 to -1.8 -0.15
    \put{$K_1(C^*(E^*\!))$}[B] at -.5 -0.25
    \arrow <0.6em> [0.15,0.45] from 0.75 -0.15 to 1.1 -0.15
    \put{$K_0(C^*(E^*\!)\times_{\gamma_E}\!\TT)$}[B] at 3 -0.25
    \arrow <0.6em> [0.15,0.45] from 5 -0.15 to 6 -0.15
    \put{$\scriptstyle 1 - \hat{\gamma}^{-1}_*$}[b] at 5.5 -0.1
    \put{$K_0(C^*(E^*\!)\times_{\gamma_E}\!\TT))$}[B] at 8 -0.25
    \arrow <0.6em> [0.15,0.45] from 10 -0.15 to 10.35 -0.15
    \put{$K_0(C^*(E^*\!))$}[B] at 11.6 -0.25
    \arrow <0.6em> [0.15,0.45] from 12.85 -0.15 to 13.5 -0.15
    \put{$0$}[B] at 13.7 -0.25
    %
    \arrow <0.6em> [0.1,0.3] from -.5 0.9 to -.5 0.1
    \put{$\scriptstyle \cong$}[l] at -0.3 0.5
    \arrow <0.6em> [0.15,0.45] from 3 0.9 to 3 0.1
    \put{$\scriptstyle \phi_E$}[l] at 3.2 0.5
    \arrow <0.6em> [0.15,0.45] from 8 0.9 to 8 0.1
    \put{$\scriptstyle \phi_E$}[l] at 8.2 0.5
    \arrow <0.6em> [0.15,0.45] from 11.7 0.9 to 11.7 0.1
    \put{$\scriptstyle \cong$}[l] at 11.9 0.5
    \endpicture
\hfill\hbox to 0pt{}}
\end{equation}

A similar commutative diagram holds for ${F^*}$, and using the standard
isomorphism of $K_*(M_m(C^*(F^*)))$ with $K_*(C^*(F^*))$, we may assemble
these two diagrams can into a single three-dimensional diagram by
connecting each term in the diagram for $E^*$ to the corresponding term in
the diagram for $F^*$ using the appropriate maps induced from $(p,
\cocycle)$. The map connecting the $K_0$-groups of the skew-product graph
algebras is induced from the connecting map in the bottom row of the
commuting diagram in Lemma~\ref{lem:skew covering}(3) by applying the
$K$-functor and using the canonical isomorphisms
\begin{gather*}
K_*(M_m(C^*(F^*) \times_{\gamma_F} \TT)) \cong K_*(C^*(F^*)
\times_{\gamma_F}
\TT) \quad\text{and} \\
M_m(C^*(F^*) \times_{\gamma_F} \TT) \cong M_m(C^*(F^*)) \times_{\id_m
\otimes \gamma_F} \TT.
\end{gather*}

Let $\eta$ denote the unlabelled inclusion $K_1(C^*(F^*)) \hookrightarrow
K_0(C^*(F^* \times_d \ZZ))$ in the bottom row of the diagram of the
form~\eqref{eq:graph alg K-th} for ${F^*}$. Notice that injectivity of the
map $m \cdot p^* : \ZZ E^0 \to \ZZ F^0$ together with the first statement
of Lemma~\ref{lem:right 2} ensures that $m \cdot p^*$ restricts to a map
from $\ker(1 - A^t)$ to $\ker(1 - B^t)$; abusing notation, we denote this
map $m \cdot p^*$ too. With this notation the diagram~\eqref{eq:thecube}
below is the left-hand cube of the three-dimensional diagram described
in the previous paragraph.
\begin{equation}\label{eq:thecube}
\parbox[c]{0.9\textwidth}{\hfill
    \beginpicture
    \setcoordinatesystem units <2.5em, 2.5em>
    \put{$\ker(1 - A^t)$} at 0 2
    \arrow <0.6em> [0.1, 0.3] from -0.5 1.75 to -2.5 0.75
    \put{$\scriptstyle m\cdot p^*$} [br] at -1.4 1.35
    \put{$\ker(1 - B^t)$} at -3 0.5
    \put{$K_1(C^*(E^*))$} at 0 -1
    \arrow <0.6em> [0.1, 0.3] from -0.5 -1.25 to -2.5 -2.25
    \put{$\scriptstyle(\iota_{p,\cocycle})_*$} [br] at -1.5 -1.65
    \put{$K_1(C^*({F^*}))$} at -3 -2.5
    \arrow <0.6em> [0.1, 0.3] from 0 1.8 to 0 -0.7
    \put{$\scriptstyle\cong$} [bl] at 0.1 0.8
    \arrow <0.6em> [0.1, 0.3] from -3 0.3 to -3 -2.2
    \put{$\scriptstyle\cong$} [bl] at -2.9 -0.7
    \arrow <0.6em> [0.1, 0.3] from 1.2 1.95 to 4.5 1.95
    \circulararc -180 degrees from 1.2 1.95 center at 1.2 2
    \arrow <0.6em> [0.1, 0.3] from 1.1 -1 to 3.6 -1
    \arrow <0.6em> [0.1, 0.3] from -1.8 0.45 to 1.5 0.45
    \circulararc -180 degrees from -1.8 0.45 center at -1.8 0.5
    \arrow <0.6em> [0.1, 0.3] from -1.9 -2.5 to 0.6 -2.5
    \put{$\scriptstyle\eta$} [b] at -0.6 -2.4
    \put{$\ZZ E^0$} at 5 2
    \arrow <0.6em> [0.1, 0.3] from 4.5 1.75 to 2.5 0.75
    \put{$\scriptstyle m \cdot p^*$} [br] at 3.6 1.35
    \put{$\ZZ F^0$} at 2 0.5
    \put{$K_0(C^*(E^*)\times_\gamma\TT)$} at 5.2 -1
    \arrow <0.6em> [0.1, 0.3] from 4.5 -1.25 to 2.5 -2.25
    \put{$\scriptstyle (\widetilde{\iota_{p, \cocycle}})_*$} [br] at
3.6 -1.65
    \put{$K_0(C^*(F^*)\times_\gamma\TT)$} at 2.2 -2.5
    \arrow <0.6em> [0.1, 0.3] from 5 1.8 to 5 -0.7
    \put{$\scriptstyle\phi_E$} [bl] at 5.1 0.8
    \arrow<0.6em> [0.1, 0.3] from 2 0.3 to 2 -2.2
    \put{$\scriptstyle\phi_E$} [bl] at 2.1 -0.7
    \endpicture
\hfill\hbox to 0pt{}}
\end{equation}
We have shown the whole cube because we prove that the left-hand face ---
which is none other than the left-hand square of~\eqref{eq:big mumma} ---
commutes by showing that the other five faces commute.

To see why this suffices, suppose that the other five faces do indeed
commute. Since $\eta$ is an injection by the exactness of the rows
of~\eqref{eq:graph alg K-th}, we just need to show that the two maps from
$\ker(1 - A^t)$ into $K_0(C^*({F^*}) \times_\gamma \TT))$ obtained from
the maps in the left-hand face of the cube followed by $\eta$ agree. A
diagram chase shows that this is the case.

It therefore remains only to show that the top, bottom, front, back and
right-hand faces of~\eqref{eq:thecube} commute. The top square commutes by
definition. The bottom square commutes by the naturality of the dual
Pimsner-Voiculescu exact sequence (see the argument at the beginning of
\cite[Section~3]{RSz}). The back and front faces commute
because~\eqref{eq:graph alg K-th} commutes.

To see that the right-hand face commutes, recall that $C^*(E^* \times_d
\ZZ)$ is AF with $K_0$-group $\varinjlim(\ZZ E^0, 1-A^t)$. Hence there is
an inclusion $\varepsilon_E : \ZZ E^0 \to K_0(C^*(E^* \times_d \ZZ))$
which takes $\delta_v$ to the $K_0$-class of the vertex projection
$s_{(v,0)}$, and likewise for $F$. Consider the map $\psi_E$ defined
in~\eqref{eq:skew-prod iso} and the map $\phi_E$ appearing
in~\eqref{eq:graph alg K-th}. It is clear that $\phi_E = (\psi_E)_* \circ
\varepsilon_E$ and similarly for $F$. So it suffices to show that the
following diagram commutes.
\begin{equation}\label{eq:two rectangles}
\parbox[c]{0.9\textwidth}{\hfill
\beginpicture
    \setcoordinatesystem units <2.5em, 2.5em>
    \put{$\ZZ E^0$}[B] at 0 2.3
    \put{$\ZZ F^0$}[B] at 5.4 2.3
    \put{$K_0(C^*(E^* \times_d \ZZ))$}[B] at 0 1
    \put{$K_0(C^*(F^* \times_d \ZZ))$}[B] at 5.4 1
    \put{$K_0(C^*(E^*) \times_{\gamma_E} \TT)$}[t] at 0 0
    \put{$K_0(C^*(F^*) \times_{\gamma_F} \TT)$}[t] at 5.4 0
    \arrow <0.6em> [0.1,0.3] from 0 2.2 to 0 1.4
    \put{$\scriptstyle \varepsilon_E$}[r] at -0.1 1.8
    \arrow <0.6em> [0.1,0.3] from 5.4 2.2 to 5.4 1.4
    \put{$\scriptstyle \varepsilon_F$}[r] at 5.3 1.8
    \arrow <0.6em> [0.1,0.3] from 0 0.9 to 0 0.1
    \put{$\scriptstyle (\psi_E)_*$}[r] at -0.1 0.5
    \arrow <0.6em> [0.1,0.3] from 5.4 0.9 to 5.4 0.2
    \put{$\scriptstyle (\psi_F)_*$}[r] at 5.3 0.5
    \arrow <0.6em> [0.1,0.3] from 0.5 2.45 to 4.9 2.45
    \put{$\scriptstyle m\cdot p^*$}[b] at 2.6 2.55
    \arrow <0.6em> [0.1,0.3] from 1.7 1.15 to 3.8 1.15
    \put{$\scriptstyle (\iota_{\tilde p, \tilde\cocycle})_*$}[b] at
2.6 1.25
    \arrow <0.6em> [0.1,0.3] from 1.7 -0.2 to 3.8 -0.2
    \put{$\scriptstyle (\widetilde{\iota_{p,\cocycle}})_*$}[b] at
2.6 -0.1
\endpicture
\hfill\hbox to 0pt{}}
\end{equation}
If one applies the $K$-functor to all terms and maps in the diagram of
Lemma~\ref{lem:skew covering}(3), and then applies the natural isomorphism
\[
K_*(M_m(C^*(E^*) \times_{\gamma_E} \TT)) \cong K_*(C^*(E^*)
\times_{\gamma_E} \TT)
\]
to the terms on the right, one obtains precisely the bottom rectangle
of~\eqref{eq:two rectangles}. The bottom rectangle of~\eqref{eq:two
rectangles} therefore commutes by naturality of the $K$-functor together
with Lemma~\ref{lem:skew covering}(3).

To see that the top rectangle of~\eqref{eq:two rectangles} commutes,
recall that $\varepsilon_E$ takes the image of the point-mass $\delta_v$
in the direct limit $\varinjlim(\ZZ E^0, A^t)$ to the class of the
projection $s_{(v,0)}$. The image of $s_{(v,0)}$ under the homomorphism
$\iota_{\tilde p, \tilde\cocycle}$ is the diagonal matrix in $M_m(C^*(F^*
\times_d \ZZ))$ whose diagonal entries are all equal to $\sum_{p(w) = v}
s_{(w,0)}$. Under the standard isomorphism of $K_0(M_m(C^*(F^* \times_d
\ZZ)))$ with $K_0(C^*(F^* \times_d \ZZ))$, we therefore obtain the
following equality in $K_0(C^*(F^* \times_d \ZZ))$:
\[
[\iota_{\tilde p, \tilde\cocycle}(s_{(v,0)})]
    = \sum_{p(w) = v} m\cdot [s_{(w,0)}]
    = m \cdot \Big(\sum_{p(w) = v} [s_{(w,0)}]\Big).
\]
Using once again the characterisation of the maps $\varepsilon_E$ and
$\varepsilon_F$, we see that this is precisely the statement that the
bottom rectangle of~\eqref{eq:two rectangles} commutes.
\end{proof}

\subsection
{Coverings of higher-rank graphs and Kasparov's spectral sequence theorem}
\label{sec:k>1K-theory}

We turn to the case where $k > 1$.  We invoke the $K$-theory
computations of \cite{e} which are based on Kasparov's spectral
sequence theorem for the computation of the $K$-theory of
crossed products by groups for which the Baum-Connes conjecture
holds (see \cite[Theorem~6.10]{k}, \cite[Lemma~3.3]{e} and
\cite{rst}). We are grateful to Gennadi Kasparov for pointing
out that the spectral sequence is natural.

The standard notation for spectral sequences is that a spectral sequence
$(E^r, d^r)$ has terms $E^r_{p,q}$ and differentials $d^r : E^r_{p,q} \to
E^r_{p-r, q+r-1}$ where $r > 0$ and $p,q \in \ZZ$. This however is
problematic in the current situation because $p$ clashes with our notation
for a covering map. To avoid this, we replace the indexing variables $p,q$
in the spectral sequence with $a,b$. That is, our spectral sequences have
terms $E^r_{a,b}$ and differentials $d^r : E^r_{a,b} \to E^r_{a-r, b+r-1}$
where $r > 0$ and $a,b \in \ZZ$.

Since each higher rank graph $C^*$-algebra $C^*(\Lambda)$ is Morita
equivalent to a crossed product by $\ZZ^k$ \cite[Theorem~5.6]{KP3},
Kasparov's result applies to give a spectral sequence which converges to
$K_*(C^*(\Lambda))$ with $E^2$ terms given by the homology of $\ZZ^k$ with
appropriately chosen coefficients. In \cite{e} Evans computes these
homology groups using a resolution related to the Koszul complex. It
follows that the above spectral sequence may be extended so that the terms
of the resolution become the terms $E^1_{a,b}$ for $b$ even.

The main result of this subsection is to show that given a
finite covering $p : \Gamma \to \Lambda$ of row-finite
$k$-graphs with no sources, a multiplicity $m$ and a cocycle
$\cocycle : \Gamma \to S_m$, there is a natural morphism of
spectral sequences defined on $E^1$ terms using $m \cdot p^* :
\ZZ \Lambda^0 \to \ZZ \Gamma^0$ which is compatible (see
\cite[p.\,126]{w}) with $(\iota_{p,\cocycle})_*$ the induced
map on $K$-theory. This result is specialised to the case $k =
2$ with a view to applications in Section~\ref{sec:examples}.

The following is an immediate Corollary of \cite[Theorem
6.10]{k} (see \cite[Lemma~3.3]{e} and \cite{rst}). For more
detail on spectral sequences used in this context, see
\cite{rst,e}.

\begin{prop}\label{prp:Kasparov}
Let $\F$ be a C*-algebra and let $\alpha : \ZZ^k \to \text{Aut\,}\F$ be an
action of $\ZZ^k$ on $\F$.  Then there is a spectral sequence $(E^r, d^r)$
with differentials $d^r : E^r_{a,b} \to E^r_{a-r, b+r-1}$ which converges
to $K_*(\F \times_\alpha \ZZ^k)$ with $E^2_{a,b} = H_a(\ZZ^k, K_b(\F))$.
Moreover, the spectral sequence is natural with respect to equivariant
maps of C*-algebras.
\end{prop}
\begin{proof}
As noted in the proof of \cite[Lemma~3.3]{e} this follows
immediately from \cite[Theorem 6.10]{k} since $\ZZ^k$ is
amenable and the Baum-Connes conjecture is known to hold for
amenable groups \cite[Theorem~1.1]{HK}, so the $\gamma$ part of
$K_*(\F \times_\alpha \ZZ^k)$ exhausts.  The naturality of the
spectral sequence with respect to equivariant maps follows from
the construction in the proof of \cite[Theorem 6.10]{k}, since
every step is functorial.
\end{proof}

Naturality means that given $\ZZ^k$ actions $\alpha_i$ on $\F_i$, a
$\ZZ^k$-equivariant map $\varphi : \F_1 \to \F_2$ induces a morphism of
spectral sequences and this morphism is compatible with
\[
\widehat\varphi_* : K_*(\F_1 \times_{\alpha_1} \ZZ^k) \to K_*(\F_2
\times_{\alpha_2} \ZZ^k)
\]
where $\widehat\varphi : \F_1 \times_{\alpha_1} \ZZ^k \to \F_2
\times_{\alpha_2} \ZZ^k$ is the natural map.

Evans applies this when $\F = \F_\Lambda$ is the crossed
product $C^*(\Lambda) \times_\gamma \TT^k$ of $C^*(\Lambda)$ by
the gauge action, and $\alpha$ is the dual action
$\hat{\gamma}$ of $\ZZ^k$. Hence, by Takai duality we have
$K_*(C^*(\Lambda)) = K_*(\F_\Lambda \times_\alpha \ZZ^k)$. In
this case we have more specific results (see
\cite[Lemma~3.3]{e})
$$
   E^2_{a,b} =
   \begin{cases}
     H_a(\ZZ^k, K_0(\F_\Lambda)) & \text { if } 0 \le a \le k
        \text{ and $b$ is even,}\\
     0 & \text{ otherwise.}
   \end{cases}
$$
In \cite[Theorem~3.14]{e}), Evans shows that these homology
groups may be computed as the homology of the complex
$D^\Lambda_* = {\bigwedge}^{*} \ZZ^k \otimes \ZZ \Lambda^0$.
That is, $D^\Lambda_a = {\bigwedge}^{a} \ZZ^k \otimes
\ZZ\Lambda^0$ for $0 \le a \le k$ and $D^\Lambda_a = 0$ for $a
> k$. For $1 \le j \le k$ let $M_j$ denote the vertex
connectivity matrix of the coordinate graph $(\Lambda^0,
\Lambda^{e_j}, r, s)$. For $1 \le a \le k$ define the
differential $\partial_a : D^\Lambda_a \to D^\Lambda_{a-1}$ by
$$
\partial_a(\e_{i_1} \wedge \cdots \wedge \e_{i_a} \otimes e_v) =
\sum_{j = 1}^a (-1)^{j+1} \e_{i_1}\wedge \cdots \wedge
\widehat{\e}_{i_j}\wedge\cdots\wedge \e_{i_a} \otimes (1 - M_j^t)e_v
$$
where $\e_1, \dots, \e_k$ constitute the canonical basis for
$\ZZ^k$, $1 \le {i_1} < \cdots < {i_a} \le k$ and $v \in
\Lambda^0$. It is straightforward to verify that $D^\Lambda_*$
is a complex. The first part of the following theorem is a
restatement of \cite[Theorem~3.15]{e}).

\begin{thm} \label{thm:k>1K-theory}
Fix $k > 1$.  Let $\Lambda$ be a row-finite $k$-graph with no sources.
With notation as above there is a spectral sequence $(E^r, d^r)$ with
differentials $d^r : E^r_{a,b} \to E^r_{a-r, b+r-1}$ which converges to
$K_*(C^*(\Lambda)) = K_*(\F_{\Lambda} \times_\alpha \ZZ^k)$ with
\[
E^1_{a,b} = D^\Lambda_a := {\bigwedge}^{a} \ZZ^k \otimes \ZZ\Lambda^0,
\]
if $0 \le a \le k$ and $b$ is even, and $0$ otherwise.  The differential
$d^1 : E^1_{a,b} \to E^1_{a-1,b}$ is given by $\partial_a$ if $b$ is even.

Let $(\Lambda,\Gamma,p,m,\cocycle)$ be a row-finite covering system of
$k$-graphs with no sources. There is a morphism $f$ of spectral sequences
compatible with $(\iota_{p, \cocycle})_* : K_*(C^*(\Lambda)) \to
K_*(C^*(\Gamma))$ such that $f^1 : D^\Lambda_a \to D^\Gamma_a$ is given by
$\id \otimes (m \cdot p^*)$.
\end{thm}
\begin{proof}
Evans computes the homology groups using a Koszul complex (see \cite[\S
4.5]{w}).  Set $G = \ZZ^k = \langle s_1, \dots s_k \rangle$, $R = \ZZ G$
and let $I$ be the ideal in $R$ generated by $\{ 1 - s_a^{-1} : 1 \le a
\le k \}$. Let $\e_1, \dots, \e_k$ constitute the canonical basis for
$R^k$. For each $a$, define $\partial_a : {\bigwedge}^a R^k \to
{\bigwedge}^{a - 1} R^k$ as follows: for $1 \le {i_1} < \cdots < {i_a} \le
k$ so that $\e_{i_1} \wedge \cdots \wedge \e_{i_a} \in {\bigwedge}^a R^k$,
define
\[
\partial_a(\e_{i_1} \wedge \cdots \wedge \e_{i_a}) =
\sum_{j = 1}^a (-1)^{j+1} (1 - {s_j}^{-1}) \e_{i_1}\wedge \cdots \wedge
\widehat{\e}_{i_j}\wedge\cdots\wedge \e_{i_a}
\]
where the symbol `` $\widehat{\cdot}$ '' denotes deletion of an element
(note that $\partial_1(\e_{i}) = 1 - s_i^{-1}$).

Then $R/I \cong \ZZ$ and the following sequence of $R$-modules is exact
(see \cite[Corollary 4.5.5]{w})
\[
0 \to {\bigwedge}^k R^k \to \cdots \to {\bigwedge}^1 R^k \to {\bigwedge}^0
R^k \to \ZZ  \to 0.
\]
Note that ${\bigwedge}^0 R^k = R$ and ${\bigwedge}^a R^k$ is a free
$R$-module with basis
\[
\{ \e_{i_1} \wedge \cdots \wedge \e_{i_a} : 1 \le {i_1} < \cdots < {i_a}
\le k \}.
\]
Hence, ${\bigwedge}^* R^k$ yields a projective resolution of $\ZZ$. Thus,
by \cite[\S III.1]{b} we have
\[
H_*(G, K_0(\F_\Lambda)) \cong H_*({\bigwedge}^* R^k\otimes_G
K_0(\F_\Lambda)).
\]
We follow Evans here but have adopted slightly different
notation to make naturality more apparent (see
\cite[Definition~3.11]{e} and following). Under the isomorphism
${\bigwedge}^a R^k\otimes_G K_0(\F_\Lambda) \cong {\bigwedge}^a
\ZZ^k\otimes K_0(\F_\Lambda)$ (as abelian groups), the boundary
map $\partial_a : {\bigwedge}^a \ZZ^k \otimes K_0(\F_\Lambda)
\to {\bigwedge}^{a - 1} \ZZ^k \otimes K_0(\F_\Lambda)$ is given
by
\[
\partial_a(\e_{i_1} \wedge \cdots \wedge \e_{i_a} \otimes x) =
\sum_{j = 1}^a (-1)^{a+1} \e_{i_1}\wedge \cdots \wedge
\widehat{\e}_{i_j}\wedge\cdots\wedge \e_{i_a} \otimes (1 - s_j)x
\]
where $1 \le {i_1} < \cdots < {i_a} \le k$ and $x \in K_0(\F_\Lambda)$.

Let $D^\Lambda_a$ be given as above. There is a natural map
$\varepsilon^\Lambda : C_0 (\Lambda^0) \hookrightarrow
\F_\Lambda$ which induces a map $\varepsilon^\Lambda_* : \ZZ
\Lambda^0 \to K_0(\F_\Lambda)$. Moreover (see
\cite[Theorem~3.14]{e}) the natural map
\[
\id \otimes\varepsilon^\Lambda_* : {\bigwedge}^{*} \ZZ^k \otimes \ZZ
\Lambda^0
 \to {\bigwedge}^{*} \ZZ^k \otimes K_0(\F_\Lambda)
\]
is a map of complexes which induces an isomorphism on homology and hence
\[
H_*(G, K_0(\F_\Lambda)) \cong H_*( {\bigwedge}^* \ZZ^k\otimes \ZZ
\Lambda^0).
\]
Therefore, setting
\[
E^1_{a,b} =
\begin{cases}
  \bigwedge^a \ZZ^k\otimes \ZZ \Lambda^0
& \text{if $0 \le a \le k$ and $b$ is even,}\\
  0 & \text{otherwise}
\end{cases}
\]
and defining $d^1 : E^1_{a,b} \to E^1_{a-1,b}$ to be $\partial_a$ if $b$
is even (and $0$ otherwise), yields
\[
E^2_{a,b} \cong
\begin{cases}
H_p(G, K_0(\F_\Lambda))
& \text{if $0 \le a \le k$ and $b$ is even,}\\
  0 & \text{otherwise}.
\end{cases}
\]
It follows by \cite[Lemma~3.3]{e} that the spectral sequence
converges to $K_*(C^*(\Lambda)) = K_*(\F_{\Lambda}
\times_\alpha \ZZ^k)$ as required.

For the second part of the theorem, fix $(\Lambda, \Gamma, p, m,
\cocycle)$. The embedding $\iota_{p,\cocycle} : C^*(\Lambda) \to
M_m(C^*(\Gamma))$ induces an embedding $\widetilde{\iota_{p, \cocycle}} :
\F_\Lambda \to M_m(\F_\Gamma)$. Functoriality yields a map of complexes
\[
\id \otimes (\widetilde{\iota_{p,\cocycle}})_* : {\bigwedge}^{*} \ZZ^k
\otimes K_0(\F_\Lambda)
 \to {\bigwedge}^{*} \ZZ^k \otimes K_0(\F_\Gamma).
\]
Since group homology is a covariant functor of its coefficient module we
obtain the functorial maps for each $n = 0, 1, \dots, k$
\[
H_n((\widetilde{\iota_{p,\cocycle}})_*) : H_n(\ZZ^k, K_0(\F_\Lambda)) \to
H_n(\ZZ^k, K_0(\F_\Gamma)).
\]

Then arguing as in Lemma~\ref{lem:right 2} with $p^* : \ZZ \Lambda^0 \to
\ZZ\Gamma^0$ defined as above we see that
\[
(1- (M^\Gamma_j)^t)(m \cdot p^*) = (m \cdot p^*)(1- (M^\Lambda_j)^t)
\]
for all $j = 1, \dots, k$.  It follows that the natural map
\[
\id \otimes (m \cdot p^*) : {\bigwedge}^{*} \ZZ^k \otimes \ZZ \Lambda^0
\to {\bigwedge}^{*} \ZZ^k \otimes \ZZ \Gamma^0
\]
is a map of complexes.

Arguing as in the proof of Theorem~\ref{thm:k=1K-theory}, we see that
$(\widetilde{\iota_{p,\cocycle}})_* \circ \varepsilon^\Lambda_* =
\varepsilon^\Gamma_* \circ (m \cdot p^*)$, so the map on homology induced
by $\id \otimes (m \cdot p^*)$ coincides with the functorial map above
(under the identifications of the homology groups induced by $\id
\otimes\varepsilon^\Lambda_*$ and $\id \otimes\varepsilon^\Gamma_*$). This
combined with the naturality of Proposition~\ref{prp:Kasparov} yields a
morphism $f$ of spectral sequences compatible with the map
\[
\big(\widehat{\widetilde{\iota_{p,\cocycle}}}\big)_*: K_*(\F_{\Lambda}
\times_\alpha \ZZ^k) \to K_*(\F_{\Gamma} \times_\alpha \ZZ^k)
\]
such that $f^1 : D^\Lambda_a \to D^\Gamma_a$ is given by $\id \otimes (m
\cdot p^*)$. Under the identifications $K_*(C^*(\Lambda)) =
K_*(\F_{\Lambda} \times_\alpha \ZZ^k)$ and $K_*(C^*(\Gamma)) =
K_*(\F_{\Gamma} \times_\alpha \ZZ^k)$, we have
$\big(\widehat{\widetilde{\iota_{p, \cocycle}}}\big)_* =
(\iota_{p,\cocycle})_*$.
\end{proof}

The following corollary is an immediate consequence of the
above theorem restricted to the case $k = 2$; for the first
assertion see \cite[Proposition~3.16]{e} and its proof (see
also \cite{rst}).

Given a $2$-graph $\Lambda$, recall that $M_1$ and $M_2$ denote the vertex
connectivity matrices of the coordinate graphs $(\Lambda^0, \Lambda^{e_1},
r, s)$ and $(\Lambda^0, \Lambda^{e_2}, r, s)$.

\begin{cor}\label{cor:k=2K-theory}
Suppose that $(\Lambda, \Gamma, p, m, \cocycle)$ is a row-finite covering
system of $2$-graphs with no sources. With the notation of
Theorem~\ref{thm:k>1K-theory}, the complex $D^\Lambda_a = {\bigwedge}^{a}
\ZZ^2 \otimes \ZZ \Lambda^0$ may be written as follows:
\begin{equation}\label{eq:complex}
0 \leftarrow \ZZ \Lambda^0 \xleftarrow{\partial_1} \ZZ \Lambda^0 \oplus
\ZZ \Lambda^0 \xleftarrow{\partial_2} \ZZ \Lambda^0 \leftarrow 0
\end{equation}
where $\partial_1 = (1 - M_1^t, 1 - M_2^t)$ and
 $\partial_2 =
\begin{pmatrix}
  M_2^t - 1 \\ 1 - M_1^t
\end{pmatrix}.$
We have $E^2_{a,b} = E^\infty_{a,b}$, and
\begin{equation}\label{eq:K-th maps}
\begin{split}
 K_0(C^*(\Lambda)) &\cong \coker \partial_1 \oplus \ker \partial_2
\\
 K_1(C^*(\Lambda)) &\cong \ker \partial_1/ \Im \partial_2 \cong
                     H_1(\ZZ^k, K_0(\F_\Lambda)).
\end{split}
\end{equation}

Moreover, the following diagram commutes
\begin{equation}\label{eq:gwion CD}
\begin{CD}
  0 @<<< \ZZ \Lambda^0 @<{\partial^\Lambda_1}<<
\ZZ \Lambda^0 \oplus \ZZ \Lambda^0 @<{\partial^\Lambda_2}<<
\ZZ \Lambda^0 @<<< 0\\
@.  @VVm \cdot p^*V  @VV m \cdot p^* \oplus m \cdot p^*V  @VV m
\cdot p^*V @. \\
  0 @<<< \ZZ \Gamma^0 @<{\partial^\Gamma_1}<<
\ZZ \Gamma^0 \oplus \ZZ \Gamma^0 @<{\partial^\Gamma_2}<<
\ZZ \Gamma^0 @<<< 0\\
\end{CD}
\end{equation}
and by naturality induces $(\iota_{p,\cocycle})_* : K_*(C^*(\Lambda)) \to
K_*(C^*(\Gamma))$.
\end{cor}

The inclusion of $\coker\partial_1$ into $K_0(C^*(\Lambda))$
obtained from~\eqref{eq:K-th maps} takes the equivalence class
(in the quotient group $\coker\partial_1 =
\ZZ\Lambda^0/\Im(\partial_1)$) of the generator $\delta_v$ of
$\ZZ \Lambda^0$ to the $K_0$-class of the vertex projection
$[s_v]$ in $C^*(\Lambda)$. The proof of this fact can be
obtained from the proof of~\cite[Proposition~4.4]{e}. We thank
Gwion Evans for pointing this out to us.

\subsection{Product coverings and the K\"unneth
formula}\label{sec:Kunneth}

In this section we consider covering systems $(\Lambda_n, p_n)$ in which
each $k$-graph $\Lambda_n$ is a cartesian product of two lower-dimensional
graphs, and the covering maps $p_n$ respect the product decomposition.

Recall from \cite[Proposition~1.8]{kp} that given a $k$-graph
$(\Lambda,d)$ and a $k'$-graph $(\Lambda',d')$, the
cartesian-product category $\Lambda \times \Lambda'$ becomes a
$(k + k')$-graph when endowed with the degree functor $d \times
d' : (\lambda,\lambda') \mapsto (d(\lambda)_1, \dots,
d(\lambda)_k, d'(\lambda')_1, \dots, d'(\lambda')_{k'})$.

\begin{prop} \label{prp:kunneth}
Fix $k,k' \in \NN\setminus\{0\}$. Let $(\Lambda,\Gamma,p,m,\cocycle)$ and
$(\Lambda',\Gamma',p',m',\cocycle')$ be row-finite covering systems of
$k$- and $k'$-graphs with no sources. Then
\[
p \times p' : \Gamma \times \Gamma' \to \Lambda \times \Lambda'
\]
is a finite covering of row-finite $(k+k')$-graphs with no
sources. Let $f : \{1, \dots, m\} \times \{1 \dots, m'\} \to
\{1, \dots, mm'\}$ denote the bijection $f(j,j') := j +
(j'-1)m$. There is a cocycle $\cocycle \times \cocycle' :
\Gamma \times \Gamma' \to S_{mm'}$ determined by
$\big((\cocycle \times \cocycle')(\alpha,\alpha')\big)f(j,j')
:= f\big(\cocycle(\alpha)j, \cocycle'(\alpha')j'\big)$.
Moreover, the following diagram commutes.
\[
\begin{CD}
  C^*(\Lambda \times \Lambda') @> \cong >> C^*(\Lambda) \otimes
    C^*(\Lambda') \\
  @VV \iota_{p \times p', \cocycle \times \cocycle'} V
    @VV \iota_{p, \cocycle} \otimes \iota_{p', \cocycle'} V \\
  M_{mm'}(C^*(\Gamma \times \Gamma')) @> \cong >>
    M_m(C^*(\Gamma)) \otimes M_{m'}(C^*(\Gamma'))
\end{CD}
\]
Suppose that at least one of $K_*(C^*(\Lambda))$, $K_*(C^*(\Lambda'))$ and
at least one of $K_*(C^*(\Gamma))$, $K_*(C^*(\Gamma'))$ are torsion-free.
Then the following diagram commutes and the horizontal connecting maps are
zero-graded isomorphisms:
\[
\begin{CD}
  K_*(C^*(\Lambda)) \otimes K_*(C^*(\Lambda'))
    @>\cong>> K_*(C^*(\Lambda \times \Lambda')) \\
  @VV (\iota_{p, \cocycle})_* \otimes (\iota_{p',\cocycle'})_* V
    @VV (\iota_{p \times p', \cocycle \times \cocycle'})_* V \\
  K_*(C^*(\Gamma)) \otimes K_*(C^*(\Gamma'))
    @>\cong>> K_*(C^*(\Gamma \times \Gamma')) \\
\end{CD}
\]
If $\Gamma^0$ and ${\Gamma'}^0$ (and hence also $\Lambda^0$ and
${\Lambda'}^0$) are finite then the $C^*$-algebras are unital, and the
horizontal isomorphisms take $[1] \otimes [1]$ to $[1]$.
\end{prop}
\begin{proof}
It is straightforward to check that $p \times p'$ is a covering using the
properties of the covering maps $p$ and $p'$ and the definition of the
cartesian-product graph. A simple calculation shows that $\cocycle \times
\cocycle'$ defines a cocycle.

Theorem~5.5 of \cite{kp} shows that $C^*(\Lambda)$,
$C^*(\Lambda')$, $C^*(\Gamma)$ and $C^*(\Gamma')$ are nuclear,
and so there is just one tensor-product $C^*$-algebra
$C^*(\Lambda) \otimes C^*(\Lambda')$. Corollary~3.5(iv) of
\cite{kp} shows that the map $s_{(\lambda,\mu)} \mapsto
s_\lambda \otimes s_\mu$ is an isomorphism of $C^*(\Lambda
\times \Lambda')$ onto $C^*(\Lambda) \otimes C^*(\Lambda')$,
and similarly for $C^*(\Gamma)$ and $C^*(\Gamma')$. It is easy
to check using the formulae for the maps $\iota_{p, \cocycle},
\iota_{p', \cocycle'}$, and $\iota_{p \times p', \cocycle
\times \cocycle'}$ and using the chain of isomorphisms
\begin{align*}
M_{mm'}(C^*(\Gamma \times \Gamma'))
 &\cong M_{mm'}(\CC) \otimes C^*(\Gamma \times \Gamma') \\
 &\cong M_m(\CC) \otimes C^*(\Gamma) \otimes M_{m'}(\CC) \otimes
C^*(\Gamma') \\
 &\cong M_m(C^*(\Gamma)) \otimes M_{m'}(C^*(\Gamma'))
\end{align*}
that the first diagram commutes.

In the presence of the additional hypothesis concerning torsion-free
$K$-groups, the K\"unneth Theorem of \cite{Sch} (see also Theorem~23.1.3
of \cite{Bla}) implies: (1) that
\[
K_*(C^*(\Lambda)) \otimes K_*(C^*(\Lambda')) \cong K_*(C^*(\Lambda)
\otimes C^*(\Lambda'))
\]
and similarly for $\Gamma, \Gamma'$; (2) that these isomorphisms are
natural and are zero-graded; and (3) that these isomorphisms take $[1]
\otimes [1]$ to $[1]$. The result therefore follows from the naturality of
the $K$-functor.
\end{proof}

Note that in general when no assumption is made about torsion, the
K\"unneth Theorem of \cite{Sch} gives a short exact sequence which is
still natural. The analogue of Proposition~\ref{prp:kunneth} still holds
and gives a (fairly complicated) commuting diagram in which the rows are
short exact sequences.

\section{Examples}\label{sec:examples}
In this section we discuss a number of examples. A recurring theme will be
supernatural numbers and the associated dimension groups, so we pause here
to establish some notation.

We will think of a supernatural number as an infinite product
$\alpha = \prod^\infty_{n=1} \alpha_n$ where each $\alpha_n$ is
an integer greater than $1$. Any two such expressions in which
the same prime factors occur with the same cardinality
correspond to the same supernatural number. Given supernatural
numbers $\alpha, \beta$, we will abuse notation and write
$\alpha\beta$ for the supernatural number
$\prod^\infty_{n=1}\alpha_n\beta_n$. We write $\alpha[1,n]$ for
the product $\prod^n_{i=1} \alpha_i$ of the first $n$ terms in
$\alpha$.

For $z_1, \dots, z_n \in \CC$, we write $\ZZ[z_1, \dots, z_n]$ for
the ring obtained by adjoining $z_1, \dots, z_n$ to $\ZZ$; we regard
$\ZZ[z_1, \dots, z_n]$ as a group under addition. Abusing notation,
for a supernatural number $\alpha$, we write
$\ZZ\big[\frac{1}{\alpha}\big]$ for the dimension group
$\varinjlim(\ZZ, \times \alpha_n)$ which we identify with the group
\[
\bigcup^\infty_{n=1} \ZZ\Big[\frac{1}{\alpha[1,n]}\Big] \subset
\QQ
\]
consisting of all fractions $p/q$ where $p,q \in \ZZ$, and $q$ is a
divisor of some $\alpha[1,n]$.

\subsection{Rank-2 Bratteli diagrams}

A rank-2 Bratteli diagram is a $2$-graph in which the blue edges form a
Bratteli diagram and the red edges determine simple cycles so that every
vertex lies on precisely one red cycle, and all vertices on a given red
cycle are at the same level in the blue Bratteli diagram.

The $C^*$-algebras of these $2$-graphs were studied in \cite{PRRS} and
provided the initial motivation for the covering construction. A rank-2
Bratteli diagram $\Lambda$ can be constructed using
Proposition~\ref{prp:multiple covering} and Corollary~\ref{cor:pasting}
precisely  when the length of each red cycle at level $n$ of $\Lambda$ is
divisible by the lengths of all the cycles at level $n-1$ to which it
connects. In particular, the $2$-graphs whose $C^*$-algebras are Morita
equivalent to the Bunce-Deddens algebras \cite[Example~6.7]{PRRS} and the
irrational rotation algebras \cite[Example~6.5]{PRRS} arise in this
fashion.

\subsection{\boldmath{Coverings of dihedral graphs $D_n$}}

\begin{figure}[ht]
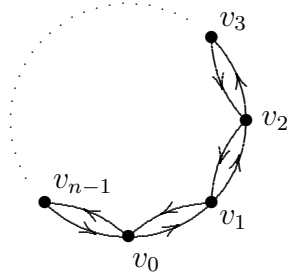

\[
    \beginpicture
    \setcoordinatesystem units <2.5em, 2.5em>
    \put{$\bullet$} at 6 0
    \put{$v_0$}[lt] at 6.05 -0.2
    \put{$\bullet$} at 7.06 0.44
    \put{$v_{1}$}[lt] at 7.16 0.34
    \put{$\bullet$} at 7.5 1.5
    \put{$v_{2}$}[l] at 7.7 1.5
    \put{$\bullet$} at 7.06 2.56
    \put{$v_{3}$}[lb] at 7.16 2.66
    \put{$\bullet$} at 4.94 0.44
    \put{$v_{n-1}$}[lb] at 5.04 0.54
    \setdots
    \circulararc 150 degrees from 6.75 2.8 center at 6 1.5
    \setsolid
    \circulararc 135 degrees from 6 0 center at 6 1.5
    \circulararc -45 degrees from 6 0 center at 6 1.5
    \setquadratic
    \plot 6 0 6.48 0.32 7.06 0.44 /
    \arrow <0.5em> [0.25,0.75] from 6.44 0.310 to 6.42 0.300
    \arrow <0.5em> [0.25,0.75] from 6.62 0.15 to 6.64 0.16
    \plot 7.06 0.44 7.18 1.02 7.5 1.5 /
    \arrow <0.5em> [0.25,0.75] from 7.20 1.05 to 7.17 0.95
    \arrow <0.5em> [0.25,0.75] from 7.385 0.95 to 7.415 1.05
    \plot 7.06 2.56 7.18 1.98 7.5 1.5 /
    \arrow <0.5em> [0.25,0.75] from 7.17 2.00 to 7.20 1.90
    \arrow <0.5em> [0.25,0.75] from 7.415 2.00 to 7.385 2.10
    \plot 6 0 5.52 0.32 4.94 0.44 /
    \arrow <0.5em> [0.25,0.75] from 5.46 0.345 to 5.44 0.355
    \arrow <0.5em> [0.25,0.75] from 5.56 0.05 to 5.58 0.04
    \endpicture
\]
\caption{The $1$-graph $D_n$} \label{fig:D_n}
\end{figure}

For $n \in \NN\setminus\{0\}$, let $D_n$ be the directed graph
with $n$ vertices $\{v_0, \dots, v_{n-1}\}$ and edges $\{x_i,
y_i : 0 \le i \le n-1\}$ where $r(x_i) = v_i = s(y_i)$ and
$s(x_i) = v_{i+1} = r(y_i)$ (throughout this section, addition
in the subscripts is understood to be evaluated modulo $n$).
More descriptively, $D_n$ is a ring of $n$ vertices, each of
which connects to both of its neighbours (see
Figure~\ref{fig:D_n}). Let $D_n^*$ be the path-category of
$D_n$, regarded as a $1$-graph.
Note that for $n \in \NN \setminus \{0\}$, the graph $D_{2n}$
is the Cayley graph for the dihedral group with $2n$ elements.

\begin{example}\label{eg:Dn}
For $n , m \ge 1$ there are $m$-fold covering maps $p_{n,mn} : D^*_{nm}
\to D^*_{n}$ as follows: for $0 \le i \le nm-1$ let $i' = i \mod n$ and
define
\[
    p_{n,mn}(v_i) := v_{i'},
    \quad p_{n,mn}(x_i) := x_{i'}
    \quad\text{and}
    \quad p_{n,mn}(y_i) := y_{i'}.
\]
Hence for each $n,m \ge 1$, we obtain a row-finite covering system
$(D^*_n, D^*_{mn}, p_{n,mn})$ of 1-graphs with no sources (see
Notation~\ref{ntn:m=1}).

Fix an infinite supernatural number $\alpha = \prod^\infty_{i=1}
\alpha_i$. Consider the sequence of covering systems $(D^*_{6
\alpha[1,n]}, D^*_{6 \alpha[1,n+1]}, p_{6 \alpha[1,n], 6
\alpha[1,n+1]})_{n=1}^\infty$ as in Notation~\ref{ntn:m=1}.

Applying Corollary~\ref{cor:tower graph}, we obtain a $2$-graph $D :=
\tgrphlim(D^*_{6 \alpha[1,n]}, p_{6 \alpha[1,n], 6 \alpha[1,n+1]})$.
\end{example}

\begin{prop}\label{prp:Dn}
Consider the situation of Example~\ref{eg:Dn}. We have $K_0(C^*(D)) =
\ZZ[\frac{1}{\alpha}] \oplus \ZZ[\frac{1}{\alpha}]$ and $K_1(C^*(D)) = \ZZ
\oplus \ZZ$. Let $P_1 := \sum_{v \in D_6^0} s_v$. Then $[P_1]$ is the 0
element of $K_0(P_1 C^*(D) P_1)$. Moreover, $C^*(D)$ is simple and purely
infinite.
\end{prop}

Before proving the proposition, we describe the $K$-theory of
$C^*(D^*_n)$ in general.

\begin{lem}\label{lem:D_n K-th}
\begin{Enumerate}
\item
$K_0(C^*(D^*_n))$ is generated by $[s_{v_0}]$ and $[s_{v_1}]$, and for
each $i$, we have $[s_{v_i}] = -[s_{v_{i+3}}]$ in $K_0(C^*(D^*_n))$.
\item
$K_1(C^*(D^*_n)) \cong \{(a_1, \dots, a_n) \in \ZZ^n : a_{i+2} = a_{i+1} -
a_i \text{ for all $i$}\}$.
\item
the following table describes the $K$-theory of each $C^*(D^*_n)$.
\[
\begin{tabular}{||c|c|c||}
\hline $\vbox to 1em{}n\ {\rm mod}\ 6$ & $K_0(C^*(D^*_n))$ &
$K_1(C^*(D^*_n))$
\\
\hline
$\vbox to 1em{} 0$ & $\ZZ^2$& $\ZZ^2$\\
$1$ & $0$ & $0$\\
$2$ & $\ZZ/3\ZZ$& $0$\\
$3$ & $\ZZ/2\ZZ \oplus \ZZ/2\ZZ$& $0$\\
$4$ & $\ZZ/3\ZZ$& $0$\\
$5$ & $0$& $0$\\
\hline
\end{tabular}
\]
\end{Enumerate}
\end{lem}
\begin{proof}
(1) The $K_0$ group is generated by the classes $[s_{v_0}], \dots,
[s_{v_{n-1}}]$ subject to the relations $[s_{v_i}] = [s_{v_{i+1}}] +
[s_{v_{i-1}}]$. This relation forces $[s_{v_{i+2}}] = [s_{v_{i+1}}] -
[s_{v_i}]$, from which we conclude first that $K_0$ is generated by
$[s_{v_0}]$ and $[s_{v_1}]$ and second that
\[
[s_{v_{i+3}}]
    = [s_{v_{i+2}}] - [s_{v_{i+1}}]
    = ([s_{v_{i+1}}] - [s_{v_i}]) - [s_{v_{i+1}}]
    = -[s_{v_i}].
\]

(2) Let $A_n$ denote the vertex connectivity matrix of $D_n$; so $A_n(i,j)
= 1$ when $i = j \pm 1\;(\textnormal{mod}\ n)$ and zero otherwise. As in
Theorem~\ref{thm:k=1K-theory}, we have $K_1(C^*(D^*_n)) \cong
\ker(1-A_n^t)$. For $m \in \ZZ^n$, $((1 - A_n^t)m)_i = - m_{i-1} + m_i -
m_{i+1}$ by definition of $A_n$, and this establishes~(2).

(3) If $E$ is a finite $1$-graph with no sinks or sources, then $C^*(E)$
is isomorphic to the Cuntz-Krieger algebra of the adjacency matrix $A_E$
of $E$ \cite{KPRR}. In particular, $K_1(C^*(E))$ is torsion-free and has
the same rank as $K_0(C^*(E))$ \cite{c2}. Hence it suffices to verify that
the first column of the table is correct. To calculate $K_0$, we use~(1)
to check by hand that the cases $n = 1, 2, \dots 6$ are as claimed. If $n
> 6$, then applying the relations we find that $[s_{v_{i+6}}] = [s_{v_i}]$
for all $i$ which accounts for all remaining cases.
\end{proof}

\begin{proof}[Proof of Proposition~\ref{prp:Dn}]
Lemma~\ref{lem:D_n K-th}(1) shows that $K_0(C^*( D^*_{6\alpha[1,n]})$ is
generated by $[s_{v^n_0}]$ and $[s_{v^n_1}]$ where the $v^n_i$ are the
vertices of $D^*_{6\alpha[1,n]}$. Fix $i \in \{1,2\}$. We have
\begin{equation}\label{eq:Dn inclusion}
(\iota_{p_n})_* [s_{v^n_i}]
 = [s_{v^{n+1}_i}]
   + [s_{v^{n+1}_{i + 6\alpha[1,n]}}]
   + \dots
   + [s_{v^{n+1}_{i + 6(\alpha_{n+1}-1)\alpha[1,n]}}].
\end{equation}
By Lemma~\ref{lem:D_n K-th}(1), each $[s_{v^{n+1}_i + 6k}] =
[s_{v^{n+1}_i}]$ in $K_0(C^*(D^*_{6\alpha[1,n+1]}))$, so~\eqref{eq:Dn
inclusion} implies $(\iota_{p_n})_* [s_{v^n_i}] =
\alpha_n\cdot[s_{v^{n+1}_i}]$. Hence $K_0(\iota_{p_n}) : \ZZ^2 \to \ZZ^2$
is multiplication by $\alpha_n$.

Fix $m \in \NN\setminus\{0\}$. By Lemma~\ref{lem:D_n K-th}(2),
$K_1(C^*(D^*_{6 m}))$ is identified with the set of sequences $(a_1,
\dots, a_{6 m})$ which satisfy $a_{i+2} = a_{i+1} - a_i$ for all $i$. By
Lemma~\ref{lem:D_n K-th}(2), this forces $a_{i+2} = a_{i+1} - a_i$ for all
$i$. Consequently, the map $a = (a_1, \dots, a_{6m}) \mapsto (a_1, a_2)$
yields an isomorphism $\zeta_m : K_1(C^*(D^*_{6m})) \to \ZZ^2$. As
$\zeta_{\alpha[1,n+1]} \circ K_1(\iota_{p_{6\alpha[1,n], 6\alpha[1,n+1]}})
= \zeta_{\alpha[1,n]}$, it follows that$K_1(D) \cong \ZZ^2$.

Recall that $D$ denotes $\tgrphlim(D^*_{6 \alpha[1,n]}, p_{6
\alpha[1,n], 6 \alpha[1,n+1]})$. By
Theorem~\ref{thm:k=1K-theory} the $K$-groups of $C^*(D)$ are as
claimed. To compute the class of the identity, let $P_1 \in
C^*(D)$ be the sum of the six vertex projections in the bottom
level. The final statement of Lemma~\ref{lem:D_n K-th}(1) shows
that the classes of the vertex projections in $K_0(C^*(D^*_6))$
cancel, so that the class of the identity in $K_0(C^*(D^*_6))$
is the zero element. It follows that the class of the identity
$P_1$ in $K_0(P_1 C^*(D) P_1)$ is also the zero element.

Each $D^*_n$ is aperiodic and cofinal (see
Definition~\ref{dfn:aperiodic+cofinal}), so we may conclude from
Corollary~\ref{cor:aperiodic components} and Lemma~\ref{lem:cofinal tower}
that $D$ is aperiodic and cofinal. Hence Proposition~4.8 of \cite{kp}
implies that $C^*(D)$ is simple. The path $x_1 y_1$ is a cycle with an
entrance (namely $y_0$) in $D^*_1$. Proposition~\ref{prp:pi tower} now
shows that $C^*(D)$ is purely infinite.
\end{proof}

\subsection{Direct limits of $\Oo_n \otimes
C(\TT)$}\label{subsec:P_n}

\begin{example}\label{eg:Pn}
Fix $n \ge 3$, and let $B_n$ be the bouquet of $n$ loops. For $m \ge 1$,
let $L_m$ denote the loop with $m$ vertices, and let $\Lambda_m$ be the
cartesian-product $2$-graph $\Lambda_m = L^*_{(n-1)^m} \times B^*_n$
obtained from the path categories of $L_{(n-1)^m}$ and $B_n$.

For each $m$, Let $p_m$ denote the obvious $(n-1)$-fold covering of
$L^*_{(n-1)^m}$ by $L^*_{(n-1)^{m+1}}$, and let $p'$ be the identity
covering of $B_n$ by $B_n$.
\end{example}

\begin{prop}
Consider the situation of Example~\ref{eg:Pn}. Let $v$ be a vertex of
$\Lambda_1$. Then $s_v C^*(\tgrphlim(\Lambda_m, p_m \times p')) s_v$ is
isomorphic to the Kirchberg algebra $\mathcal{P}_n$ (see \cite{Bla2})
whose $K$-theory is opposite to that of $\Oo_n$.
\end{prop}
\begin{proof}
Since $C^*(B_n)$ is generated by $n$ isometries whose range projections
sum to the identity, $C^*(B_n)$ is canonically isomorphic to $\Oo_n$
\cite{c1}. Hence
\[
C^*(\Lambda_m) \cong C^*(L^*_{(n-1)^m}) \otimes \Oo_n
\]
by \cite[Corollary~3.5(iv)]{kp}. As in \cite[Lemma~2.4]{aHR},
$C^*(L^*_{(n-1)^m}) \cong M_{(n-1)^m}(C(\TT))$ for all $m$, and
in particular, $K_*(C^*(L^*_{(n-1)^m})) \cong (\ZZ,\ZZ)$. Since
$K_*(\Oo_n) = (\ZZ/(n-1)\ZZ,0)$ \cite{c2}, the K\"unneth
theorem implies that $K_*(C^*(\Lambda_m)) \cong (\ZZ/(n-1)\ZZ,
\ZZ/(n-1)\ZZ)$.

A special case of \cite[Equation~(4.7)]{PRRS} implies that the covering
map $p_m$ induces multiplication by $n-1$ from $K_0(C^*(L^*_{(n-1)^m}))$
to $K_0(C^*(L^*_{(n-1)^{m+1}}))$, and the identity homomorphism from
$K_1(C^*(L^*_{(n-1)^m}))$ to $K_1(C^*(L^*_{(n-1)^{m+1}}))$. Clearly $p'$
induces the identity map on $K_*(\Oo_n)$.

Let $\Lambda = \tgrphlim(\Lambda_m, p_m \times p')$.
Theorem~\ref{thm:direct limit} and
Proposition~\ref{prp:kunneth} combine to show that
\[
K_*(C^*(\Lambda)) \cong \varinjlim((\ZZ/(n-1)\ZZ, \ZZ/(n-1)\ZZ), (\times (n-1),
\id)).
\]
Since multiplication by $n-1$ is the 0 homomorphism from
$\ZZ/(n-1)\ZZ$ to $\ZZ/(n-1)\ZZ$, it follows that
$K_*(C^*(\Lambda)) \cong (0, \ZZ/(n-1)\ZZ)$.

Lemma~\ref{lem:cofinal tower} proves that $\Lambda$ is cofinal.
For an infinite path $y \in \Lambda^\infty$, Lemma~\ref{lem:LPF
analogue} combined with the observation that the cycles in the
$L^*_{(n-1)^m}$ grow with $m$ shows that if $a,b \in \NN^3$ and
$\sigma^a(y) = \sigma^b(y)$, then $a$ and $b$ differ only in
their first coordinates. It follows from
Proposition~\ref{prp:aperiodic tower} that the aperiodicity of
$\Lambda$ is implied by the well-known aperiodicity of $B_n$.
Hence $C^*(\Lambda)$ is simple by \cite[Proposition~4.8]{kp}.
Moreover, since every vertex of $\Lambda$ hosts a cycle with an
entrance, $C^*(\Lambda)$ is also purely infinite (see
\cite[Proposition~4.9]{kp}, \cite[Proposition~8.8]{Si2}). The
result therefore follows from the Kirchberg-Phillips
classification theorem \cite{Phi}.
\end{proof}

\subsection{Higher-rank Bunce-Deddens algebras}\label{subsec:Delta_2}

In this subsection we describe a class of simple A$\TT$ algebras with
real-rank 0 which arise from sequences of covering systems of $2$-graphs
and which cannot in general be obtained from the construction of
\cite{PRRS} (see Example~\ref{eg:HRBD} and Theorem~\ref{thm:3-graph
K-th}). We indicate in Remark~\ref{rmk:HRBD} why we think of these
algebras as higher-rank analogues of the Bunce-Deddens algebras.

For $k \ge 1$, let $\Delta_k$ be the $k$-graph with vertices $\ZZ^k$,
morphisms $\{(m,n) \in \ZZ^k \times \ZZ^k : m \le n\}$ where $r(m,n) = m$,
$s(m,n) = n$ and $d(m,n) = n-m$. There is a free action of $\ZZ^k$ on
$\Delta_k$ given by translation; that is $m \cdot (p,q) = (p+m, q+m)$ for
$m \in \ZZ^k$ and $(p,q) \in \Delta_k$.

Given a finite-index subgroup $H$ of $\ZZ^k$, we denote by
$\Delta_k/H$ the quotient of $\Delta_k$ by the action of $H$.
That is, for $q \in \NN^k$, $(\Delta_k/H)^q = \{[g, g + q] : g
\in \ZZ^k \}$; in particular, $(\Delta_k/H)^0 = \{[g, g] : g
\in \ZZ^k \}$, and we henceforth identify $(\Delta_k/H)^0$ with
$\ZZ^k/H$ via the map $[g, g] \mapsto [g]$ where $[g]$ denotes
the class $g + H$ of $g$ in $\ZZ^k/H$. The range and source
maps in $\Delta_k/H$ are then given by $r([g,g+q]) = [g]$ and
$s([g,g+q]) = [g + q]$. If $H' \subset H$ is a finite-index
subgroup of $H$, then it also has finite index in $\ZZ^k$, and
there is a natural surjection $p : \ZZ^k/H' \to \ZZ^k/H$ which
induces a finite covering map, also denoted $p$ of $\Delta_k/H$
by $\Delta_k/H'$.

Most of the remainder of this section is concerned with the following
example of a sequence of covering systems.

\begin{example}\label{eg:HRBD}
Let $H_1 \supset H_2 \supset H_3 \supset \dots$ be a chain of finite-index
subgroups of $\ZZ^2$. For each $n$, let $p_n : \Delta_2/H_{n+1} \to
\Delta_2/H_n$ be the canonical covering induced by the quotient maps
described above, and let $\cocycle_n : \Delta_2/H_{n+1} \to S_1$ be the
trivial cocycle. This data specifies a sequence $(\Delta_2/H_n,
\Delta_2/H_{n+1}, p_n)^\infty_{n=1}$ of row-finite covering systems of
$2$-graphs with no sources. Applying Corollary~\ref{cor:tower graph}, we
obtain a $3$-graph $\tgrphlim(\Delta_2/H_n; p_n))$. As always, $P_1$
denotes $\sum_{v \in (\Delta_2/H_1)^0} s_v \in C^*(\Delta_2/H_1) \subset
C^*(\tgrphlim(\Delta_2/H_n; p_n))$.
\end{example}

\begin{thm}\label{thm:3-graph K-th}
Consider the situation of Example~\ref{eg:HRBD}.
\begin{Enumerate}
\item We have
\begin{equation*}\label{eq:Delta tower K_0}
 K_0(P_1 C^*(\tgrphlim(\Delta_2/H_n; p_n)) P_1) \cong \varinjlim(\ZZ, \times
[H_n : H_{n+1}])\oplus \ZZ,
\end{equation*}
and this isomorphism takes $[P_1]$ to $(g,0)$ where $g$ is the image of
$[\ZZ^2 : H_1]$ in the direct limit $\varinjlim(\ZZ, \times [H_n :
H_{n+1}])$.

\item For each $n$ the homomorphism from $\ZZ^2$ to $\ZZ^2$
    determined by coordinatewise multiplication by the integer
    $[H_n : H_{n+1}]$ restricts to a homomorphism
    $m_{H_n,H_{n+1}} : H_n \to H_{n+1}$. Moreover,
\begin{equation*}\label{eq:Delta tower K_1}
K_1(P_1 C^*(\tgrphlim(\Delta_2/H_n; p_n)) P_1) \cong \varinjlim(H_n,
m_{H_{n}, H_{n+1}}).
\end{equation*}

\item $C^*(\tgrphlim(\Delta_2/H_n; p_n))$ is simple if and only if
$\bigcap_{n = 1}^\infty H_n = \{0\}$, and is an A$\TT$ algebra with
real-rank 0 when it is simple.
\end{Enumerate}
\end{thm}

The proof of this result will occupy the bulk of this section.
Before presenting it, we state a Corollary and use it to
formulate some concrete examples.

\begin{cor}\label{cor:matrix realisation}
Consider the situation of Example~\ref{eg:HRBD}. There are
sequences $(h^n_1)^\infty_{n=1}$ and $(h^n_2)^\infty_{n=1}$ in
$\ZZ^2$ such that: (1) for each $n$, the elements $h^n_1$ and
$h^n_2$ generate $H_n$; and (2) the matrix $M_n = \left(
\begin{smallmatrix}
m^n_{1,1} & m^n_{1,2} \\
m^n_{2,1} & m^n_{2,2}
\end{smallmatrix}
\right)$ satisfying $h^{n+1}_1 = m^n_{1,1} h^n_1 + m^n_{1,2} h^n_2$ and
$h^{n+1}_2 = m^n_{2,1} h^n_1 + m^n_{2,2} h^n_2$ has positive determinant
for all $n$. Moreover, if $M^{\rm ca}_n$ denotes the classical adjoint
$\left(
\begin{smallmatrix}
m^n_{2,2} & -m^n_{1,2} \\
-m^n_{2,1} & m^n_{1,1}
\end{smallmatrix}
\right)$ of $M_n$ for each $n$, and if we regard these matrices as
homomorphisms of $\ZZ^2$, then
\begin{equation}\label{eq:alt Delta tower K_1}
K_1(P_1 C^*(\tgrphlim(\Delta_2/H_n; p_n)) P_1) \cong \varinjlim(\ZZ^2,
M_n^{\rm ca}).
\end{equation}
\end{cor}
\begin{proof}
That we can choose the $h^n_i$ so that the matrices $M_n$ all
have positive determinant follows from an inductive argument
based on the observation that replacing $h^{n+1}_i$ with
$-h^{n+1}_i$ reverses the sign of $\det(M_n)$.

For each $n$, let $\psi_n$ be the isomorphism of $\ZZ^2$ onto $H_n$
satisfying $\psi_n(e_i) = h^n_i$, and let $m_{H_n, H_{n+1}} : H_n \to
H_{n+1}$ be the homomorphism described in Theorem~\ref{thm:3-graph
K-th}(2). We claim that $\psi_{n+1} \circ M^{\rm ca}_n = m_{H_n, H_{n+1}}
\circ \psi_n$.

To see this, observe that $m_{H_n, H_{n+1}}$ is multiplication by the
determinant of $M_n$. Hence, as rational transformations,
$m_{H_n,H_{n+1}}^{-1} \circ M_n^{\rm ca} = M_n^{-1}$. Since $m_{H_n,
H_n+1}$ commutes with $\psi_{n+1}$, the desired equality $\psi_{n+1} \circ
M^{\rm ca}_n = m_{H_n, H_{n+1}} \circ \psi_n$ is therefore equivalent to
$\psi_{n+1} = \psi_n \circ M_n$, which follows from the definitions of the
maps involved. This establishes the claim.

The claim guarantees that $\varinjlim(H_n, m_{H_n, H_{n+1}}) \cong
\varinjlim(\ZZ^2, M^{\rm ca}_n)$, and~\eqref{eq:alt Delta tower K_1} then
follows from Theorem~\ref{thm:3-graph K-th}(2).
\end{proof}

\begin{examples}\label{egs:quotients}
\begin{Enumerate}
\item Let $\alpha$ and $\beta$ be supernatural numbers. For $n
    \in \NN \setminus\{0\}$, let $\phi_n$ be the homomorphism of $\ZZ^2$
    determined by the diagonal matrix $M_n :=
    \left(\begin{smallmatrix}
\alpha_n & 0 \\
0 & \beta_n
\end{smallmatrix}\right)$.

For each $n$, let
\[\textstyle
H_n := \alpha[1,n]\ZZ \times \beta[1,n]\ZZ = \phi_n(\ZZ^2) \subset \ZZ^2.
\]
We deduce from Theorem~\ref{thm:3-graph K-th} that
\[\textstyle
K_*(P_1 C^*(\tgrphlim(\Delta_2/H_n; p_n)) P_1) = \Big(
\ZZ\big[\frac{1}{\alpha\beta}\big] \oplus \ZZ,\
\ZZ\big[\frac{1}{\alpha}\big] \oplus \ZZ\big[\frac{1}{\beta}\big] \Big),
\]
that the position of the unit in $K_0$ corresponds to the element
$(\alpha_1,0)$, and that $P_1 C^*(\tgrphlim(\Delta_2/H_n; p_n)) P_1$ is a
simple A$\TT$ algebra of real-rank 0.

We claim that this is an example of an A$\TT$ algebra which cannot be
realised using a rank-2 Bratteli diagram as in \cite{PRRS}. To see this,
suppose otherwise. Then \cite[Theorem~6.1]{PRRS} implies that there exists
an injective homomorphism $\phi : \ZZ\big[\frac{1}{\alpha}\big] \oplus
\ZZ\big[\frac{1}{\beta}\big] \to \ZZ\big[\frac{1}{\alpha\beta}\big] \oplus
\ZZ$ such that each element of $\coker(\phi)$ has finite order. Hence
there exists $(x,y) \in \ZZ\big[\frac{1}{\alpha}\big] \oplus
\ZZ\big[\frac{1}{\beta}\big]$ such that $\phi(x,y) = (z,m)$ with $m \not=
0$. Since $\ZZ\big[\frac{1}{\alpha}\big] \oplus
\ZZ\big[\frac{1}{\beta}\big]$ is generated by elements of the form $(x,0)$
and $(0,y)$, we may in fact assume without loss of generality that there
is an element $x \in \ZZ\big[\frac{1}{\alpha}\big]$ such that $\phi(x,0) =
(z,m)$. Since $\alpha$ is infinite, there exist $n > m$ and $x' \in
\ZZ\big[\frac{1}{\alpha}\big]$ such that $n \cdot x' = x$, and this forces
$n \cdot\phi(x',0) = (z,m)$ which is impossible by our choice of $n$.

Since each $\Delta_2/H_n \cong L^*_{\alpha[1,n]} \times L^*_{\beta[1,n]}$,
the $K$-theory calculations for this example can also be verified using
the K\"unneth formula (Theorem~\ref{thm:direct limit} and
Proposition~\ref{prp:kunneth}).

\item Let $\phi$ be the homomorphism of $\ZZ^2$ determined by the
integer matrix $M := \left(\begin{smallmatrix} a & b \\ c & d
\end{smallmatrix}\right)$.
Suppose that $M$ is diagonalisable as a real $2 \times 2$ matrix, and that
its eigenvalues are greater than $1$ in modulus. Let $D := ad - bc$ be the
determinant of $M$.  For $n \ge 1$, let $H_n := M^n\ZZ^2$ and $\Lambda_n
:= \Delta_2/H_n$. Our assumption regarding the eigenvalues of $M$ ensures
that $\bigcap_{n=1}^\infty H_n = \{0\}$, so Theorem~\ref{thm:3-graph K-th}
and Corollary~\ref{cor:matrix realisation} imply that
$C^*(\tgrphlim(\Delta_2/H_n; p_n))$ is a simple A$\TT$ algebra of real
rank zero with
\[\textstyle
K_*(P_1 C^*(\tgrphlim(\Delta_2/H_n; p_n)) P_1) \cong \Big(
\ZZ\!\left[\frac{1}{D}\right] \oplus \ZZ,\ \varinjlim\left(\ZZ^2, \left(
  \begin{smallmatrix}
    \phantom{-}d & -b \\
    -c & \phantom{-}a
  \end{smallmatrix}
  \right)\right)\Big).
\]
In particular, let $M = \left(\begin{smallmatrix} a & -b \\ b &
\phantom{-}a \end{smallmatrix}\right)$ with $a^2 + b^2 > 1$. We may
identify $\ZZ^2$ with the group of Gaussian integers $\ZZ[i]$ by $(m,n)
\mapsto m + i n$, and then the group homomorphism of $\ZZ^2$ obtained from
multiplication by $M$ coincides with the group homomorphism of $\ZZ[i]$
obtained from multiplication by $a + i b$. Likewise $M^{\rm ca}$
implements multiplication by the conjugate $a - i b$. With $D := a^2 +
b^2$ and $\zeta := \frac{1}{a - ib} = \frac{a + ib}{a^2 + b^2}$, we have
\[\textstyle
  K_*(P_1 C^*(\tgrphlim(\Delta_2/H_n; p_n)) P_1) \cong \Big(
   \ZZ\!\big[\frac{1}{D}\big] \oplus \ZZ,\
   \ZZ\!\big[i, \frac{1}{\zeta}\big]\Big).
  \]
by Theorem~\ref{thm:3-graph K-th}.

\item More generally, a sequence of Gaussian integers $\zeta_j := a_j
+ b_j i$ with $|\zeta_j| > 1$ for all $j$ gives rise to a natural notion
of a Gaussian supernatural number $\zeta = \prod^\infty_{j=1} \zeta_j$.
Generalising the construction of the latter part of example~(2) above, let
$H_n := (\prod^n_{j=1} \overline{\zeta_j}) \ZZ[i]$ for each $n$, and
identify $\ZZ[i]$ with $\ZZ^2$ as a group to obtain a decreasing chain of
subgroups of $H_n$ of $\ZZ^2$ with trivial intersection.

Let $\alpha$ be the supernatural number $\alpha = \prod^\infty_{j=1}
|\zeta_j|^2$. Then
  \[\textstyle
  K_*(P_1 C^*(\tgrphlim(\Delta_2/H_n; p_n)) P_1) \cong \Big(
   \ZZ\!\big[\frac{1}{\alpha}\big] \oplus \ZZ,\
   \ZZ\big[i, \frac{1}{\zeta}\big]\Big).
  \]
by Theorem~\ref{thm:3-graph K-th} and Corollary~\ref{cor:matrix
realisation}.
\end{Enumerate}
\end{examples}

We now turn to the proof of Theorem~\ref{thm:3-graph K-th}; in particular,
we adopt the notation and conventions of Example~\ref{eg:HRBD}. Our first
step is to describe explicitly the $K$-theory of $C^*(\Delta_2/H_n)$ for a
fixed $n \in \NN\setminus\{0\}$. We do this using the results of
Section~\ref{sec:k>1K-theory}.

For $q \in \ZZ^k$ we write $q_+$ and $q_-$ for the positive and negative
parts of $q$. That is to say that $q_+$ and $q_-$ are the unique elements
of $\NN^k$ whose coordinate-wise minimum $q_+ \wedge q_-$ is equal to $0$,
and which satisfy $q = q_+ - q_-$.

For $q \in \ZZ^k$, a \emph{cycle of degree $q$} in a $k$-graph $\Lambda$
is a pair $(\mu,\nu)$ where $\mu \in \Lambda^{q_+}$ and $\nu \in
\Lambda^{q_-}$ such that $r(\mu) = r(\nu)$ and $s(\mu) = s(\nu)$. When $q
\in \NN^k$, $q = q_+$ and $q_- = 0$, so $\nu$ is a vertex, and $\mu$ is a
cycle in the usual sense: a path whose range and source coincide.

Let $H \subset \ZZ^2$ be a finite-index subgroup of $\ZZ^2$. Let $G =
\ZZ^2/H$. We view the ring $\ZZ G$ as the collection of functions $f : G
\to \ZZ$. For $X \subseteq G$ we denote the indicator function of the set
$X$ by $1_X$. We denote the point-mass at $g \in G$ by $\delta_g$.

Let $\Lambda := \Delta_2/H$. Let $E$ be the skeleton of $\Lambda$. That is
$E$ is the directed graph with the same vertices as $\Lambda$, and edges
$\Lambda^{e_1} \cup \Lambda^{e_2}$, with range and source inherited from
$\Lambda$. The degree map from $\Lambda$ restricts to a map from $E^1$ to
$\{e_1, e_2\}$. As in \cite{rsy, PRRS} we describe edges in $E$ as
\emph{blue} when they are of degree $e_1$ in $\Lambda$, and as \emph{red}
when they are of degree $e_2$. We often blur the distinction between
concatenation of edges in $E$ and the corresponding factorisation of a
path in $\Lambda$.

Recall that we are identifying $\Lambda^0$ with $G = \ZZ^2/H$. Hence,
given a path $\alpha = a_0 a_1\cdots a_n$ in $E$, we define functions
$f^b_\alpha$ and $f^r_\alpha$ in $\ZZ G$ by
\begin{eqnarray*}
f^b_\alpha (g) & =& \# \{ 0 \le j \le n :  r(a_j)  = g, d(a_j) =
e_1\} \\
f^r_\alpha (h) & =& \# \{ 0 \le k \le n :  r(a_k) = h, d(a_k) = e_2\}.
\end{eqnarray*}
The idea is that $f^b_\alpha(g)$ counts the number of blue edges in
$\alpha$ whose range is $g$, and $f^r_\alpha(g)$ does the same thing for
red edges.

We define $f_\alpha \in \ZZ G \oplus \ZZ G$ by $f_\alpha = f^b_\alpha
\oplus f^r_\alpha$. For a vertex $g \in \Lambda^0 = G$, we define $f^b_g$
and $f^r_g$ to be the zero element of $\ZZ G$, and $f_g = f^b_g \oplus
f^r_g$ is then the zero element of $\ZZ G \oplus \ZZ G$.

As $\Lambda = \Delta_2/H$, for each $g \in \Lambda^0 = G$ there is a
unique path $[g, g + (1,1)]$ of degree $(1,1)$ with range $g$. Using the
factorisation property, we can express this path as $b_g r_{g + [e_1]} =
r_g b_{g + [e_2]}$ (for $n \in \ZZ^2$, $[n]$ denotes the class of $n$ in
the quotient group $G = \ZZ^2/H$). We write $z_g$ for the function
$(\delta_{g + [e2]} - \delta_g) \oplus (\delta_g - \delta_{g + [e_1]})$ in
$\ZZ G \oplus \ZZ G$.

Given paths $\alpha  = a_0\cdots a_m$ and $\beta = b_0\cdots b_n$ in the
skeleton $E$ of $\Lambda$ such that $r(a_0) = r(b_0)$ and $s(a_m) =
s(b_n)$, let $f_{\alpha,\beta} := f_\alpha - f_\beta \in \ZZ G \oplus \ZZ
G$. Fix generators $h_1, h_2$ for $H$; so $[h_i] = [0]$ in $G$. By
definition of $\Lambda$, there are unique paths $\mu_1^+ \in
\Lambda^{(h_1)_+}$ and $\mu_1^- \in \Lambda^{(h_1)_-}$ with $r(\mu_1^\pm)
= 0$. Fix factorisations $\alpha_1^\pm$ of $\mu_1^\pm$ into edges from the
skeleton $E$. Since
\[
   s(\mu_1^+) = [(h_1)_+] = [(h_1)_-] = s(\mu_1^-)
\]
in $G$, the pair $(\mu_1^+, \mu_1^-)$ is a cycle of degree
$h_1$ in $\Lambda$ with range $[0]$. The same construction for
$h_2$ gives a cycle $(\mu_2^+, \mu_2^-)$ of degree $h_2$ with
range $[0]$ and fixed factorisations $\alpha_2^\pm$ of
$\mu_2^\pm$ into edges from the skeleton $E$.

\begin{lem}\label{lem:K-th of quotient}
With the notation established in the preceding paragraphs, the chain
complex~\eqref{eq:complex} can be described as follows:
\begin{Enumerate}
\item for each $g \in G$, $\partial_1 (\delta_g \oplus 0) =
    \delta_g - \delta_{g +
    [e_1]}$, $\partial_1 (0 \oplus \delta_g) = \delta_g -
    \delta_{g + [e_2]}$, and
  \[
  \partial_2(\delta_g) = (\delta_{g + [e_2]} - \delta_g) \oplus
  (\delta_g - \delta_{g + [e_1]}) = z_g.
  \]
\item
$\coker(\partial_1) \cong \ZZ$ is generated by $\delta_0 +
\Im(\partial_1)$;
\item
$\ker(\partial_2) \cong \ZZ$ is generated by $1_G$;
\item
For each $h \in G$, the set $\{z_g : g\in G\setminus\{h\}\}$ is a basis
for $\Im(\partial_2) \cong \ZZ^{|G|-1}$.
\item
Fix any two factorisations $\alpha$ and $\beta$ of a path $\mu$ in
$\Lambda$ into edges from $E$. Then $f_\alpha - f_\beta \in
\Im(\partial_2)$, and $\partial_1(f_\alpha) =
\partial_1(f_\beta) = \delta_{r(\alpha)} - \delta_{s(\alpha)}$.
\item
$\ker(\partial_1)$ is the subgroup of $\ZZ G \oplus \ZZ G$ generated by
the elements $f_{\alpha,\beta}$ where $\alpha$ and $\beta$ are paths in
the skeleton $E$ with $r(\alpha) = r(\beta)$ and $s(\alpha) = s(\beta)$.
\item There is an isomorphism $\psi$ of $H$ onto
    $\ker(\partial_1)/\Im(\partial_2)$ which takes $d(\mu) -
    d(\nu)$ to $f_{\alpha,\beta} + \Im(\partial_2)$ for each
    cycle $(\mu,\nu)$ in $\Lambda$ and pair of factorisations
    $\alpha$ of $\mu$ and $\beta$ of $\nu$. In particular, for any
    basis $B$ for $\Im(\partial_2)$, the set $B \cup
    \{f_{\alpha_1^+, \alpha_1^-}, f_{\alpha_2^+, \alpha_2^-}\}$
    is a basis for $\ker(\partial_1) \cong \ZZ^{|G|+1}$ (where
    $\alpha_i^{\pm}$ are the fixed factorisations of the paths
    $\mu_i^{\pm}$ of degree $(h_i)_{\pm}$ described above).
\end{Enumerate}
In particular, $K_*(C^*(\Lambda)) \cong (\ZZ^2, H)$ where the class of the
identity in $K_0$ is identified with the element $(|G|,0)$ of $\ZZ^2$.
\end{lem}
\begin{proof}
(1) The adjacency matrix $M_1$ associated to $(\Lambda^0, \Lambda^{e_1},
r, s)$ is the permutation matrix determined by translation by $[e_1]$ in
$G$ and similarly for $M_2$. The first statement then follows from the
formulae for $\partial_1$ and $\partial_2$ in terms of $M_1$ and $M_2$.

(2) The formulae for $\partial_1(\delta_g \oplus 0)$ and $\partial_1(0
\oplus \delta_g)$ show that $\delta_g + \Im(\partial_1) = \delta_{g +
[e_i]} + \Im(\partial_1)$ in $\coker(\partial_1)$ for $i = 1,2$ and $g \in
G$. Since the action of $\ZZ^2$ on $G$ by translation is transitive, this
establishes~(2).

(3) Using the formula for $\partial_2$ established in~(1), one can see
that for $f \in \ZZ G$, $\partial_2(f) = f_1 \oplus f_2$ where
\[
f_1(g) = -f(g) + f(g - [e_1])\qquad\text{and}\qquad
f_2(g) = f(g) - f(g - [e_2])).
\]
Hence $f \in \ker(\partial_2)$ if and only if $f(g) = f(g - [e_1]) = f(g -
[e_2])$ for all $g \in G$, and since the action of $\ZZ^2$ on $G$ is
transitive, this establishes~(3).

(4) Part~(1) establishes that $\Im(\partial_2)$ is generated by $\{z_g : g
\in G\}$. A simple calculation shows that $\sum_{g \in G} z_g = 0$ in $\ZZ
G \oplus \ZZ G$, and it follows that for any $h \in G$, the set $\{z_g :
g\in G\setminus\{h\}\}$ generates $\Im(\partial_2) \cong \ZZ^{|G|-1}$.
Since $\ker(\partial_2)$ has rank 1, the rank of its image is $|G| - 1$,
establishing~(4).

(5) By part~(4), the image of $\partial_2$ is generated by elements of the
form $f_\alpha - f_\beta$ where $\alpha$ and $\beta$ are the two possible
factorisations of a path in $\Lambda^{(1,1)}$. Since $f_{\alpha\beta} =
f_\alpha + f_\beta$ when $\alpha$ and $\beta$ are paths in $E$ which can
be concatenated, this establishes the first claim. The second statement
follows from a straightforward calculation using that
\begin{equation}\label{eq:p1 formula}
\partial_1(f^b \oplus f^r)(g) = f^b(g) - f^b(g - [e_1]) + f^r(g)
- f^r(g - [e_2]).
\end{equation}

(6) If $\alpha,\beta$ are paths in the skeleton with $r(\alpha) =
r(\beta)$ and $s(\alpha) = s(\beta)$ then $f_{\alpha,\beta}$ belongs to
$\ker(\partial_1)$ by ~(5).

We must show that every $f \in \ker(\partial_1)$ can be written as a
$\ZZ$-linear combination of elements of the form $f_{\alpha,\beta}$. First
note that it suffices to treat the case where $f$ takes only nonnegative
values (this is because $1_G \oplus 1_G$ can be so expressed). So suppose
that $f$ takes nonnegative values, and write $f = f^b \oplus f^r$. Let
$E_f$ be the directed graph with vertices $G$ and which contains $f^b(g)$
parallel copies of the blue edge in $E$ with range $g$ and $f^r(g)$ copies
of the red edge in $E$ with range $g$. If $E_f$ contains a \emph{terminal
vertex} $g$ which receives at least one edge but emits no edges at all,
then $f^b(g) + f^r(g) \not= 0$, but $f^b(g - [e_1]) = f^r(g - [e_2]) = 0$,
and~\eqref{eq:p1 formula} shows that $\partial_1(f)(g) \not= 0$. Hence
$E_f$ contains no such vertex, and therefore must either contain a cycle
$\alpha$ or contain no edges at all. In the latter case, the claim is
trivial, and in the former case, $f \ge f_\alpha$, and removing the cycle
$\alpha$ from $E_f$ produces the graph $E_{f - f_\alpha}$ for the function
$f - f_\alpha$. After finitely many such steps, we must obtain a forest
with no terminal vertex. The only such forest is the empty graph which
corresponds to the function $0 \oplus 0$. That is $f - \sum_{\alpha \in L}
f_\alpha = 0 \oplus 0$ for some collection $L$ of cycles, and this
proves~(6).

(7) Suppose that $(\mu,\nu)$ is a cycle in $\Lambda$. Then
\[
s(\mu) - [d(\mu)] = r(\mu) = r(\nu) = s(\nu) - [d(\nu)] = s(\mu)
-
[d(\nu)]
\]
in $G = \Lambda^0 = \ZZ^2/H$, so $d(\mu) - d(\nu) \in H$. It is
clear from the definition of $\Lambda$ that each element of $H$
arises as $d(\mu) - d(\nu)$ for some cycle $(\mu,\nu)$ in
$\Lambda$.

To see that the assignment $d(\mu) - d(\nu) \mapsto
f_{\alpha,\beta} + \Im(\partial_2)$ is well defined, we must
show two things. First that for distinct factorisations
$\alpha$ and $\alpha'$ of $\mu$ and distinct factorisations
$\beta$ and $\beta'$ of $\nu$, the difference $f_{\alpha,\beta}
- f_{\alpha',\beta'}$ lies in the image of $\partial_2$. This
follows from~(5). Second, we must show that if $(\mu,\nu)$ and
$(\mu',\nu')$ are cycles in $\Lambda$ with $d(\mu) - d(\nu) =
d(\mu') - d(\nu')$, then there exist factorisations $\alpha$ of
$\mu$, $\beta$ of $\nu$, $\alpha'$ of $\mu'$, and $\beta'$ of
$\nu'$ such that $f_{\alpha,\beta} - f_{\alpha',\beta'}$ is in
$\Im(\partial_2)$. To see this, first note that by factorising
$\mu = \mu'\tau$ and $\nu = \nu'\tau$ where $d(\tau) = d(\mu)
\wedge d(\nu)$, we can reduce to the case where $d(\mu) \wedge
d(\nu) = 0$. Next we claim that it suffices to consider the
case where $r(\mu) = r(\nu) = r(\mu') = r(\nu') = [0]$. To see
this, fix $\eta$ in $[0] \Lambda r(\mu)$ and note that the
cycle $(\eta\mu, \eta\nu)$ corresponds to the same class as
$(\mu,\nu)$ in $\ker(\partial_1)/\Im(\partial_2)$. Factorise
$\eta\mu = \xi\rho$ and $\eta\nu = \omega\sigma$ where $d(\xi)
= d(\mu)$, $d(\omega) = d(\nu)$ and $d(\rho) = d(\sigma) =
d(\eta)$. Since each $g\Lambda^n$ is a singleton and since
$\ZZ^2$ acts on $\Lambda$ by translation, $(\xi,\omega)$ is a
cycle with range $[0]$, and $\rho = \sigma$. Hence the cycle
$(\xi,\omega)$ corresponds to the same class in
$\ker(\partial_1)/\Im(\partial_2)$ as $(\mu,\nu)$. After
shifting $(\mu', \nu')$ in a similar way we may assume that
both cycles have range $[0]$. We now have cycles $(\mu,\nu)$
and $(\mu',\nu')$ with range $[0]$ and such that $d(\mu) -
d(\nu) = d(\mu') - d(\nu')$ and $d(\mu) \wedge d(\nu) = 0 =
d(\mu') \wedge d(\nu')$. Since $[0]\Lambda^n$ is a singleton
for any $n \in \ZZ^2$, this forces $\mu = \mu'$ and $\nu =
\nu'$. This completes the proof that $d(\mu)-d(\nu) \mapsto
f_{\alpha,\beta} + \Im(\partial_2)$ is well defined.

That $f_{\alpha\beta} = f_{\alpha} + f_{\beta}$ ensures that
$\psi(g + h) = \psi(g) + \psi(h)$, and that $f_{\beta,\alpha} =
-f_{\alpha,\beta}$ shows that $\psi(-g) = -\psi(g)$. Hence
$\psi$ is a homomorphism. By part~(6), to see that $\psi$ is
surjective, we just need to show that each $f_{\alpha,\beta} +
\Im(\partial_2)$ is in the range of $\psi$. This is clear
because $f_{\alpha,\beta} + \Im(\partial_2)$ is precisely
$\psi(d(\mu) - d(\nu))$ where $\mu$ factorises as $\alpha$ and
$\nu$ factorises as $\beta$. Finally, to see that $\psi$ is
injective, note that if $f_{\alpha,\beta} \in \Im(\partial_2)$,
then $d(\mu) = d(\nu)$ where $\mu$ factorises as $\alpha$ and
$\nu$ factorises as $\beta$. This completes the proof that
$\psi : H \to \ker(\partial_1)/\Im(\partial_2)$ is an
isomorphism. The remaining statement follows from~(4) and that
$(\mu_1^+, \mu_1^-)$ and $(\mu_2^+, \mu_2^-)$ are cycles whose
degrees form a basis for $H$. This proves~(7).

The final statement of the Lemma follows from~\eqref{eq:K-th maps}.
\end{proof}

We now consider two consecutive graphs in the sequence of
covering systems described in Example~\ref{eg:HRBD}, and
describe the homomorphism of $K$-invariants obtained from
Proposition~\ref{prp:embedding}(6).

\begin{thm}\label{thm:two quotients}
Consider the situation described in Example~\ref{eg:HRBD}, and fix $n \in
\NN\setminus\{0\}$. For $i = n, n+1$, let $\Lambda_i := \Delta_2/H_i$, and
consider the commuting diagram
\[
\begin{CD}
  0 @<<< \ZZ \Lambda_n^0 @<{\partial^{\Lambda_n}_1}<<
\ZZ \Lambda_n^0 \oplus \ZZ \Lambda_n^0 @<{\partial^{\Lambda_n}_2}<<
\ZZ \Lambda_n^0 @<<< 0\\
@.  @VVp_n^*V  @VV p_n^* \oplus p_n^*V  @VVp_n^*V @. \\
  0 @<<< \ZZ \Lambda_{n+1}^0 @<{\partial^{\Lambda_{n+1}}_1}<<
\ZZ \Lambda_{n+1}^0 \oplus \ZZ \Lambda_{n+1}^0
@<{\partial^{\Lambda_{n+1}}_2}<<
\ZZ \Lambda_{n+1}^0 @<<< 0\\
\end{CD}
\]
\begin{Enumerate}
\item The right-hand vertical map $p_n^* : \ZZ\Lambda^0_n \to
    \ZZ\Lambda^0_{n+1}$ restricts to a homomorphism
    $p_n^*|_{\ker(\partial_2^{\Lambda_n})} :
    \ker(\partial_2^{\Lambda_n}) \to
    \ker(\partial_2^{\Lambda_{n+1}})$ characterised by
    $p_n^*|_{\ker(\partial_2^{\Lambda_n})}(1_{G_n}) =
    1_{G_{n+1}}$.
\item
The left-hand vertical map $p_n^* : \ZZ\Lambda^0_n \to \ZZ\Lambda^0_{n+1}$
induces a homomorphism $\widetilde{p_n^*} : \coker(\partial^{\Lambda_n}_1)
\to \coker(\partial^{\Lambda_{n+1}}_1)$ characterised by
   \[
   \widetilde{p_n^*}(\delta_0 + \Im(\partial^{\Lambda_n}_1))
   = [H_n : H_{n+1}] \cdot \delta_0 + \Im(\partial^{\Lambda_{n+1}}_1).
   \]
\item
The middle vertical map $p_n^* \oplus p_n^* : \ZZ\Lambda^0_n \oplus
\ZZ\Lambda^0_n \to \ZZ\Lambda^0_{n+1} \oplus \ZZ\Lambda^0_{n+1}$ induces a
homomorphism $(p_n^* \oplus p_n^*)^{\sim} :
\ker(\partial^{\Lambda_n}_1)/\Im(\partial^{\Lambda_n}_2) \to
\ker(\partial^{\Lambda_{n+1}}_1)/\Im(\partial^{\Lambda_{n+1}}_2)$ such
that the following diagram commutes.
  \[
  \begin{CD}
    {H_n} @>{\psi_n}>>
    {\ker(\partial_1^{\Lambda_n})/\Im(\partial_2^{\Lambda_n})} \\
    @VV{m_{H_n,H_{n+1}}}V @VV{(p_n^* \oplus p_n^*)^{\sim}}V \\
    {H_{n+1}} @>{\psi_{n+1}}>>
    {\ker(\partial_1^{\Lambda_{n+1}})/\Im(\partial_2^{\Lambda_{n+1}})} \\
  \end{CD}
  \]
  where $\psi_n$ and $\psi_{n+1}$ are the isomorphisms obtained from
  Lemma~\ref{lem:K-th of quotient}(7), and $m_{H_n, H_{n+1}}$ is as in
  Theorem~\ref{thm:3-graph K-th}(2).
\end{Enumerate}
Under the isomorphism
\[
K_*(C^*(\Lambda_i)) \cong \big(\coker(\partial^{\Lambda_i}_1) \oplus
\ker(\partial_2^{\Lambda_i}),\:
\ker(\partial^{\Lambda_i}_1)/\Im(\partial^{\Lambda_i}_2)\big)
\]
obtained from Corollary~\ref{cor:k=2K-theory}, the maps
described in~(1), (2)~and~(3) determine the map
$(\iota_{p_n})_* : K_*(C^*(\Lambda_n)) \to
K_*(C^*(\Lambda_{n+1}))$ obtained from
Proposition~\ref{prp:embedding}(6).
\end{thm}
\begin{proof}
Lemma~\ref{lem:K-th of quotient}(3) ensures that $1_{G_i}$
generates $\ker(\partial_2^{\Lambda_i})$ for $i = n, n+1$. The
formula for $p_n^*$ shows that $p_n^*(1_{G_n}) = 1_{G_{n+1}}$,
which gives~(1). Statement~(2) follows from the formula for
$p_n^*$ combined with the observation that for $i = n, n+1$,
the $\delta_g$, $g \in G_i$ are all equivalent modulo
$\Im(\partial_1^{\Lambda_i})$.

It remains only to prove~(3). We first consider the case where
$H_n = \ZZ^2$, so $G_n = \{0\}$ and $\Lambda_n$ is a copy of
the $2$-graph $T_2 \cong \NN^2$ (as a category) with one vertex
and one morphism $\lambda_m$ of each degree $m \in \NN^2$. In
this case, $\psi_n$ is just the identity map from $\ZZ^2$ to
$\ZZ \oplus \ZZ$. Let $h_1, h_2$ be a pair of generators for
$H_2$. It suffices to show that $(p_n^* \oplus p_n^*)(h_i) =
[\ZZ^2 : H_{n+1}] \cdot h_i \in H_{n+1}$ for $i = 1,2$. We just
argue that this happens for $i = 1$ (the case $i=2$ follows
from a symmetric argument). Writing $h_1 = (x,y)$ where $x,y
\in \ZZ$, the formula for $p_n^*$ then ensures that $(p_n^*
\oplus p_n^*)$ takes $h_1$ to $x 1_{G_2} \oplus y 1_{G_2}$. To
see that this is $[\ZZ^2 : H_{n+1}] \cdot h_1$, let $f :=
f_{\alpha_1^+,\alpha_1^-} = \psi(h_1)$ be the function in $\ZZ
G_{n+1} \oplus \ZZ G_{n+1}$ obtained from Lemma~\ref{lem:K-th
of quotient}(7). By definition of $f$, $f = f_b + f_r$ where
the entries of $f_b$ sum to $x$ and the entries of $f_r$ sum to
$y$. For $g \in G_2$, let $g \cdot f_b$ be the function
determined by $g\cdot f_b(h) = f_b(h - g)$, and similarly for
$f_r$. Since $G_{n+1}$ acts freely and transitively on
$\Lambda_{n+1}^0 = G_{n+1}$, it follows that
\begin{equation}\label{eq:g-shift sum}
\sum_{g \in G_{n+1}} g\cdot f = x 1_{G_{n+1}} \oplus y 1_{G_{n+1}}.
\end{equation}
The proof of statement~(7) in Lemma~\ref{lem:K-th of quotient}
shows that each $g\cdot f := g\cdot f_b \oplus g\cdot f_r$
represents the same class as $f$ in
$\ker(\partial_1^{\Lambda_{n+1}})/\Im(\partial_2^{\Lambda_{n+1}})$.
Hence the left-hand side of~\eqref{eq:g-shift sum} has the same
class in
$\ker(\partial_1^{\Lambda_{n+1}})/\Im(\partial_2^{\Lambda_{n+1}})$
as $|G_{n+1}| \cdot h_1$ as required.

For the general case, let $p_{[1,n]} := p_1 \circ \dots \circ
p_{n-1}$ and $p_{[1,n+1]} := p_1 \circ \dots \circ p_n$ be the
coverings of $T_2$ by $\Lambda_n$ and $\Lambda_{n+1}$ obtained
by composing the first $n$ and $n+1$ levels of the covering
system; we may apply the argument of the previous paragraph to
these coverings. Then $p_{[1,n+1]} = p_{[1,n]} \circ p_n$, so
$p_{[1,n+1]}^* \oplus p_{[1,n+1]}^* = (p_{[1,n]}^* \oplus
p_{[1,n]}^*)\circ(p_n^* \oplus p_n^*)$, and since these maps
induce homomorphisms between
$\ker(\partial_1^{T_2})/\Im(\partial_2^{T_2})$ and
$\ker(\partial_1^{\Lambda_{n+1}})/\Im(\partial_2^{\Lambda_{n+1}})$
which are rational isomorphisms, it follows that $(p_n^* \oplus
p_n^*)$ behaves as claimed.

The final statement follows from Corollary~\ref{cor:k=2K-theory}.
\end{proof}

We are now ready to prove Theorem~\ref{thm:3-graph K-th}.

\begin{proof}[Proof of Theorem~\ref{thm:3-graph K-th}]
Proposition~\ref{prp:embedding} shows that $P_1$ is full so that
compression by $P_1$ induces an isomorphism on $K$-theory. The formulae
for the $K$-groups in statements (1)~and~(2) follow from
Lemma~\ref{lem:K-th of quotient} and Theorem~\ref{thm:two quotients} and
the continuity of the $K$-functor.

Since each $v (\Delta_2/H_n) w \not= \emptyset$ for all $n \in
\NN\setminus\{0\}$, and $v,w \in \Delta^0_2/H_n$, the $3$-graph
$\tgrphlim(\Delta_2/H_n, p_n)$ is cofinal. Moreover a given
infinite path $x$ in $\tgrphlim(\Delta_2/H_n, p_n)$ is periodic
with period $m \in \ZZ^2$ if and only if every infinite path in
$\tgrphlim(\Delta_2/H_n, p_n)$ is periodic with period $m$,
which in turn is equivalent to the condition that $m \in
\bigcap_{n=1}^\infty H_n$. It follows from Lemma~\ref{lem:LPF
analogue} that $\tgrphlim(\Delta_2/H_n, p_n)$ is simple if and
only if $\bigcap H_n = \{0\}$; moreover, in this case, the
argument of the second part of \cite[Section~5]{PRRS} shows
that $C^*(\tgrphlim(\Delta_2/H_n, p_n))$ has unique trace.

We next claim that each $C^*(\Delta_2/H_n) \cong
M_{[\ZZ^2:H_n]}(C(\TT^2))$. To see this, one checks that $h \mapsto
s_{[(0,h_+)]} s^*_{[(0,h_-)]}$ is a group isomorphism $H_n \to
\mathcal{U}(s_{[0]} C^*(\Delta_2/H_n) s_{[0]})$ for each $n$. The standard
argument used in \cite[Lemma~3.9]{PRRS} shows that each $s_{[(0,h_+)]}
s^*_{[(0,h_-)]}$ has full spectrum. One can then deduce that $s_{[0]}
C^*(\Delta_2/H_n) s_{[0]} \cong C^*(H_n) \cong C^*(\ZZ^2) \cong C(\TT^2)$.
For $m \in \ZZ^2/H_n$, define $V_m := s^*_{[0,m]} \in C^*(\Delta_2/H_n)$.
Applying Lemma~\ref{lem:matrix alg} to these partial isometries with $p =
s_{[0]}$ and $q = 1_{C^*(\Delta_2/H_n)}$proves that $C^*(\Delta_2/H_n)
\cong M_{[\ZZ^2 : H_n]}(C(\TT^2))$.

It now follows from \cite[Theorem~1.3]{BBEK} that
$C^*(\tgrphlim(\Delta_2/H_n; p_n))$ has real-rank 0. The
classification of such algebras of
D\u{a}d\u{a}rlat-Elliott-Gong (see \cite[Theorem~3.3.1]{Ror}),
and the $K$-theory calculations above complete the proof.
\end{proof}

\begin{rmk}\label{rmk:HRBD}
\emph{Higher-rank Bunce-Deddens algebras and generalised
odometer actions.} We consider a slightly more general version
of the situation described in Example~\ref{eg:HRBD}. Let $H_1
\supset H_2 \supset H_3 \supset \dots$ be a chain of
finite-index subgroups of $\ZZ^k$ such that $\bigcap_n H_n = \{
0 \}$. For each $n$, let $p_n : \Delta_k/H_{n+1} \to
\Delta_k/H_n$ be the canonical covering induced by the quotient
maps described above, and let $\cocycle_n : \Delta_k/H_{n+1}
\to S_1$ be the trivial cocycle. This data specifies a sequence
$(\Delta_k/H_n, \Delta_k/H_{n+1}, p_n)^\infty_{n=1}$ of
row-finite covering systems of $k$-graphs with no sources.
Applying Corollary~\ref{cor:tower graph}, we obtain a
$(k+1)$-graph $\tgrphlim(\Delta_k/H_n; p_n)$.

We claim that $P_1 C^*(\tgrphlim(\Delta_k/H_n; p_n)) P_1$ can be thought
of as a higher-rank Bunce-Deddens algebra. We justify this with a
description of $P_1C^*(\tgrphlim(\Delta_k/H_n; p_n)) P_1$ as a crossed
product by a generalised odometer action. We assume here that $H_1 =
\ZZ^k$ so that $\Delta_k/H_1$ is a copy of the $k$-graph $T_k \cong \NN^k$
(as a category) with one vertex and one morphism $\lambda_m$ of each
degree $m \in \NN^k$.

One way to realise the Bunce-Deddens algebras is as crossed
products of algebras of continuous functions on Cantor sets by
generalised odometer actions. Given a supernatural number
$\alpha = \alpha_1\alpha_2 \cdots$, let $G_n :=
\ZZ/\alpha[1,n]\ZZ$ for all $n$. Then for each $n$, since
$\alpha[1,n+1]\ZZ \supset \alpha[1,n]\ZZ$, there is a natural
surjective group homomorphism from $G_{n+1}$ to $G_n$. Hence,
we may form the projective limit group $\varprojlim(G_n, p_n)$.
The automorphism $\tau(g_1, g_2, \dots) = (g_1 + [1], g_2 +
[1], \dots)$ for $(g_1, g_2, \dots) \in \varprojlim(G_n, p_n)$
can then naturally be regarded as an odometer action on
$\varprojlim(G_n, p_n)$. The Bunce-Deddens algebra of type
$\alpha$ is the crossed product $C(\varprojlim(G_n, p_n))
\rtimes_{\tilde{\tau}} \ZZ$ where $\tilde{\tau}$ is the
automorphism of $C(\varprojlim(G_n, p_n))$ induced by $\tau$
(see \cite[Examples~1(3)]{Riedel1982}).

There is an analogous realisation of
$P_1C^*(\tgrphlim(\Delta_k/H_n, p_n)) P_1$ as follows. Let
$\Lambda := \tgrphlim(\Delta_k/H_n, p_n)$. Let $F$ denote the
fixed-point algebra of $C^*(\Lambda)$ for the gauge action
$\gamma$ of $\TT^{k+1}$.  Note that by
Remark~\ref{rmk:restricted-gauge}, the restriction of the gauge
action to $P_1C^*(\Lambda) P_1$ is trivial on the last
coordinate of $\TT^{k+1}$ and therefore becomes an action by
$\TT^{k}$ denoted $\tilde{\gamma}$. Recall that
$\Lambda^\infty$ denotes the collection of infinite paths in
$\Lambda$ (see Notation~\ref{ntn:inf paths}). It is not hard to
see that $P_1 F P_1$ is canonically isomorphic to
$C(v\Lambda^\infty)$ where $v$ is the unique vertex of
$\Delta_k/H_1 \cong T_k$. Let $G_n := \ZZ^k/H_n$ for each $n$,
and let $p_n : G_{n+1} \to G_n$ be the induced map $p_n(m +
H_{n+1}):= m + H_n$. Observe that $G = \varprojlim(G_n, p_n)$
is a compact abelian group. By functoriality of the projective
limit the quotient maps $\ZZ^k \to \ZZ^k/H_n$ induce a
homomorphism $j: \ZZ^k \to G$; injectivity of $j$ follows from
the fact $\bigcap_n H_n = \{ 0 \}$. There is an action $\tau$
of $\ZZ^k$ on $G$ given by $\tau_m(g_1, g_2, \dots) = (g_1 +
[m], g_2 + [m], \dots)$, which generalises the odometer action
discussed above. Since there is just one infinite path in
$T_k$, the arguments of Section~\ref{sec:simplicity} show that
$v\Lambda^\infty \cong G$ as a topological space. Note that for
every $m \in \NN^k$, the generator $s_{\lambda_m}$ associated
to the unique path $\lambda_m \in T_k^m$ is a unitary in
$P_1C^*(\Lambda) P_1$ and that under the identification of $P_1
F P_1$ with $C(v\Lambda^\infty) = C(G)$ conjugation by
$s_{\lambda_m}$ implements the automorphism induced by the
homeomorphism $\tau_m$ of $G$. It follows that the reduction of
the path groupoid (see \cite[Section~2]{kp}) of $\Lambda$ to
$v\Lambda^\infty$ is isomorphic to the semidirect product
groupoid $G \rtimes_{\tau} \ZZ^k$. Therefore, standard
arguments show that
\[
P_1C^*(\Lambda) P_1 \cong  C(G) \rtimes_{\tilde{\tau}} \ZZ^k
\]
where $\tilde{\tau}$ is the action induced by $\tau$. Note that under this
identification the restricted gauge action $\tilde{\gamma}$ coincides with
the dual action of $\TT^{k} = \widehat{\ZZ^k}$.

The action of $G$ on $C(G)$ induced by translation in $G$
yields an action of $G$ on $C(G) \rtimes_{\tilde{\tau}} \ZZ^k$
which commutes with the dual action of $\TT^{k} =
\widehat{\ZZ^k}$. Thus we obtain an action $\alpha$ by the
compact abelian group $G \times \TT^{k}$ with fixed point
algebra isomorphic to $\CC$. Hence, $C(G)
\rtimes_{\tilde{\tau}} \ZZ^k$ (and thus $P_1C^*(\Lambda) P_1$)
admits an ergodic action of a compact abelian group.  Such
ergodic actions have been classified in \cite[4.5, 6.1]{opt};
the invariant is a symplectic bicharacter $\chi_\alpha$ on
$\widehat{G} \times \ZZ^{k}$, the dual of $G \times \TT^{k}$.
This gives rise to an alternative description of the
$C^*$-algebra as a twisted group $C^*$-algebra with the group
$\widehat{G} \times \ZZ^{k}$ and a $2$-cocycle associated to
the bicharacter $\chi_\alpha$ (only its cohomology class is
determined by the bicharacter). It follows that
$$C(G) \rtimes_{\tilde{\tau}} \ZZ^k
\cong C(\TT^k)\rtimes \widehat{G}$$ where the action of $\widehat{G}$ on
$C(\TT^k)$ arises by translation from the embedding $\widehat{G} \to
\TT^k$ dual to $j: \ZZ^k \to G$.
\end{rmk}

\end{document}